\documentclass[review]{elsarticle}

\usepackage{lineno,hyperref}
\modulolinenumbers[5]

\makeatletter
\def\ps@pprintTitle{%
 \let\@oddhead\@empty
 \let\@evenhead\@empty
 \def\@oddfoot{}%
 \let\@evenfoot\@oddfoot}
\makeatother

\usepackage{amsmath, amssymb, amsfonts}

\usepackage[ruled,vlined]{algorithm2e}

\begin{document}
\begin{frontmatter}
\title{A diffusion generated method for \\ computing Dirichlet partitions}
\tnotetext[mytitlenote]{B. Osting is partially supported by NSF DMS 16-19755.}

\author{Dong Wang}
\ead{dwang@math.utah.edu}
\author{Braxton Osting\corref{mycorrespondingauthor}}
\cortext[mycorrespondingauthor]{Corresponding author}
\ead{osting@math.utah.edu}
\address{Department of Mathematics, University of Utah, Salt Lake City, UT}

\begin{abstract} 
A Dirichlet $k$-partition of a closed $d$-dimensional surface is a collection of $k$ pairwise disjoint open subsets such that the sum of their first Laplace-Beltrami-Dirichlet eigenvalues is minimal. In this paper, we develop a simple and efficient diffusion generated method to compute Dirichlet $k$-partitions for $d$-dimensional flat tori and spheres. 
For the $2d$ flat torus, for most values of $k=3$--9,11,12,15,16, and 20, we obtain hexagonal honeycombs. 
For the $3d$ flat torus and $k=2,4,8,16$, we obtain the rhombic dodecahedral honeycomb, the Weaire-Phelan honeycomb, and Kelvin's tessellation by truncated octahedra. 
For the $4d$ flat torus, for $k=4$, we obtain a constant extension of the rhombic dodecahedral honeycomb along the fourth direction and for $k=8$, we obtain a 24-cell honeycomb. 
For the $2d$ sphere, we also compute Dirichlet partitions for $k=3$--7,9,10,12,14,20.
Our computational results agree with previous studies when a comparison is available. As far as we are aware, these are the first published results for Dirichlet partitions of the $4d$ flat torus. 
\end{abstract}

\begin{keyword}
Dirichlet partition \sep diffusion generated method \sep honeycomb  \sep Weaire-Phelan  structure \sep Kelvin structure \sep 24-cell
\MSC[2010] 
49Q10 \sep 
35R01 \sep 
05B45 
\end{keyword}

\end{frontmatter}


\clearpage

\section{Introduction} \label{s:intro}
For $d\geq 2$, let $U$ be either an open bounded domain in $\mathbb R^d$ with Lipschitz boundary or a closed, smooth, $d$-dimensional manifold. For $k\geq 2$ fixed, the \emph{Dirichlet $k$-partition problem} for
$U$ is to choose a $k$-partition, {\it i.e.}, $k$ disjoint quasi-open sets
$U_1, U_2, \ldots, U_k \subseteq U$, that attains
\begin{equation} \label{eq:ContDirPart}
\min_{U = \cup_{\ell \in [k]} U_\ell} \ \sum_{\ell \in [k]} \lambda_1(U_\ell)
\end{equation} 
where 
\begin{equation}
  \lambda_1(U) := \min_{\substack{u \in H^1_0(U)\\ \|u\|_{L^2(U)}=1}} E(u) 
  \qquad \textrm{and} \qquad 
  E(u) := \begin{cases} 
    \int_U |\nabla u|^2 \  dx & u \in H^1_0(U)\\ 
    \infty & \text{ otherwise}
\end{cases}. 
\end{equation}
Here, $E$ is the Dirichlet energy and $\lambda_1(U)$ is the first
Dirichlet eigenvalue of the  Laplace-Beltrami operator, $-\Delta$, on $U$ with Dirichlet boundary conditions imposed on $\partial U$.  
We refer to any  $k$-partition that attains the minimum in \eqref{eq:ContDirPart} as a \emph{Dirichlet $k$-partition of $U$}, or simply a \emph{Dirichlet partition} when $k$ and $U$ are understood. 
Observe that by the monotonicity of Dirichlet eigenvalues, any Dirichlet partition satisfies
$\overline{U} = \cup_{i=1}^k \overline{U_i}$, which justifies the use of the word ``partition'' in the name. 
The existence of optimal partitions in the class
of quasi-open sets was proved in \cite{bucur1998existence} and, subsequently,
several papers have investigated properties of
optimal partitions including the regularity of the partition interfaces and the asymptotic behavior  as
$k \to \infty$ \cite{caffarelli2007optimal,Helffer2010b,bourdin2010optimal}. Dirichlet partitions arise in the study of Bose-Einstein condensates \cite{bao2004ground,bao2004computing,chang2004segregated}, 
models for interacting agents \cite{conti2002nehari,conti2003optimal,chang2004segregated,cybulski2005,cybulski2008}, 
and have recently been studied in the context of data analysis \cite{osting2013minimal,Zosso2015,Osting_2017}.

\subsection{Results}
In this paper, \emph{we develop an efficient diffusion generated method for computing Dirichlet partitions of $d$-dimensional flat tori and spheres}; see Algorithm~\ref{a:MBO}. 
The method is best motivated by a mapping formulation of Dirichlet partitions that we review in Section \ref{s:num}. 
The method is very simple, consisting  of iterating the following three steps: 
(i) Evolve $k$ functions on $U$ by the diffusion equation until time $\tau$. 
(ii) At each point of $U$, find which of the $k$ functions is largest and set  the other functions to zero. 
(iii) Renormalize each of the $k$ functions. 
 This method is implemented using the Fast Fourier Transform (FFT) and Spherical Harmonic Transform (SHT), as described in Section \ref{s:comp}. 

in Section \ref{s:numeircs}, we present results of extensive numerical experiments. 
For the $2d$ flat torus, for most values of $k=3$--9,11,12,15,16, and 20, we obtain hexagonal honeycombs. 
For the $3d$ flat torus and $k=2,4,8,16$, we obtain the rhombic dodecahedral honeycomb, the Weaire-Phelan honeycomb, and Kelvin's tessellation by truncated octahedra. 
For the $4d$ flat torus and $k=4$, we obtain a constant extension of the rhombic dodecahedral honeycomb along the fourth direction and for $k=8$, we obtain a 24-cell honeycomb. 
For the $2d$ sphere, we also compute Dirichlet partitions for $k=3$--7,9,10,12,14,20.
Our results agree with previous studies when a comparison is available. As far as we are aware, these are the first published results for Dirichlet partitions of  the $4d$ flat torus. 

For each of the flat tori considered, we have fixed a periodic box and the value $k$ and approximate the optimal partition. This is an easier problem than determining the optimal partition as $k\to \infty$. It has been observed that, for two-dimensional domains, as $k\to \infty$, a regular tiling of hexagons is optimal \cite{bourdin2010optimal}. In four dimensions,  our computational study suggests that, as $k\to \infty$, a regular 24-cell honeycomb is a good candidate minimizer. 

\section{A diffusion generated method for computing  Dirichlet partitions} \label{s:num}
In this section we first describe a mapping reformulation of the Dirichlet partitioning problem, \eqref{eq:ContDirPart}. 
Motivated by the formulation of the problem, we introduce an efficient diffusion generated method for computing Dirichlet partitions; see Section \ref{s:DiffGenMeth}. 
A brief comparison of our method with previous methods is given in Section \ref{s:previous}

\subsection{Mapping reformulation of Dirichlet partitions} \label{s:MappingReform}
Let $\Sigma_k$ denote the coordinate axis in $\mathbb R^k$ and define the Sobolev space 
$$ 
H^1_0(U; \Sigma_k) = \{\mathbf{u} \in H^1_0(U; \mathbb{R}^k) \colon \mathbf{u}(x) \in \Sigma_k \text{ a.e.}\} . 
$$ 
Since at most one component of a vector $\mathbf{v} \in \Sigma_k$ is non-zero, it follows that if $\mathbf{u} \in H^1_0(U; \Sigma_k)$ is continuous, then the sets $U_\ell = u_\ell^{-1}\left( \mathbb R \setminus \{0\} \right)$ partition $U$.  

The Dirichlet partition problem for $U$ is equivalent to the mapping problem
\begin{equation}\label{eq:MapProb}
  \min
  \left\{\mathbf{E}(\mathbf{u}) \colon \mathbf{u} = (u_1,  \ldots, u_k) \in
    H^1_0(U; \Sigma_k), \int_U u_\ell^2(x)  \ dx = 1 \ \forall \   \ell \in [k]\right\}, 
 \end{equation}
where 
$\mathbf{E}(\mathbf{u}) := \sum_{\ell=1}^k\int_{U}|\nabla u_\ell|^2 \ dx$
 is the Dirichlet energy of $\mathbf{u}$ \cite{caffarelli2007optimal}. We refer to a solution of \eqref{eq:MapProb} as a \emph{ground state of $U$}, which, without loss of generality, we may assume to be nonnegative. In particular, if $\mathbf{u}$ is a quasi-continuous representative of a ground state such that each component function $u_\ell$ assumes only nonnegative values, then a Dirichlet partition $U = \amalg_\ell U_\ell$ is given by $U_\ell = u_\ell^{-1}(0,\infty)$ for $\ell=1,\ldots, k$. Likewise, the first Dirichlet eigenvectors $u_\ell$ of a Dirichlet partition $\amalg_\ell U_\ell$ may be assembled into a function $\mathbf{u} \in H^1_0(U ;\Sigma_k)$ that solves the mapping problem \eqref{eq:MapProb}.

\subsection{Computational methods for Dirichlet partitions} \label{s:DiffGenMeth}
We consider the mapping formulation for Dirichlet partitions, \eqref{eq:MapProb}, for which there are basically three ingredients: 
(i) the Dirichlet energy,
(ii) the constraint that $u(x) \in \Sigma_k$, and
(iii) the constraint that $\int_U u_\ell^2 = 1$. 
Algorithm~\ref{a:MBO} iteratively handling these three ingredients. We begin with an initial vector valued function $\mathbf{u}^0 \in H^1_0(U;\mathbb R^k)$. 
Since $\Sigma_k \subset \mathbb R^k$, we can consider the unconstrained gradient flow of the Dirichlet energy until time $\tau$, which is exactly the evolution by the diffusion equation, given in the Diffusion Step of Algorithm~\ref{a:MBO}. Let $\tilde u_\ell(x) = u_\ell(\tau,x)$ denote the solution at time $\tau>0$. 
Next, for each point $x \in U$, we consider the closest value in $\Sigma_k$ to $\tilde u(x)$. This is exactly the Projection Step of Algorithm~\ref{a:MBO}. In this step, a rule should be devised to break any ties, but in practice we do not observe any. 
Finally, we renormalize each component of the vector to satisfy the $L^2(U)$ constraint as in the Renormalization Step of Algorithm~\ref{a:MBO}. 
These three steps are iterated until the condition that the partitions memberships didn't change in the previous iteration, i.e.,
\begin{equation} \label{e:ConvCrit}
\sum_{\ell \in [k]} \| \chi_{\{u_\ell^s>0\}} -\chi_{\{u_\ell^{s-1}>0\}}\| = 0
\end{equation}
where $\chi_{\{\cdot\}}$ denotes the indicator function. 

We refer to this algorithm as ``diffusion generated'' as it contains a diffusion step, similar to the Merriman-Bence-Osher (MBO) diffusion generated motion for approximating mean curvature flow \cite{merriman1992diffusion,merriman1994motion,MBO1993}. This method has subsequently been extensively analyzed and extended to more general contexts; see \cite{Ruuth_2001,esedoglu2015threshold,osting2017generalized}.

\begin{algorithm}[t!]
\DontPrintSemicolon
 \KwIn{Let $U$ be a $d-$dimensional Euclidean subset or a closed surface, $\tau > 0$ be a time-parameter, and $\mathbf{u}^0 \in H^1(U;\mathbb{R}^k)$ be an initial condition.}
 \KwOut{ An approximate ground state, $\mathbf{u}^s \in H^1(U;\mathbb{R}^k)$, satisfying \eqref{eq:MapProb}. }
 \For{ $s=1,2,\ldots$}{
{\bf 1.  Diffusion Step.} Solve the initial value problem for the diffusion equation  until time $\tau$ with initial value given by each of the components of $\mathbf{u}^{s-1}(x)$:
\begin{align*}
&\partial_t u_\ell(t,x) = \Delta u_\ell(t,x) \\
&A(0,x) = u^{s-1}_\ell(x).
\end{align*}
Let $\tilde u_\ell(x) = u_\ell(\tau,x)$. \;
{\bf 2. Projection Step.} Set \begin{align*}
u^{*}_\ell(x) =  
\begin{cases}
 \tilde u_\ell(x) & \text{if}  \  \tilde{u}_\ell(x) = \max\limits_{j \in [k] } \tilde{u}_j(x) \\
0 & \text{otherwise} 
\end{cases}. 
\end{align*}
{\bf 3. Renormalization Step.} Set $u_\ell^s(x) = \frac{u^{*}_\ell(x)}{\| u^{*}_\ell(x) \|}$ where $\| \cdot \|$ denotes the $L^2(U)$ norm. \;
\If{ \eqref{e:ConvCrit} is satisfied,} {STOP\;}
}
\caption{A diffusion generated algorithm for computing Dirichlet partitions. } 
\label{a:MBO}
\end{algorithm}

\subsection{Comparison with other methods for computing Dirichlet partitions} \label{s:previous}
There are variety of approaches to computing Dirichlet partitions, which we organize by the way in which the energy \eqref{eq:ContDirPart}, or equivalently \eqref{eq:MapProb}, is relaxed. 

One relaxation of the constraint $\mathbf{u}(x) \in \Sigma_k$ is the following. Consider the function $f\colon \mathbb R^k \to \mathbb R$, given by 
$f(x) = \sum_{i \neq j}^k x_i^2x_j^2$. It is not difficult to see that $f(x) \geq 0$ and $\Sigma_k =  f^{-1}(0)$. For $\varepsilon > 0$, we can consider the relaxation of \eqref{eq:MapProb}, given by 
\begin{equation}\label{eq:MapProb2}
  \min
  \left\{\mathbf{E}^\varepsilon(\mathbf{u}) \colon \mathbf{u} = (u_1,  \ldots, u_k) \in
    H^1_0(U; \mathbb{R}^k), \int_U u_\ell^2(x)  \ dx = 1 \ \forall \   \ell \in [k]\right\}, 
 \end{equation}
where the relaxed energy is given by 
$\mathbf{E}^\varepsilon(\mathbf{u}) :=\mathbf{E}(\mathbf{u}) +  \frac{1}{2 \varepsilon^2} \int_U  f( \mathbf{u}(x) ) \ dx$.
Properties of this relaxation can be found in \cite{bao2004computing,caffarelli2007optimal} and it was used to devise computational methods  in \cite{bao2004computing,Du_2008,elliott2015computational}.

In particular, in \cite{Du_2008}, Q. Du and F. Lin introduce a three-step diffusion generated motion similar to the one considered in Algorithm~\ref{a:MBO}. However, in the second step, rather than taking the closest point in $\Sigma_k$, the following system of ODEs is solved until time $\tau$,
$$
\frac{d}{dt} \tilde{u}_\ell = \frac{1}{\varepsilon^2} \left( \sum_{j \neq \ell} \tilde{u}^2_j \right) \tilde{u}_\ell, 
\qquad \qquad   
\ell \in [k]. 
$$
This is precisely the gradient flow of the second term of the relaxed energy $\mathbf{E}^\varepsilon$. 
Numerically, this system is solved using the  Gauss-Seidel method. However, the small parameter $\varepsilon$ here restricts the mesh size and fats the interface between any two partitions. Also, the authors only considered 2-dimensional case there.

Another approach, developed in first \cite{bourdin2010optimal}, is based on a Schr\"odinger operator relaxation of \eqref{eq:ContDirPart} and was further used in \cite{osting2013minimal,Zosso2015,Bogosel_2016,bogosel2017efficient}. 

Other related ideas based on a stochastic interpretation can be found in  \cite{cybulski2005,cybulski2008}.

\section{Implementation of Algorithm~\ref{a:MBO}} \label{s:comp}
In this section, we describe a numerical implementation of Algorithm~\ref{a:MBO} for $d$-dimensional flat tori and spheres. Although Algorithm~\ref{a:MBO} could in principle be implemented in more generality, our implementation relies on the Fast Fourier Transform (FFT) or Spherical Harmonic Transform (SHT).

\subsection{Implementation for flat tori}\label{ss:flat}
In this section, we consider the implementation of Algorithm~\ref{a:MBO} on the computational domain $\Omega = [-1,1]^d$ ($d = 2, 3,4$) with edges identified (periodic boundary conditions). 

The diffusion step in Algorithm~\ref{a:MBO} for partition $\ell$ is to solve
\begin{subequations} \label{e:PeriodicHeat}
\begin{align}
&\partial_t u_\ell(t,x) = \Delta u_\ell(t,x)  && x\in \Omega, \ t\geq 0,\\
&u_\ell(0,x) = u_\ell^{s-1}(x) && x \in \Omega \\
& u_\ell \textrm{ satisfies periodic boundary conditions on } \partial \Omega. 
\end{align}
\end{subequations}
It is well-known that the solution for the diffusion equation for a scalar function on $\mathbb R^d$ at time $t=\tau$ can be expressed as the convolution of the  heat kernel, 
\begin{equation*} \label{e:HeatKer}
G_{\tau}^d(x) = ( 4 \pi \tau )^{-\frac d 2 } \exp\left( -\frac{|x|^2}{4 \tau} \right), 
\end{equation*}
 and the initial condition, $u_\ell^{s-1}(x)$. For our periodic domain, $\Omega \subset \mathbb R^d$, we denote by $G_{p,\tau}$ the periodic heat kernel, given by 
  \[ G_{p,\tau}(x) = \sum_{\alpha \in \mathbb Z^d }  G_\tau^d(x- \alpha). \]
 The solution, $\tilde{u}_\ell(x) = u_\ell(\tau,x)$ to \eqref{e:PeriodicHeat} at time  $t=\tau$ has matrix components given by $\tilde{u}_\ell=G_{p,\tau} * u_\ell^{s-1}$, where $*$ denotes the convolution. 

We denote the Fourier transform and its inverse by $\mathcal{F}$ and $\mathcal{F}^{-1}$, respectively. Using the convolution property that 
$\mathcal{F}(G_{p,\tau}*u^{s-1}_\ell) =\mathcal{F}(G_{p,\tau}) \ \mathcal{F}(u^{s-1}_\ell) $, we can express the solution to \eqref{e:PeriodicHeat} as 
\[ 
\tilde{u}_{\ell} = \mathcal{F}^{-1} \left( \ \mathcal{F}(G_{p,\tau}) \ \mathcal{F}( u^{s-1}_\ell) \ \right). 
\] 

In our numerical implementation, due to the periodic boundary condition, can efficiently compute an approximation to the Fourier transform and its inverse using the  fast Fourier transform (FFT) and inverse fast Fourier transform (iFFT). That is, an approximate solution to \eqref{e:PeriodicHeat} is evaluated via 
\begin{align*}
\tilde{u}_\ell =  \textrm{iFFT} \left( \ \textrm{FFT}(G_{p,\tau}) \ \textrm{FFT}(u^{s-1}_\ell) \ \right).
\end{align*}

It is well known that the computational complexity of the FFT is $O( n^d \log n)$ where $n$ is the number of grid points in each direction. The total computational complexity of this Algorithm \ref{a:MBO} is then 
$$
(\# \textrm{steps}) \cdot k \cdot  O( n^d \log n). 
$$
 
\subsection{Implementation on a spherical surface} \label{ss:sphericalsurface}

In this section, we consider the implementation on the computational domain $\Omega = S^2$ which is a spherical surface with radius $1$. Here is is  understood that $\Delta$ is the Laplace-Beltrami operator on the spherical surface. We parameterize  $S^2$ in spherical coordinates, 
\begin{equation} \label{e:SphereCoord}
(x,y,z)= (\sin \theta \sin \phi, \sin \theta \cos \phi, \cos \theta),
\end{equation}
where  $\theta\in[0,\pi]$ is the inclination and $\phi\in[0,2\pi]$ is the azimuth. 
It is well known that the eigenfunctions of the Laplace Beltrami operator on the spherical surface are the spherical harmonic functions, $Y_{l}^{m}$, where 
\[
Y_{l}^{m}(\theta,\phi)  =  \sqrt{\frac{(2l+1)}{4\pi}\frac{(l-m)!}{(l+m)!}}P_{l}^{m}(\cos(\theta))e^{im\phi}, 
\qquad \quad l \in \mathbb N, \ m \in \{-l, \ldots, l\}. 
\]
with the corresponding eigenvalues being $-l(l+1)$.
Denote $SHT$ as the spherical harmonic transform and $iSHT$ as the inverse spherical harmonic transform. For each partition $\ell$, the initial condition $u_\ell(t=0,x)$ can be decomposed by $n^2$ spherical harmonic functions:
\[u_\ell(t=0,x) = \sum_{l=0}^n\sum_{m=-l}^{l} s_{l,m}^{\ell} Y_{l}^{m}.\]
Using the spherical harmonic functions to express the solution of the surface diffusion equation at $t=\tau$, the coefficients are given by $s_{l,m}^{\ell} e^{-l(l+1)\tau}$. The solution to the diffusion equation can be computed by the inverse spherical harmonic transform, 
\[u_\ell(\tau,x) = iSHT\left( SHT(u_\ell(0,x)) \ e^{-l(l+1)\tau} \right).\]

\section{Numerical results} \label{s:numeircs}
In this section, we use the implementation of Algorithm \ref{a:MBO}, described in Section \ref{s:num}, to compute approximate Dirichlet partitions. 
The algorithms are implemented in MATLAB. For the results in two, three, and four dimensional periodic space, we used fast Fourier transform (FFT) to solve the heat diffusion equation; see in Sections~\ref{s:2d}, \ref{s:3d}, and \ref{s:4d}. For the results on the sphere, we used the spherical harmonic transform to solve the surface diffusion equation on a spherical surface; see Section~\ref{s:sphere}. 
 For all numerical results, we initialize the algorithm by computing the Voronoi tessellation for a random point set in $U$ and use the normalized indicator functions for this tessellation. Below, we simply refer to this as \emph{initializing using a random tessellation}.  All reported results were obtained on a laptop with a 2.7GHz Intel Core i5 processor and 8GB of RAM. 

To compare the energies between configurations and for different size domains and values of $k$, we consider the normalized energy 
\begin{equation} \label{e:Energy}
E(k, U) := \min_{U = \cup_\ell U_\ell} \ \  
\frac{|U|^{\frac{2}{d}}}{k^{1+\frac{2}{d}}} \sum_{\ell=1}^k \lambda (U_\ell). 
\end{equation}
This quantity is invariant under homothety, {\it i.e.}, $E(k,\alpha U) = E(k,U)$ and has the property that for $m \in \mathbb N$, 
$$
E(m^d k, m U) = 
\min_{\cup_\ell U_\ell} \ \frac{ m^2 |U|^{\frac{2}{d}}}{m^{d+2} k^{1 + \frac{2}{d}}} \sum_{\ell=1}^{m^d k}  \lambda (U_\ell) 
= \min_{\cup_\ell U_\ell} \frac{ |U|^{\frac{2}{d}}}{ k^{1 + \frac{2}{d}}} \frac{1}{m^d} \sum_{\ell=1}^{m^d k}  \lambda (U_\ell) 
\leq E(k,U).
$$
where the last inequality comes from repeating the $k$-Dirichlet partition on $U$---$m$ times in each direction--- to form a $m^d k $-Dirichlet partition on $m U$. 
We report values for an approximation of $E$ in \eqref{e:Energy}, given by 
\begin{align}
\label{e:aEnergy}
\tilde{E}(k, U) 
&:= \frac{|U|^{\frac{2}{d}}}{k^{1 + \frac{2}{d}}} \frac{1}{\tau} \left(k- \sum_{\ell=1}^k \langle u_\ell, e^{\Delta \tau} u_\ell \rangle \right) \\
\nonumber
&\approx \frac{|U|^{\frac{2}{d}}}{k^{1 + \frac{2}{d}}} \sum_{\ell=1}^k \langle u_\ell, -\Delta u_\ell \rangle, 
\end{align}
where the $\{u_\ell\}_{\ell \in [k]}$ have unit $L^2(U)$ norm. 
See \cite{esedoglu2015threshold,osting2017generalized} for more intuition on this approximate energy.

\subsection{2d flat torus} \label{s:2d}
It was proven by T. Hales that the regular hexagon tessellation is the equal-area partition that minimizes surface area  \cite{Hales_2001}. In two-dimensional Euclidean space, it has been conjectured that this tessellation is also a Dirichlet partition  \cite{caffarelli2007optimal}. Computationally the problem of partitioning 2D rectangles, either with periodic or Dirichlet boundary conditions, has been addressed in 
 \cite{cybulski2005,Du_2008,bourdin2010optimal,osting2013minimal,bogosel2017efficient} and embedded tori have been studied in  \cite{elliott2015computational,bogosel2017efficient}. In all of these studies, for large values of $k$, regular hexagons are ubiquitous. 

\begin{figure}[ht]
\centering
\includegraphics[scale=0.15,clip,trim= 7cm 1cm 7cm 1cm]{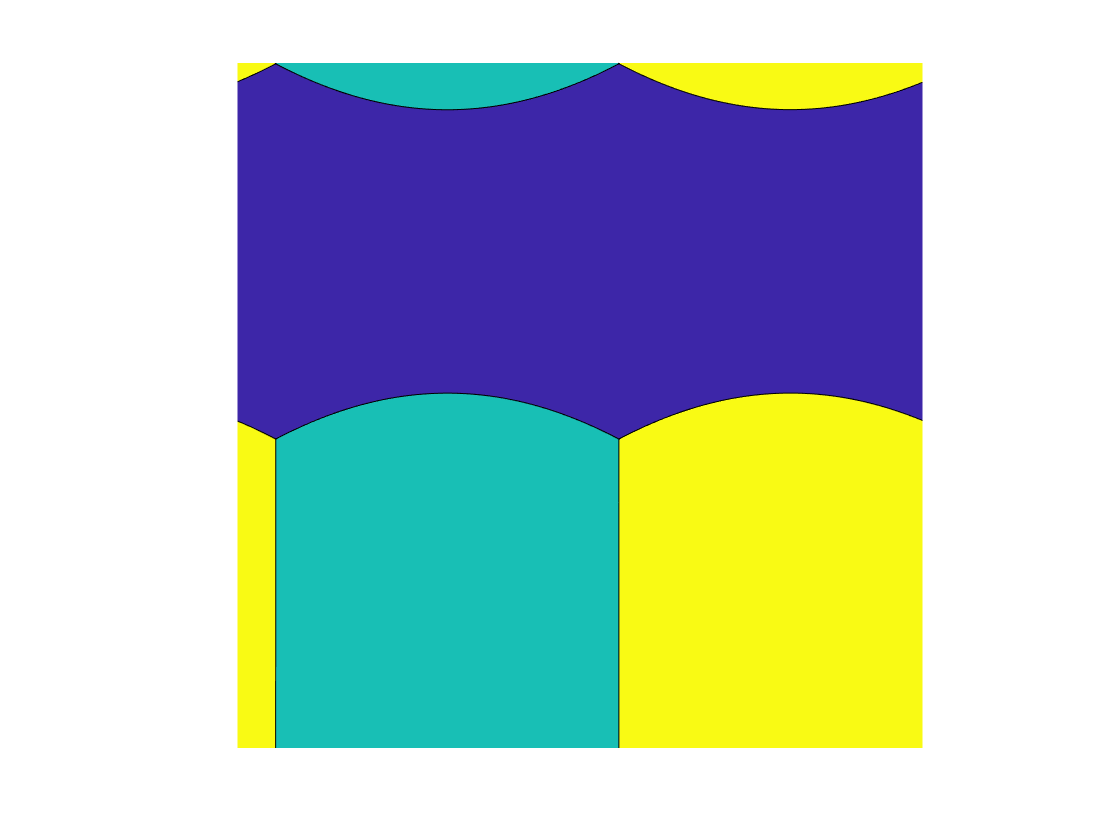}
\includegraphics[scale=0.15,clip,trim= 7cm 1cm 7cm 1cm]{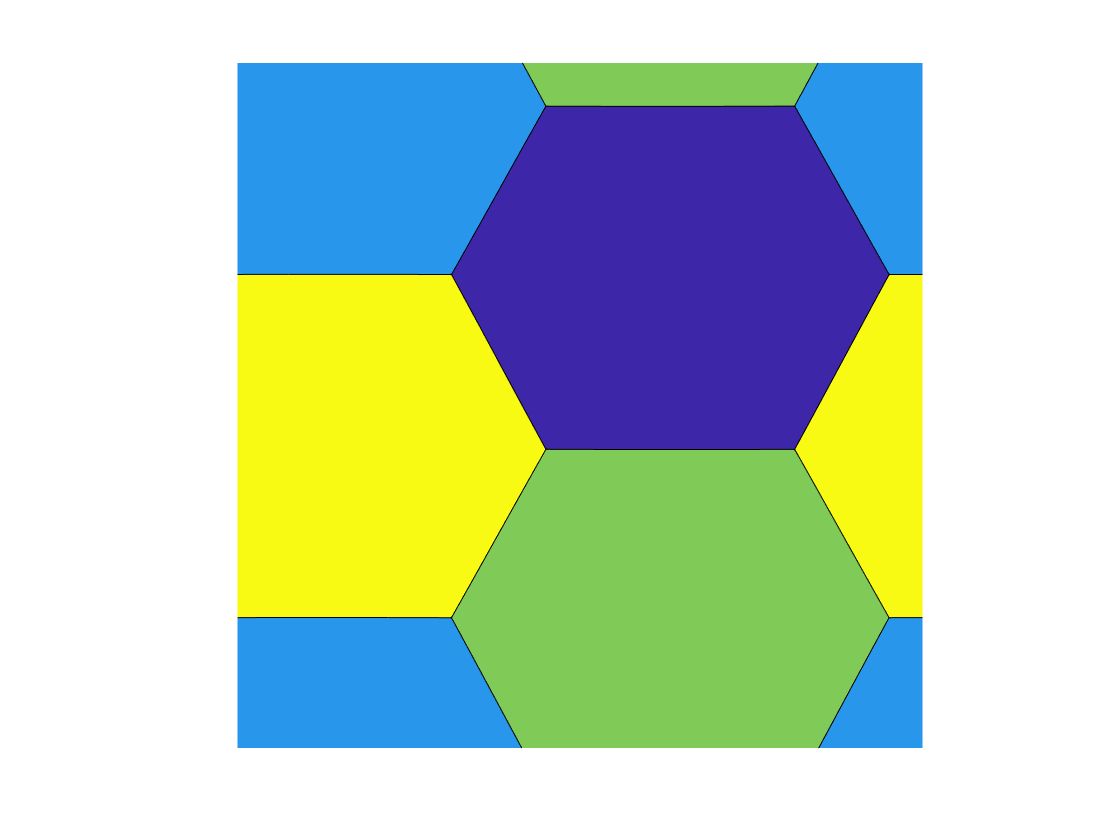}
\includegraphics[scale=0.15,clip,trim= 7cm 1cm 7cm 1cm]{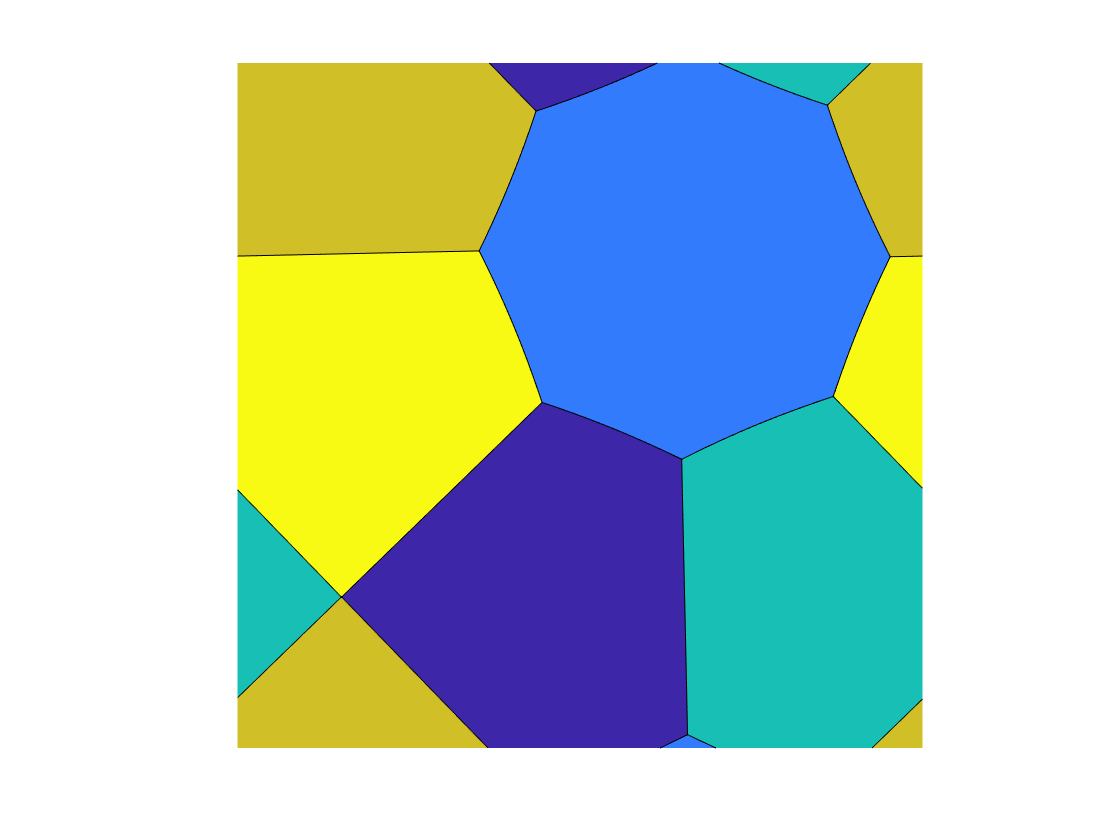}\\
\includegraphics[scale=0.15,clip,trim= 7cm 1cm 7cm 1cm]{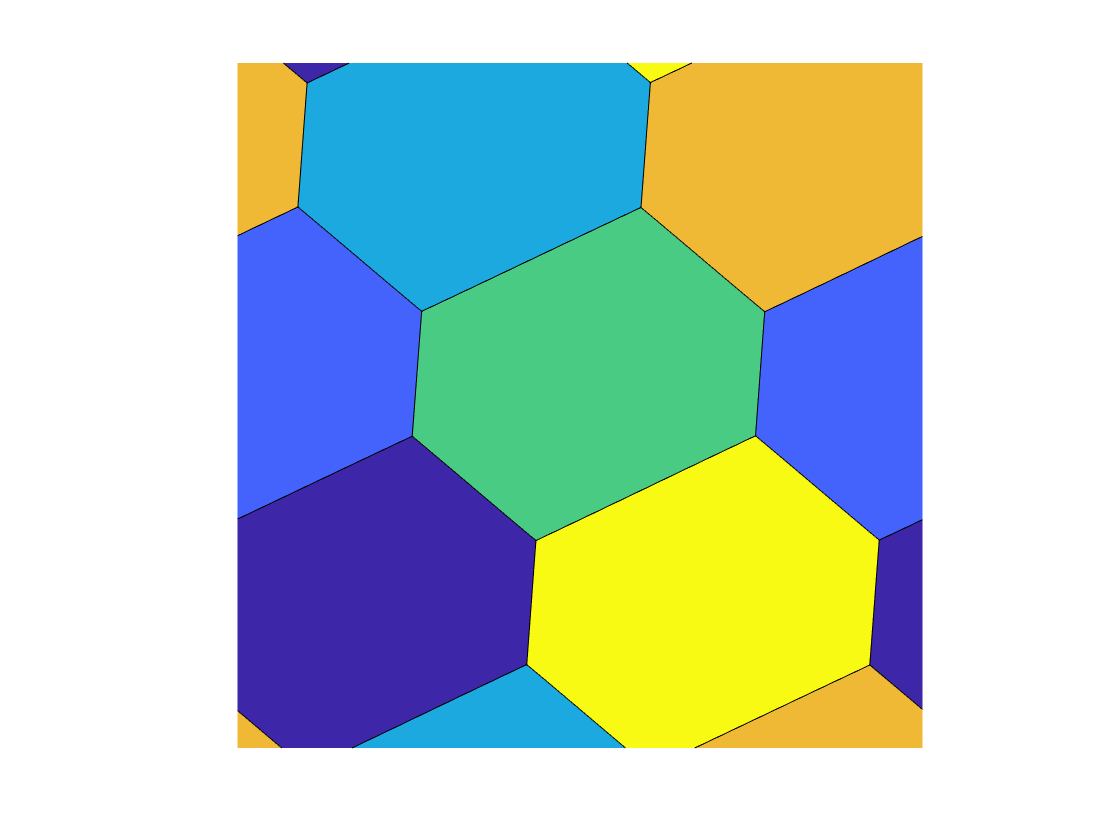}
\includegraphics[scale=0.15,clip,trim= 7cm 1cm 7cm 1cm]{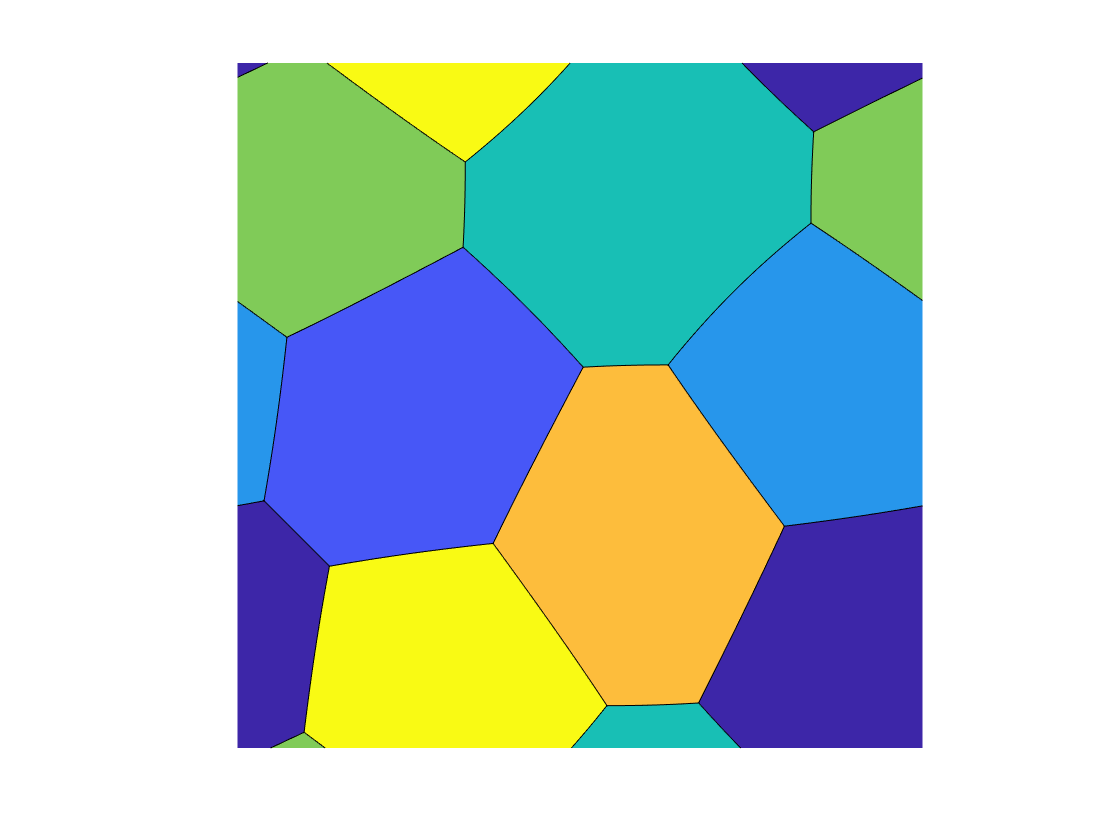}
\includegraphics[scale=0.15,clip,trim= 7cm 1cm 7cm 1cm]{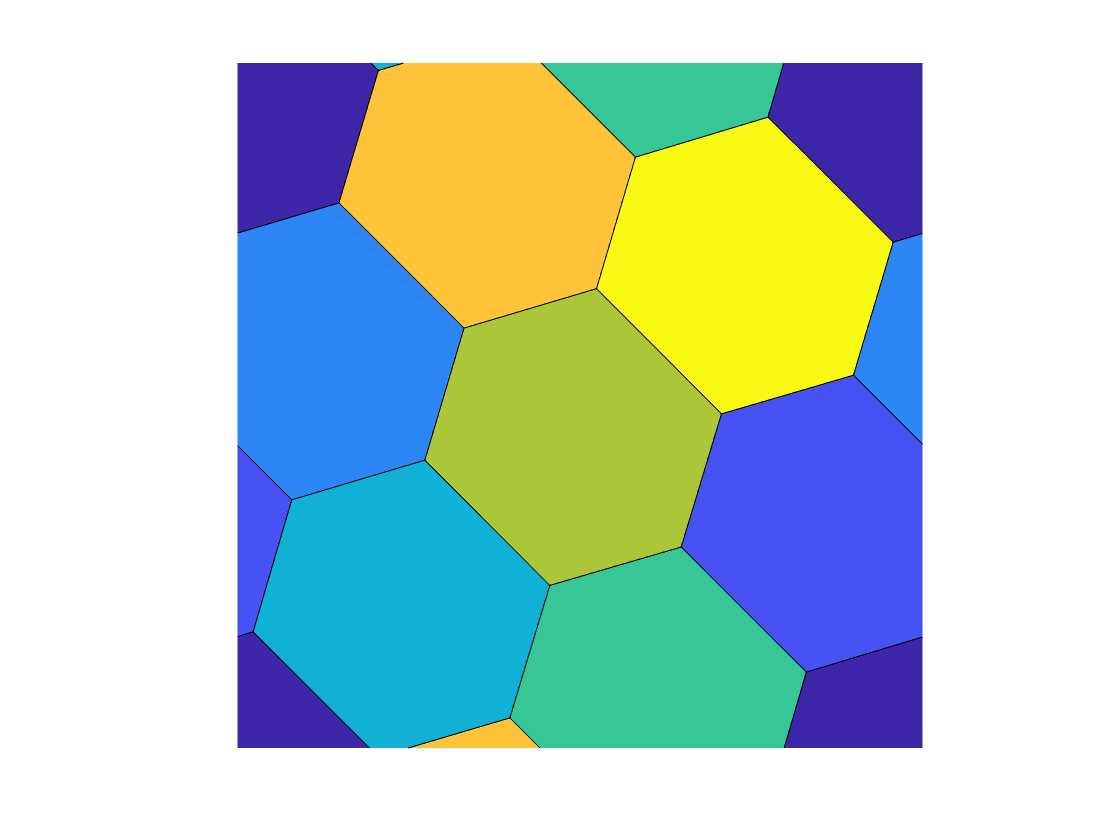}\\
\includegraphics[scale=0.15,clip,trim= 7cm 1cm 7cm 1cm]{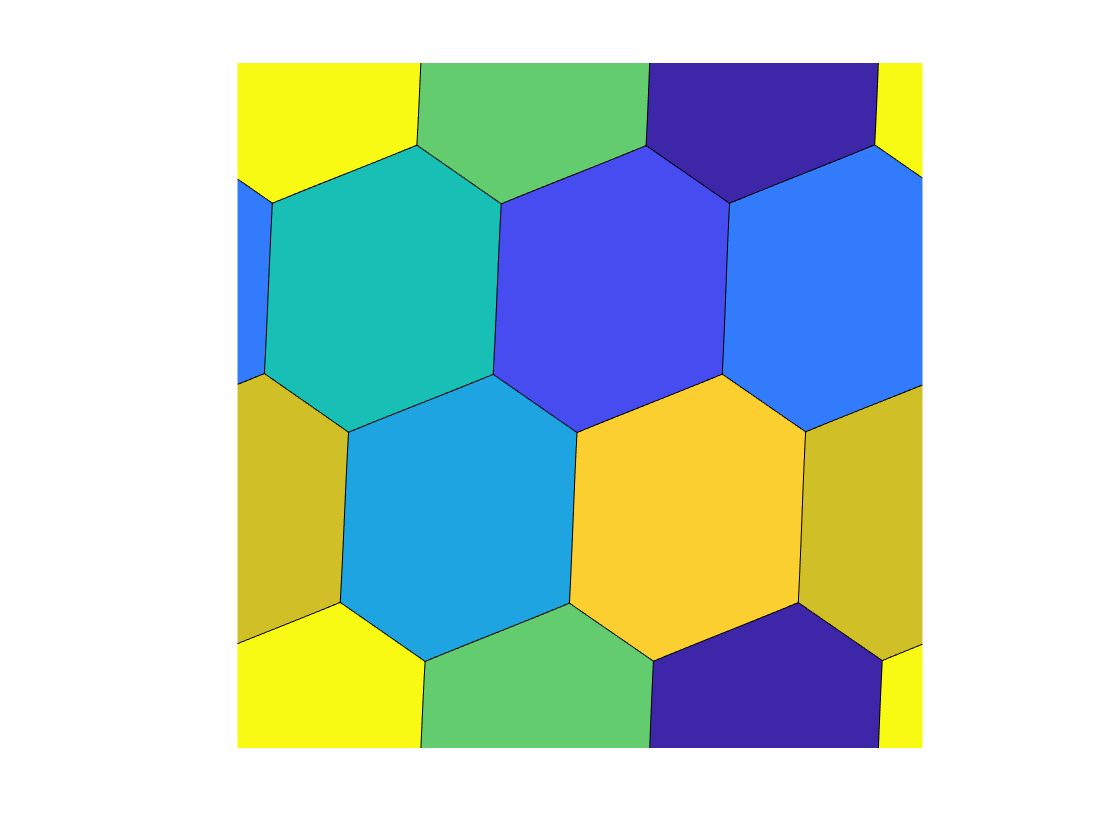} \includegraphics[scale=0.15,clip,trim= 7cm 1cm 7cm 1cm]{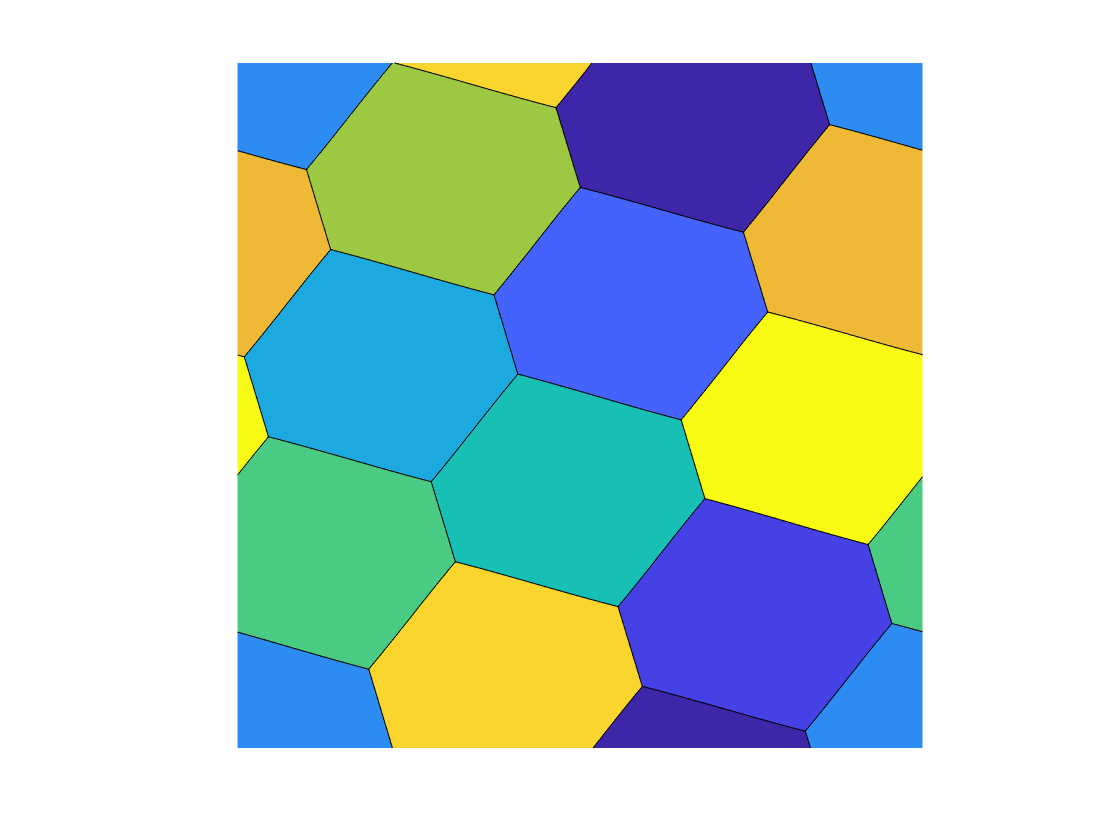}    \includegraphics[scale=0.15,clip,trim= 7cm 1cm 7cm 1cm]{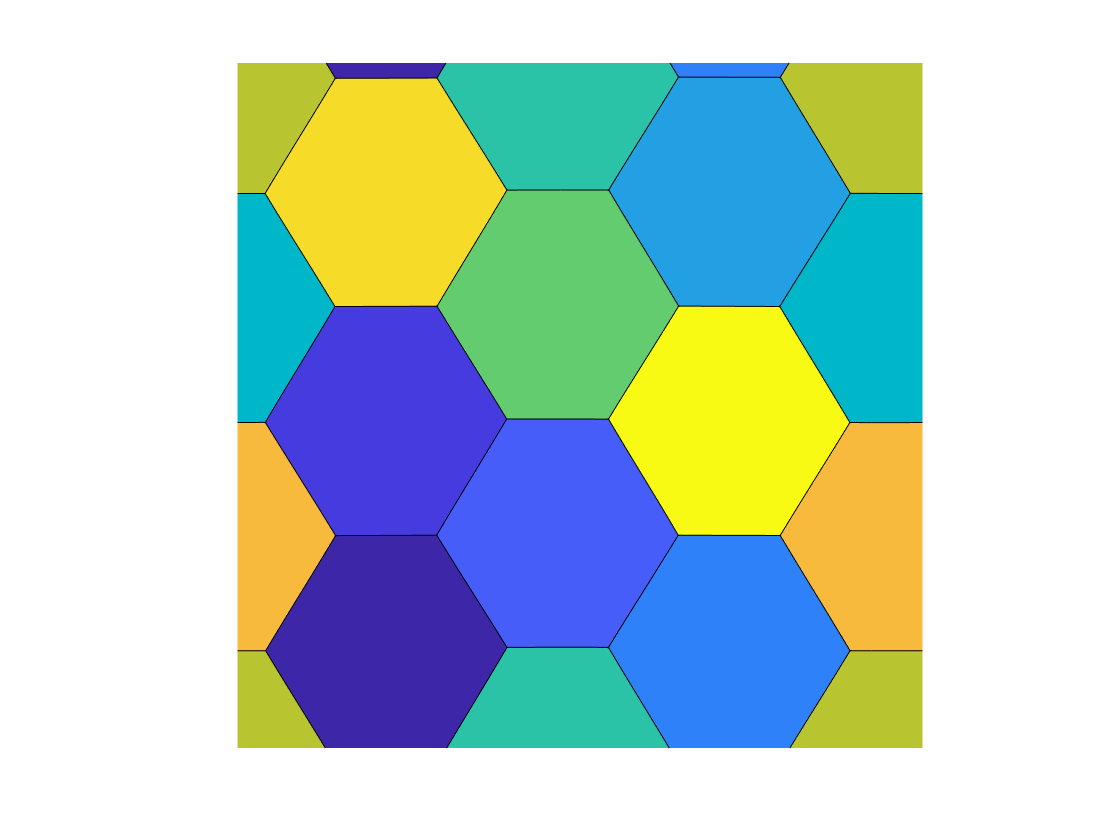}\\
\includegraphics[scale=0.15,clip,trim= 7cm 1cm 7cm 1cm]{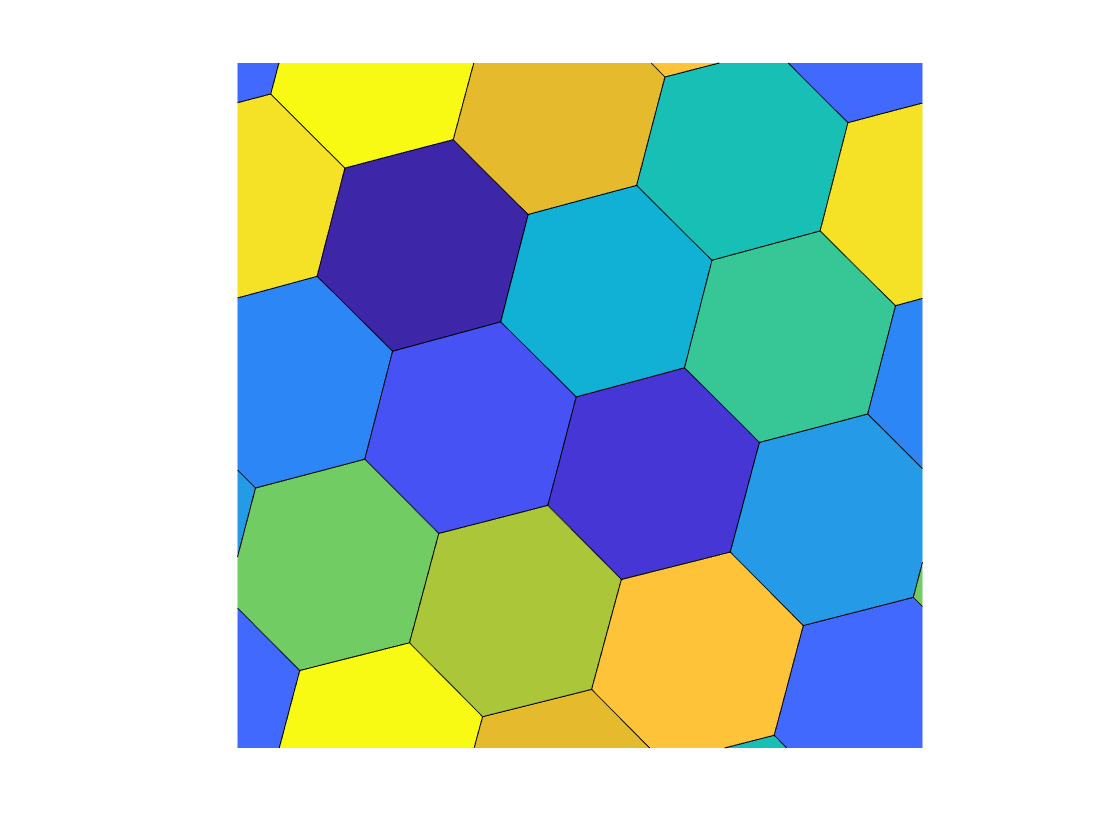}
\includegraphics[scale=0.15,clip,trim= 7cm 1cm 7cm 1cm]{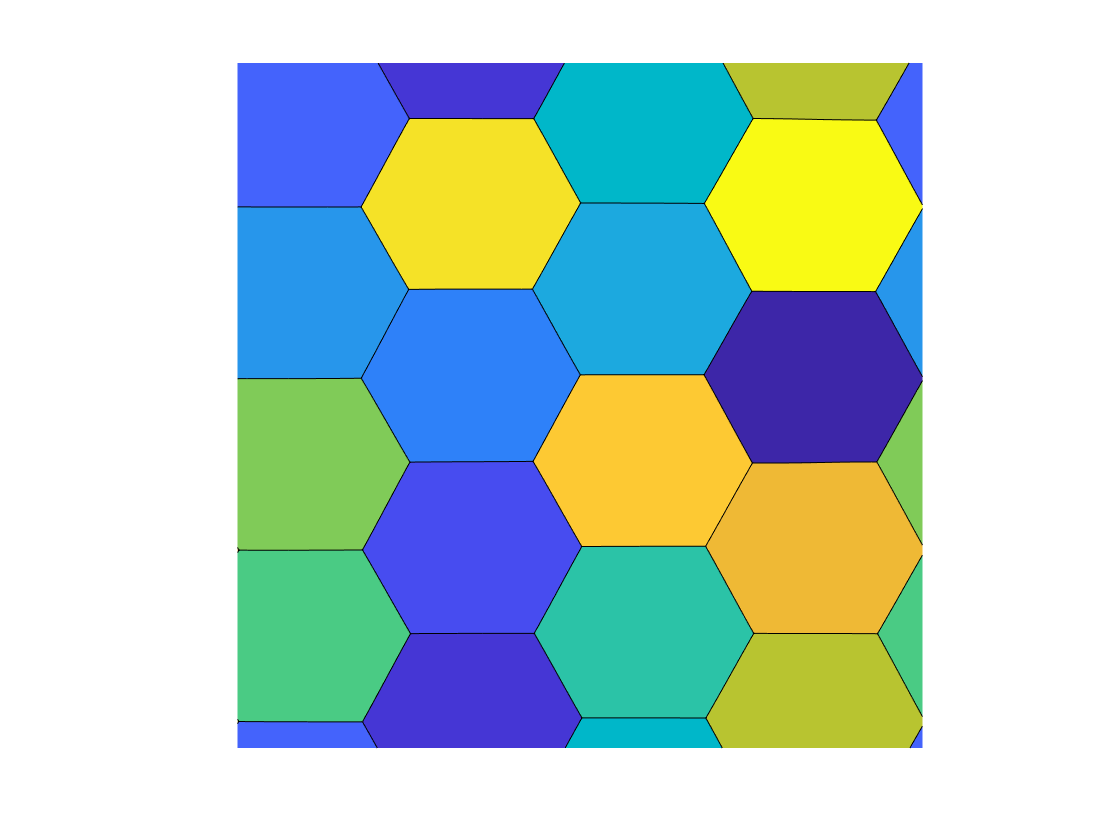}
\includegraphics[scale=0.15,clip,trim= 7cm 1cm 7cm 1cm]{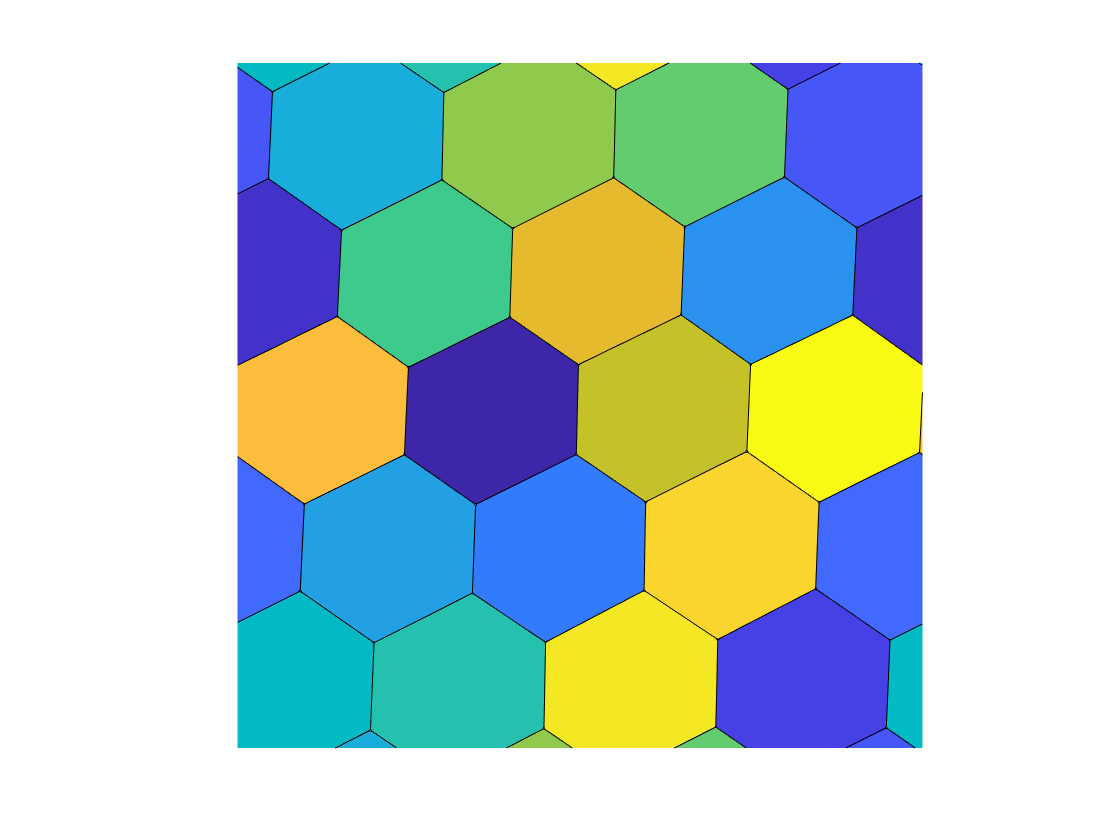}
\caption{From left to right and top to bottom: Dirichlet partitions on the $[-1,1]^2$ periodic domain discretized by $256^2$ uniform grid points with $k=$3--9,11,12,15,16, and 20. The last one is computed using $\tau = 0.0625$ while others are all computed using $\tau = 0.125$. The average CPU time for each case is $3.02$, $1.89$, $5.09$, $3.49$, $6.89$, $6.36$, $9.89$, $11.02$, $8.42$, $16.18$, $21.45$, and $35.38$ seconds respectively.}
\label{fig:2d_1}
\end{figure}

 In Figure~\ref{fig:2d_1}, we display Dirichlet partitions for the $[-1,1]^2$ periodic domain discretized by $256^2$ uniform grid points with $k=3-9,11,12,15,16$ and $20$. The code was executed several times initialized using random $k$-tessellations. For these values of $k$, the algorithm always converges to the same pattern, but for larger values of $k$, we observe local minima.  In this experiment, we use $\tau = 0.0625$ for $k=20$ and $\tau = 0.125$ for all other values of $k$. In Table~\ref{tab:2d}, we display the average CPU time for each value of $k$. Here, the average CPU time is calculated by averaging the CPU time for each of the $10$ experiments (with random initial conditions).  
 
\begin{figure}[ht]
\centering
\includegraphics[scale=0.23,clip,trim= 7cm 1cm 7cm 1cm]{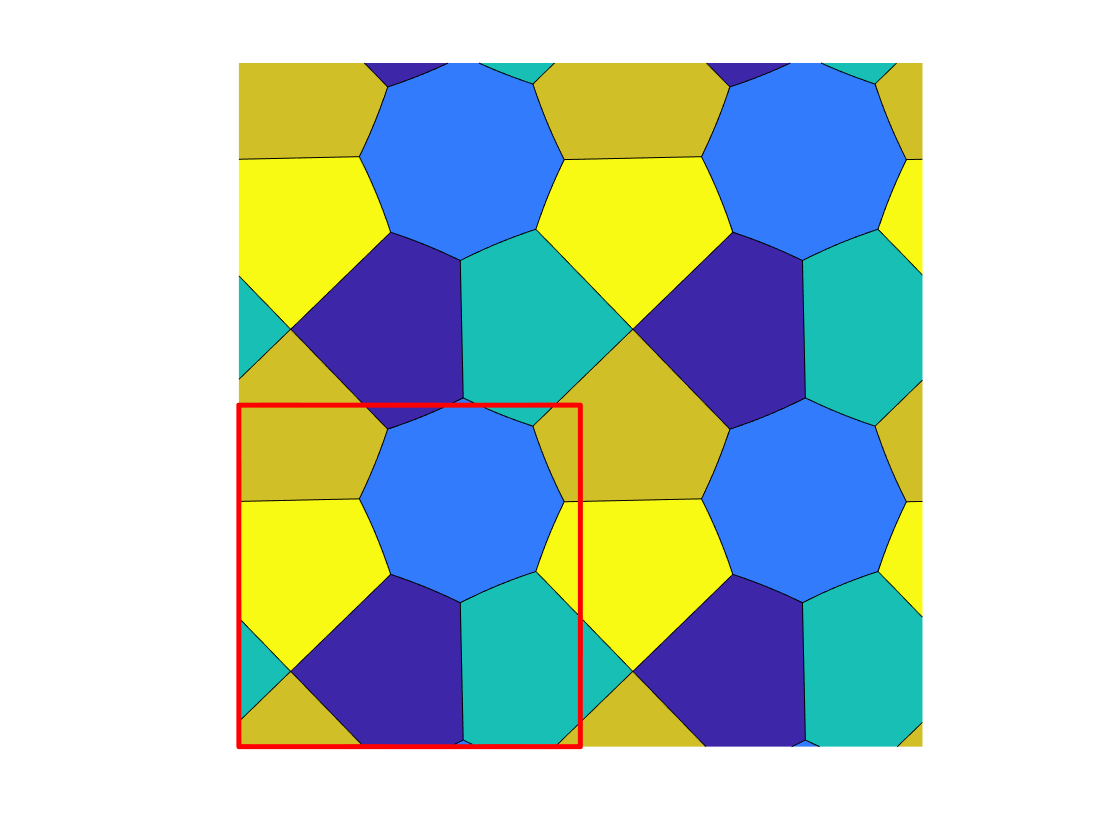}
\includegraphics[scale=0.23,clip,trim= 7cm 1cm 7cm 1cm]{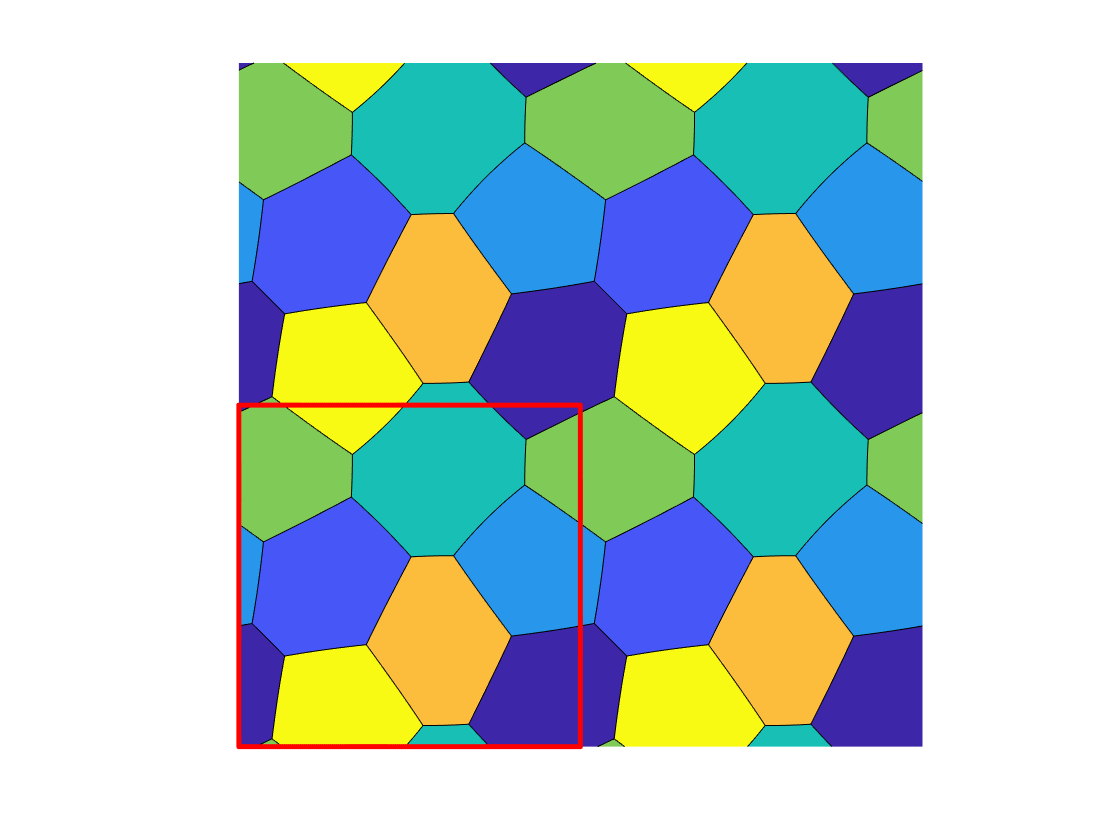}
\caption{ The periodic extension of the $k=5$ (left) and $k=7$ (right) Dirichlet partitions to a larger domain. 
In both panels, the red lines are the boundary of $[-1,1]^2$. See Figure \ref{fig:2d_1}. } 
\label{fig:2d_2}
\end{figure}

The partitions obtained are similar to those found previously. Since the domain has aspect ratio equal to one, regular hexagons cannot be used to tile the domain, so the hexagons are slightly distorted. To better see the irregular Dirichlet partitions for $k=5$ and $k=7$,  in Figure~\ref{fig:2d_2}, we plot their periodic extensions. These numerical results demonstrate that, although Algorithm \ref{a:MBO} is simple, it is efficient and stable. In Table~\ref{tab:2d}, we also tabulate the values of $\tilde E$ in \eqref{e:aEnergy} for different values of $k$.

\begin{table}[ht]
\centering
\caption{Values of  $\tilde E$ in \eqref{e:aEnergy} and the average CPU time for different values of $k$.} \label{tab:2d}
\begin{tabular}{|c|c|c|c|c|c|c|c|c|c|c|c|c|}
\hline
$k$ & 3 & 4& 5& 6 &7&8 \\
\hline 
$\tilde E$  & 2.39 & 2.13 & 2.23 & 2.18 & 2.17 & 2.09 \\
\hline
Average CPU time (s) & 3.02 & 1.89 & 5.09 & 3.49 & 6.89 & 6.36 \\
\hline
\hline
 $k$ &9&11& 12 & 15 &16&20\\
 \hline 
 $\tilde E$& 2.11 & 2.09 & 2.03 & 1.97 & 1.99 & 1.70 \\
 \hline
 Average CPU time (s) & 9.89 & 11.02 & 8.42 & 16.18 & 21.45 & 35.38\\
 \hline
\end{tabular}
\end{table}

\subsection{3d flat torus} \label{s:3d}
In three dimensions, the minimal total surface area partition is unknown. Lord Kelvin conjectured that a packing of  truncated octahedra was optimal \cite{thompson1887}. However, R. Weaire and D. Phelan discovered another structure comprised of two polyhedra which has a slightly smaller surface area  \cite{WeairePhelan}. For the three-dimensional Dirichlet partitioning problem, as far as we are aware, very little is known analytically and only a few papers have investigated the problem computationally \cite{cybulski2008,bogosel2017efficient}. Interestingly, both the Kelvin and the Weaire-Phelan structures appear as Dirichlet partitions, depending on the domain and value $k$. In this section, we compute Dirichlet partitions using Algorithm~\ref{a:MBO} for the  periodic cube, $[-1,1]^3$ and $k=2,4,8,16$. 

\begin{figure}[ht]
\centering
\includegraphics[scale=0.28,clip,trim= 8cm 2cm 8cm 2cm]{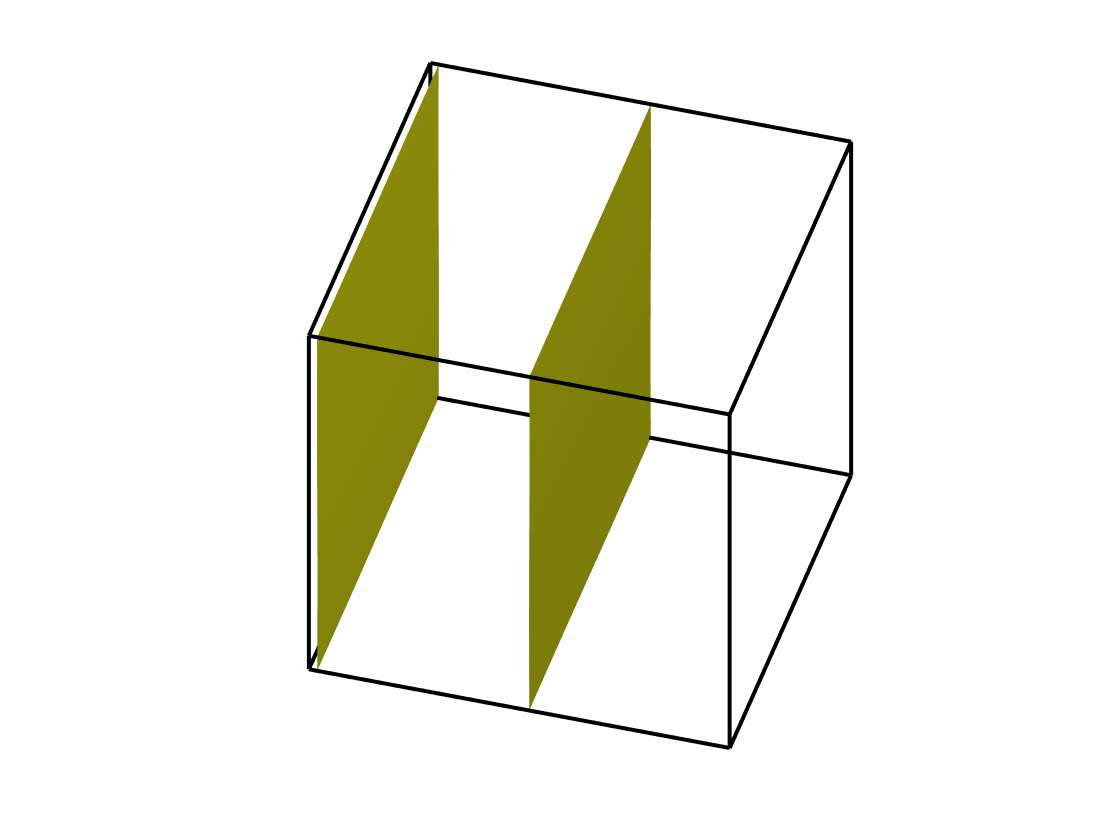}
\includegraphics[scale=0.28,clip,trim= 11cm 5cm 10cm 5cm]{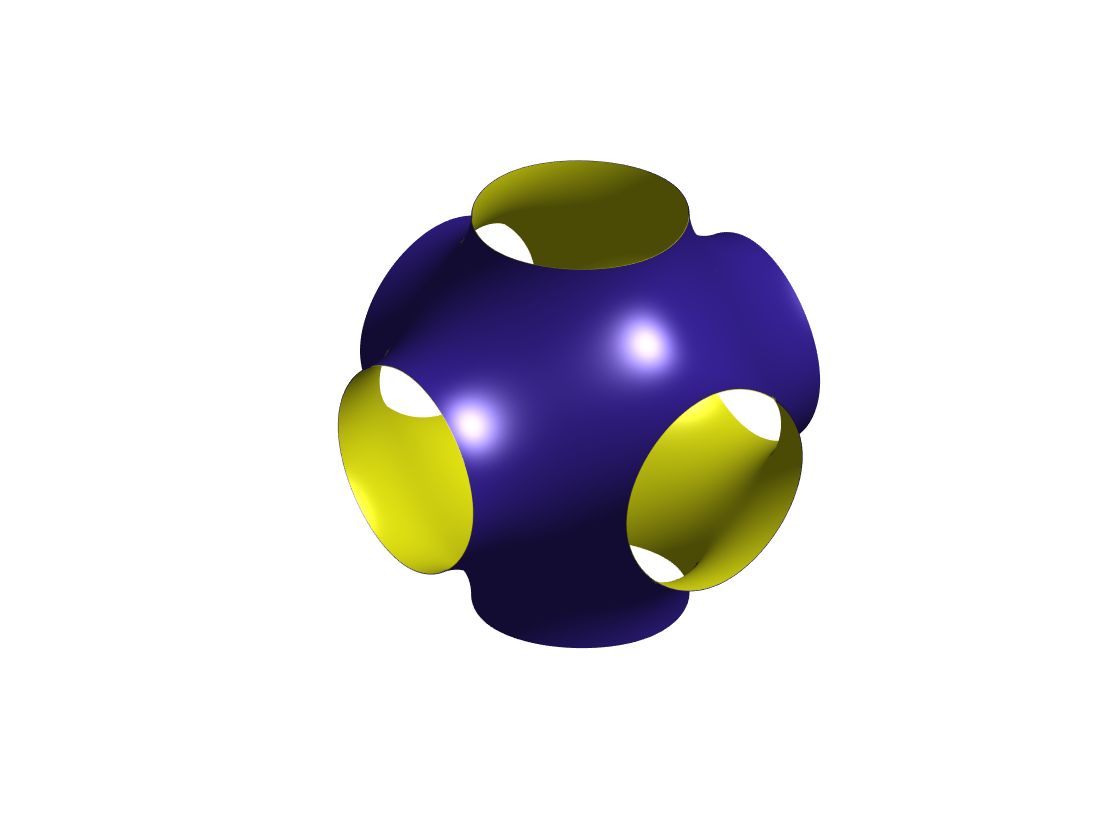}
\caption{{\bf (left)} A $k=2$ Dirichlet partition of the periodic cube $[-1,1]^3$ with interface given by parallel planes. 
{\bf (right)} The periodic cube $[-1,1]^3$ is partitioned into two components by a surface that is similar to the Schwarz P surface. 
In this experiment, the cube is discretized by $128^3$ uniform grid points and $\tau = 0.25$.  
The CPU time for the left case is $26$ seconds while the CPU time for the right case is $3$ seconds.}
\label{fig:3d_1}
\end{figure}

For $k=2$ and for every initialization using a random tessellation we tried, we obtained a partition given by a slab, which is shown in the left panel of Figure~\ref{fig:3d_1}. If we choose an initial condition so that the interface is the implicit equation of the surface, $\cos(x) + \cos(y) + \cos(z) = 0$, we obtain a partition that has interface that is similar to the Schwarz P surface, displayed in the right panel of Figure~\ref{fig:3d_1}.  These partitions are similar to ones reported in \cite{cybulski2008,Zosso2015}. In this experiment, the cube is discretized by $128^3$ uniform grid points and $\tau = 0.25$. The CPU time for the first one is $26$ seconds and the CPU time for the second one is $3$ seconds.

\begin{figure}[ht]
\centering
\includegraphics[scale=0.3,clip,trim= 12cm 5cm 10cm 5cm]{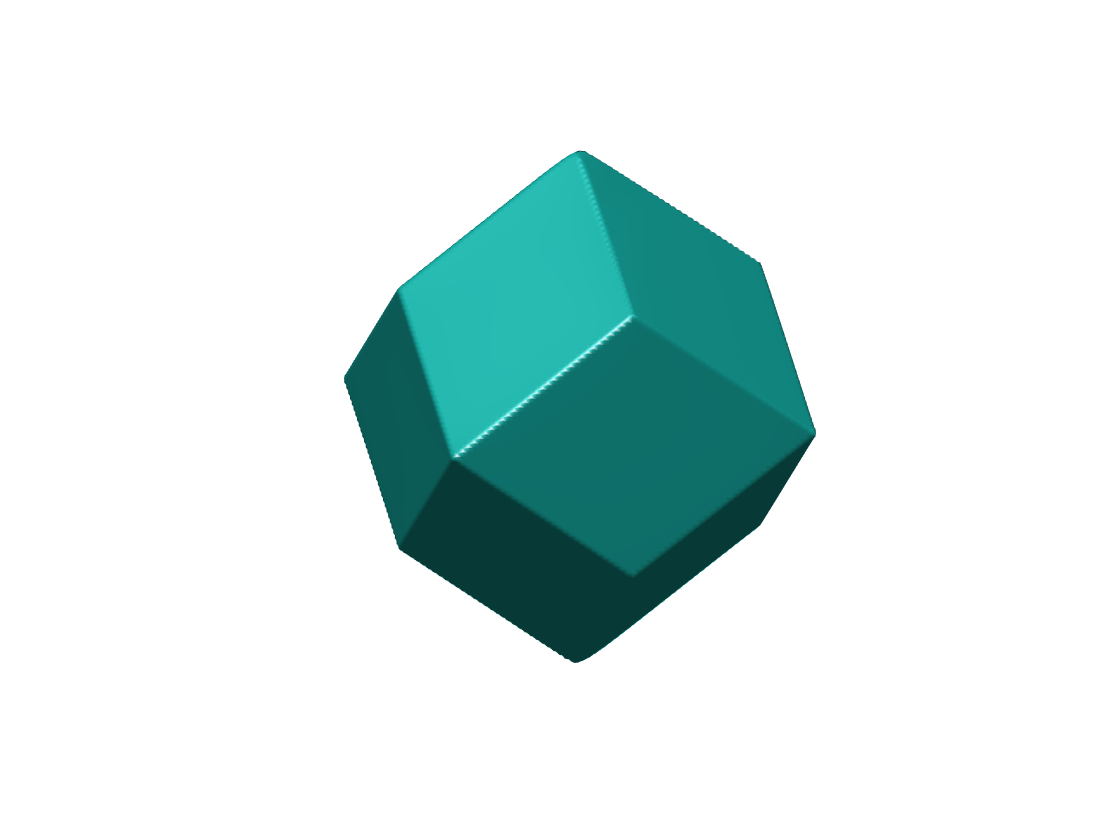}
\includegraphics[scale=0.2,clip,trim= 12cm 5cm 14cm 5cm]{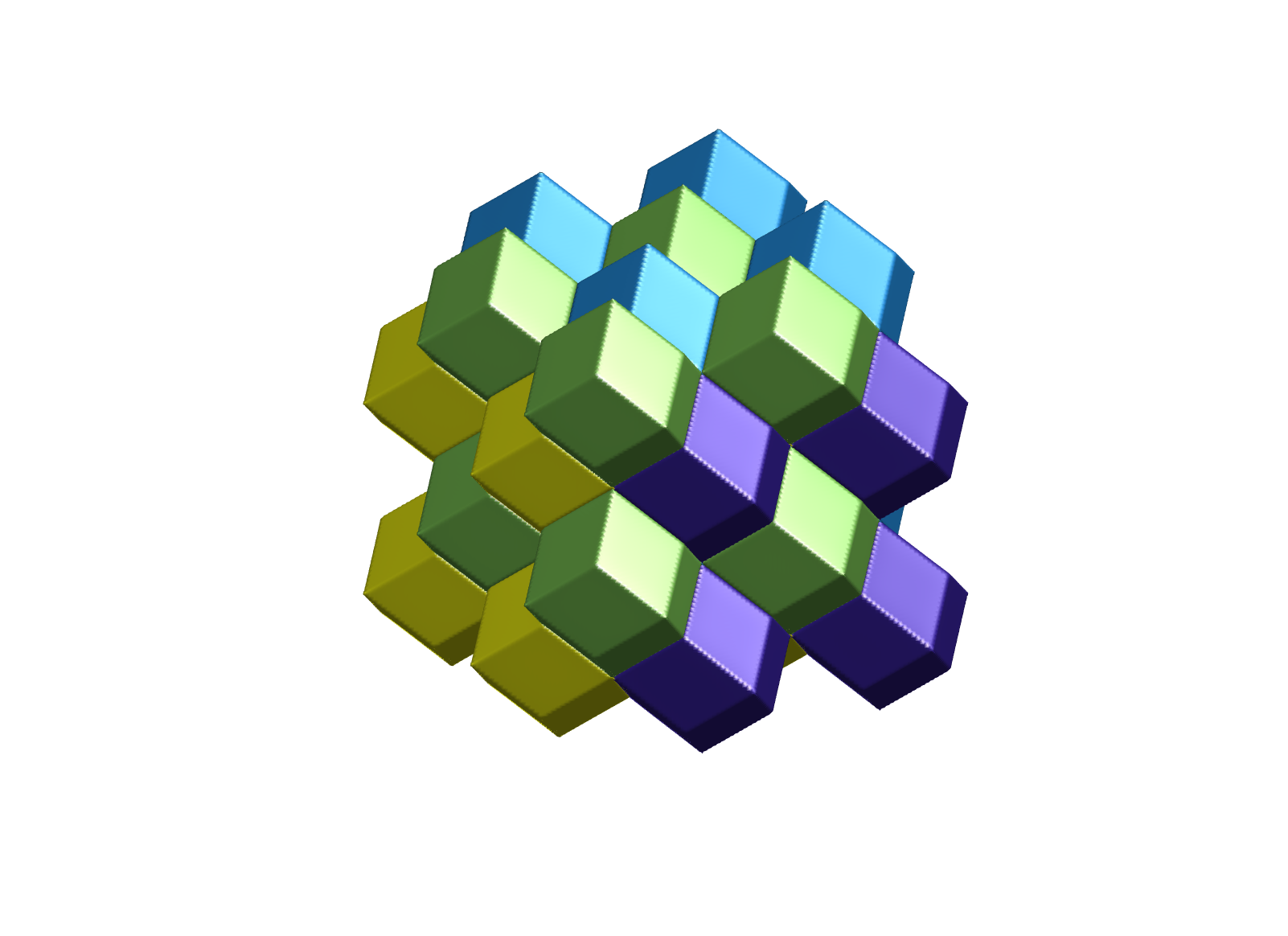}
\caption{A $k=4$ Dirichlet partition of the periodic cube, $[-1,1]^3$ consisting of 
rhombic dodecahedra  (left). On the right, we periodically extend the obtained partition to show how the rhombic dodecahedra fit together.  In this experiment, the cube is discretized by $128^3$ grid points and  $\tau = 0.125$. The CPU time for this experiment is $112$ seconds. } 
\label{fig:3d_2}
\end{figure}

For $k=4$ and  initialization using a random tessellation, we obtain a partition of the cube by four identical rhombic dodecahedron structures which is displayed in Figure~\ref{fig:3d_2}. In this experiment, the cube is discretized by $128^3$ uniform grid points and $\tau = 0.125$. The CPU time for this experiment is $112$ seconds.

\begin{figure}[ht]
\centering
\includegraphics[scale=0.25,clip,trim= 8cm 5cm 8cm 5cm]{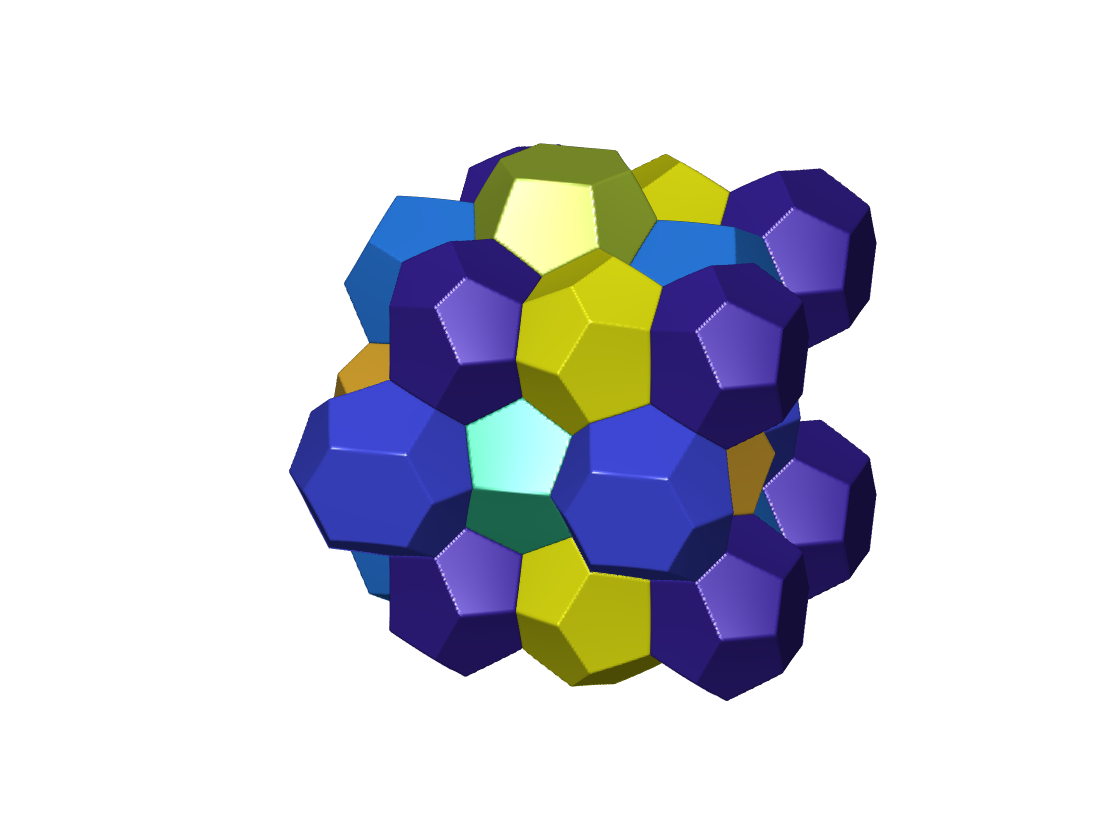}
\includegraphics[scale=0.2,clip,trim= 8cm 2cm 6cm 2cm]{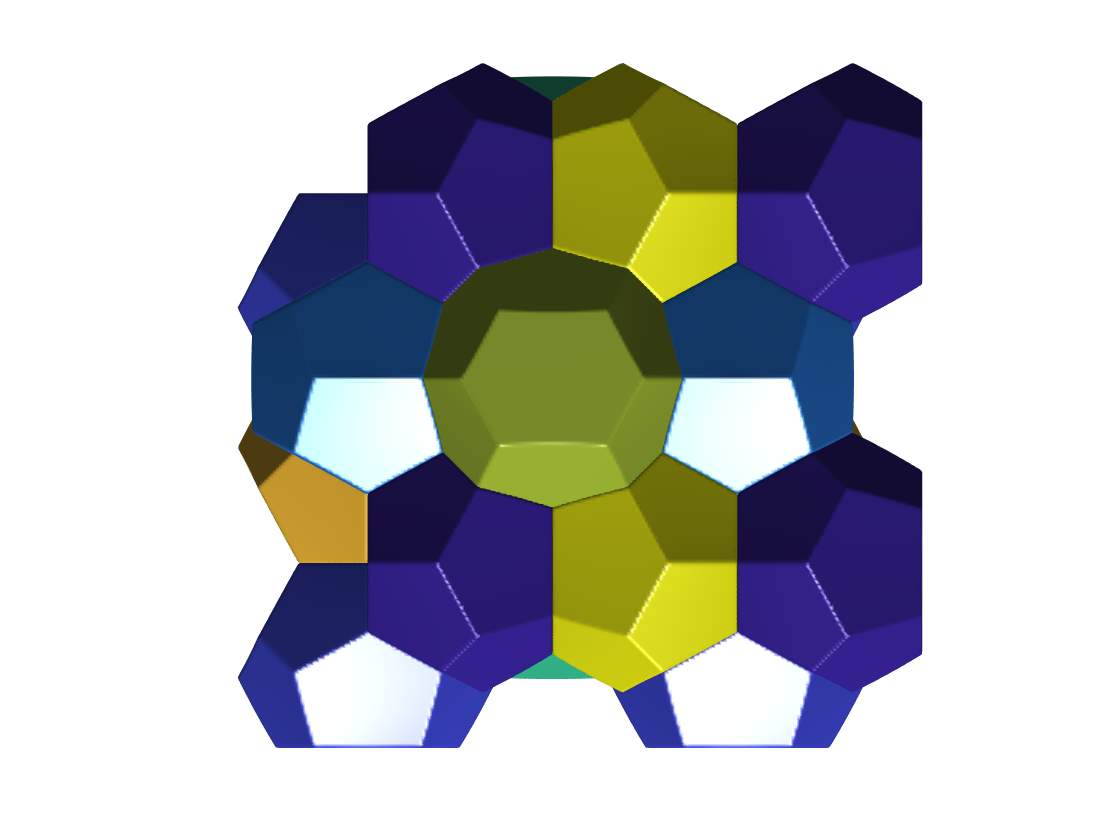}
\includegraphics[scale=0.2,clip,trim= 8cm 2cm 6cm 2cm]{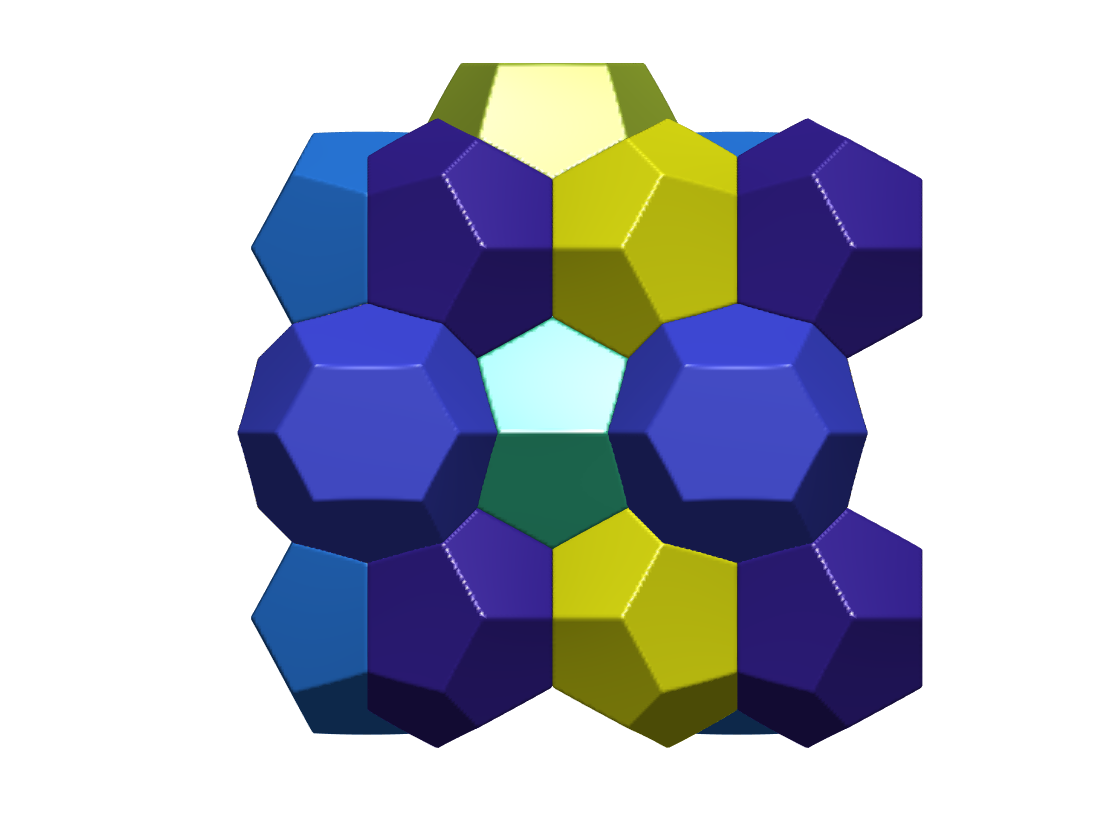}
\includegraphics[scale=0.2,clip,trim= 8cm 2cm 6cm 2cm]{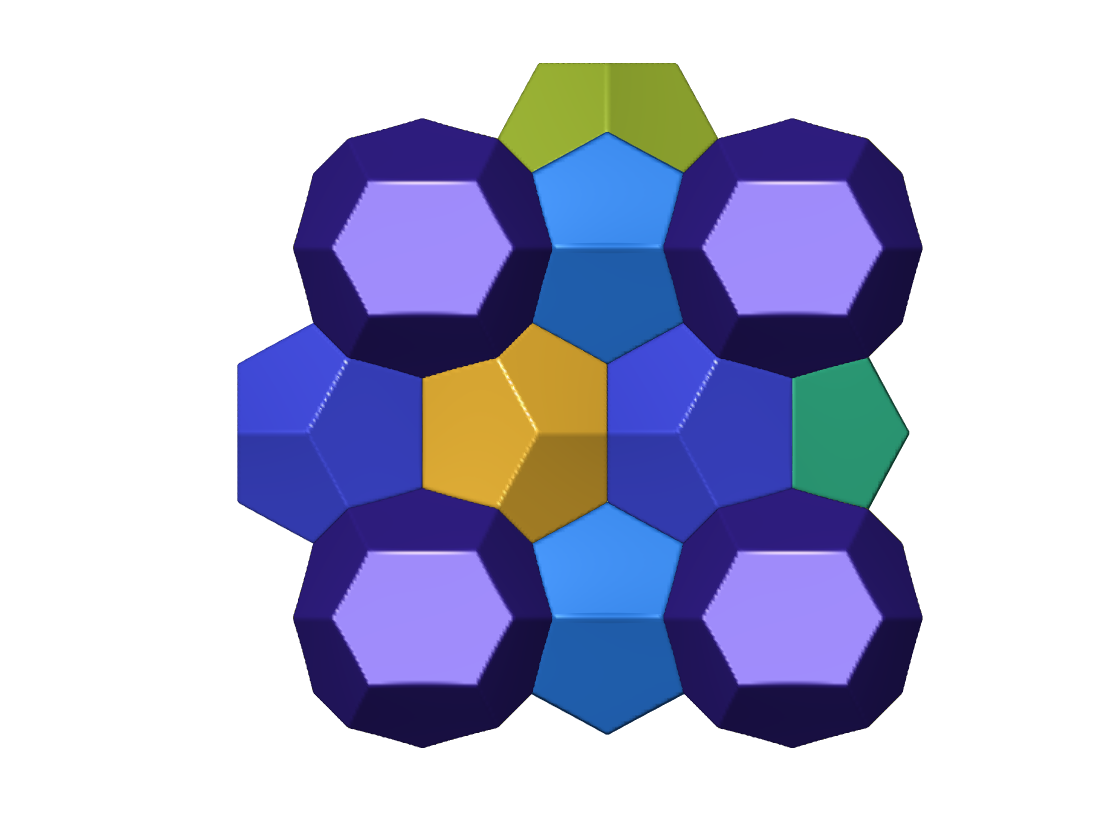}
\caption{A $k=8$ Dirichlet partition of the periodic cube, $[-1,1]^3$, which is similar to the 
Weaire-Phelan structure. The different panels show 
a 3d view (top left), 
a vertical view (top right),
a front view (bottom left), 
and a side view (bottom right). 
There are $6$ type--one Weaire-Phelan structures and $2$  type--two Weaire-€"Phelan structures in the partition; see Figures \ref{fig:3d_3_1} and \ref{fig:3d_3_2} for plots of these structures. In each panel, we have extended the partition periodically, so that it is easier to see how the structures fit together.  In this experiment, the cube is discretized by $128^3$ uniform grid points and $\tau=0.0625$. The CPU time for this experiment is $ 1200$ seconds.}
\label{fig:3d_3_0}
\end{figure}

\begin{figure}[ht]
\centering
\includegraphics[scale=0.18,clip,trim= 12cm 5cm 10cm 5cm]{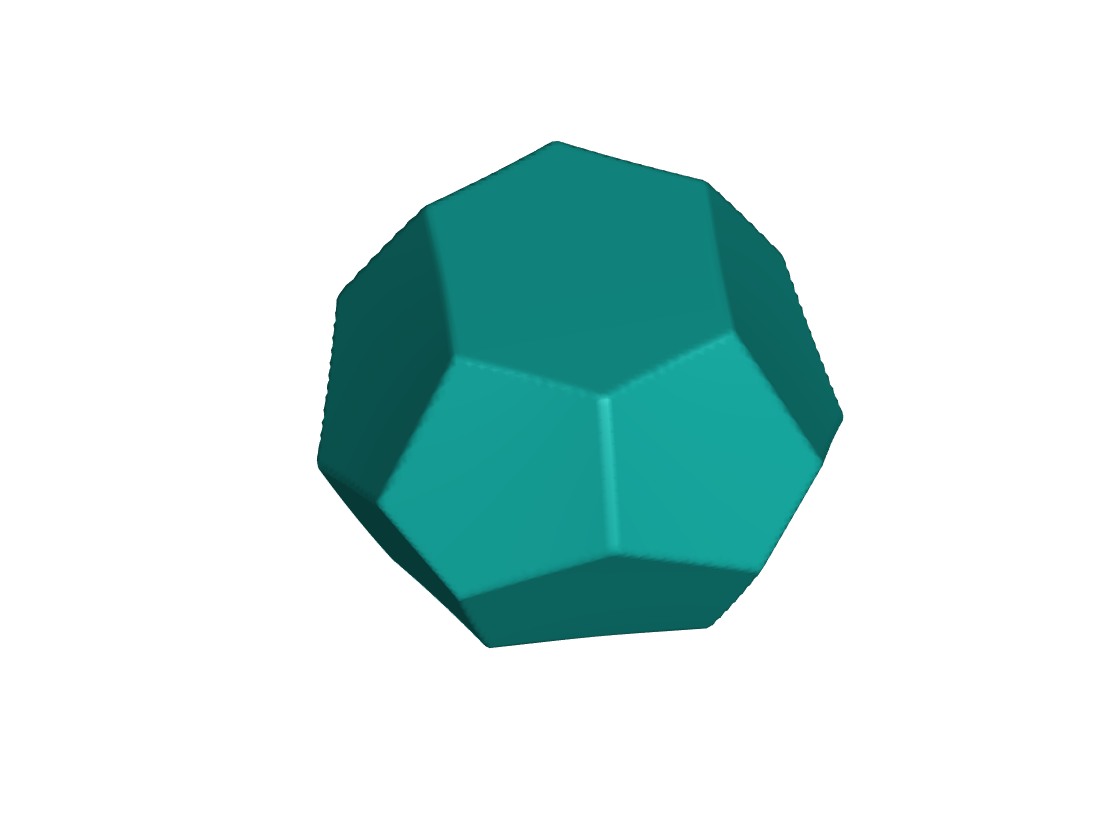}
\includegraphics[scale=0.14,clip,trim= 6cm 3cm 6cm 1cm]{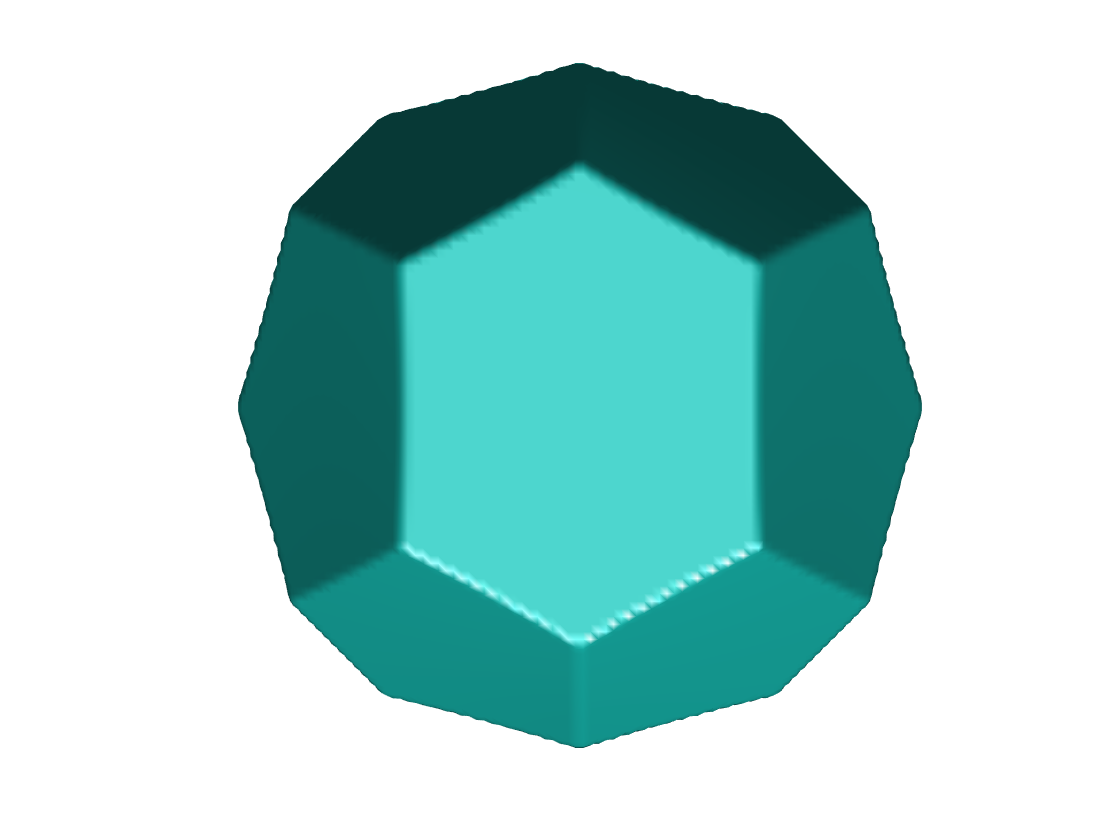}
\includegraphics[scale=0.12,clip,trim=3cm 3cm 3cm 1cm]{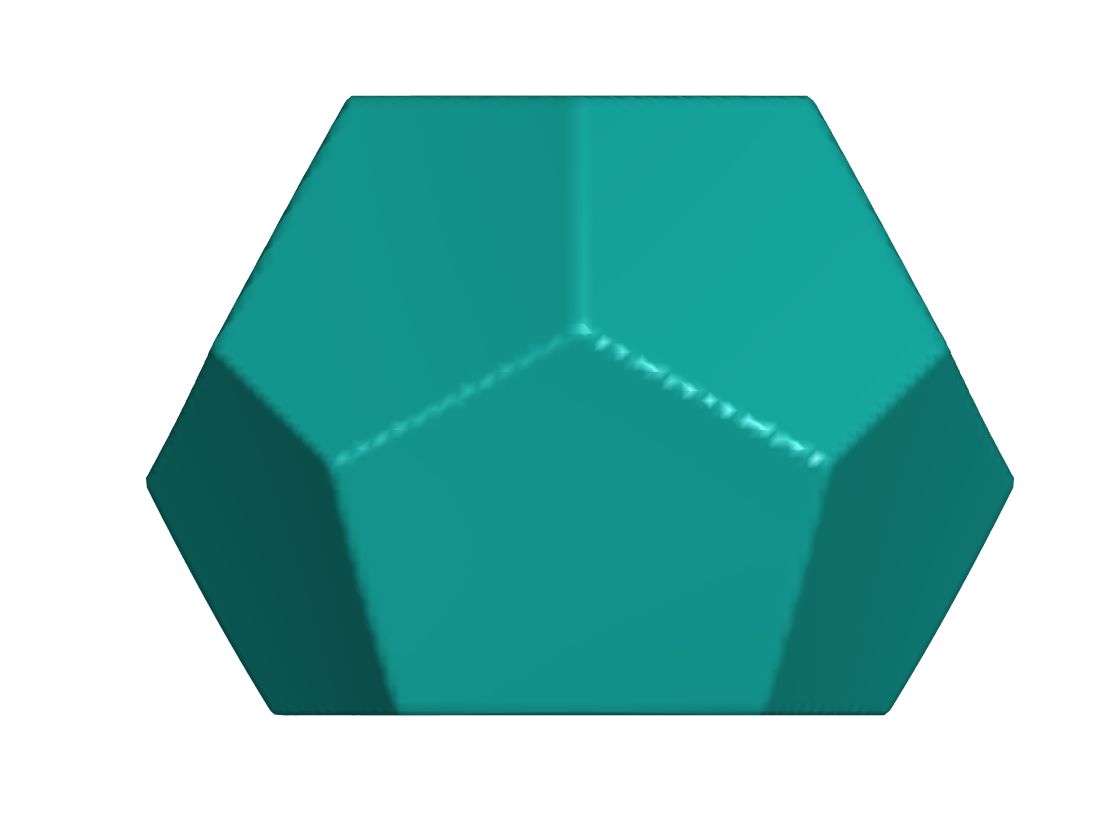}
\caption{ {\bf (left)} A type--one Weaire-€"Phelan structure, {\bf (center)} a vertical view, and {\bf (right)} a front view. The side view is same as the front view.} 
\label{fig:3d_3_1}
\end{figure}

\begin{figure}[ht]
\centering
\includegraphics[scale=0.18,clip,trim= 12cm 5cm 10cm 5cm]{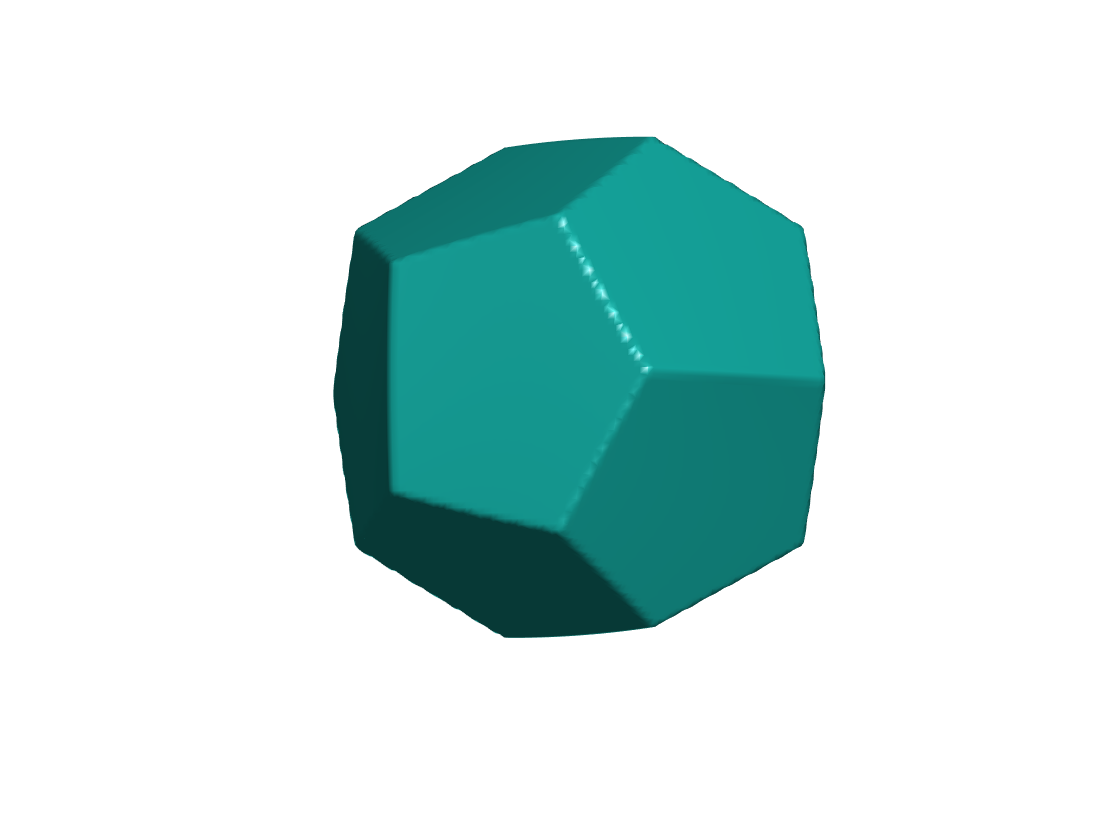}
\includegraphics[scale=0.15,clip,trim= 6cm 1cm 10cm 1 cm]{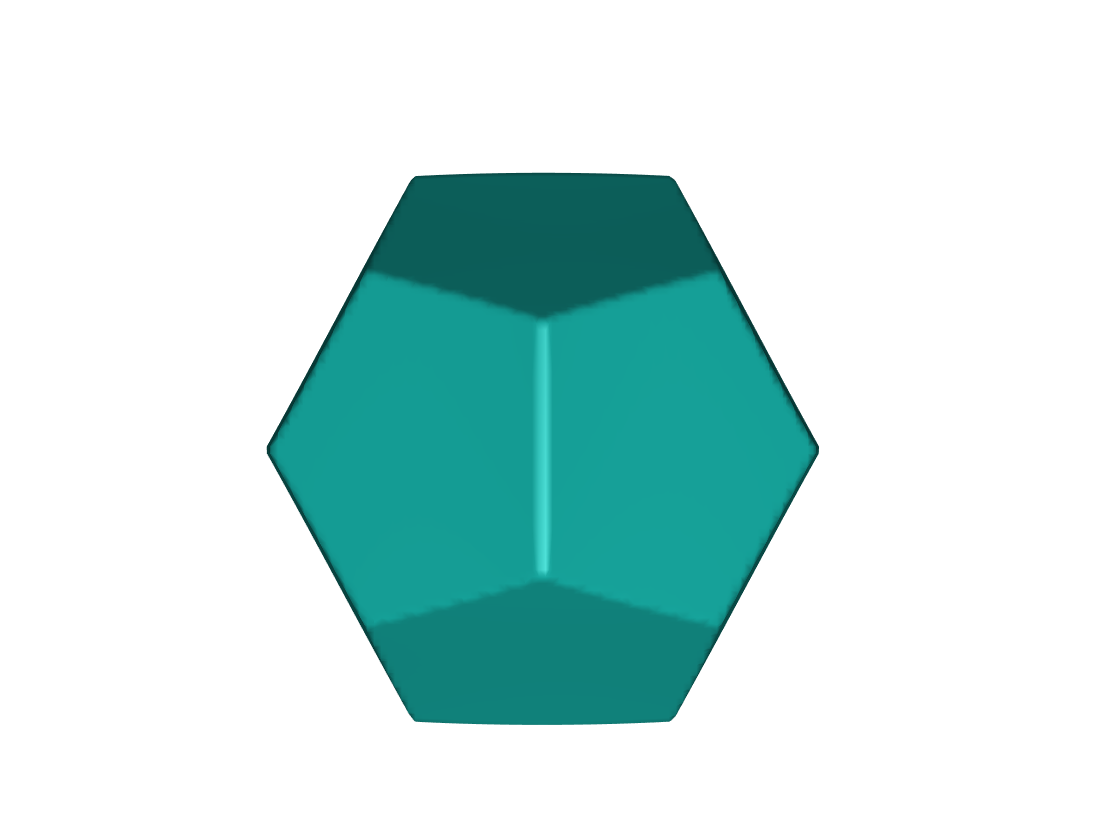}
\includegraphics[scale=0.15,clip,trim= 6cm 3cm 6cm 1cm]{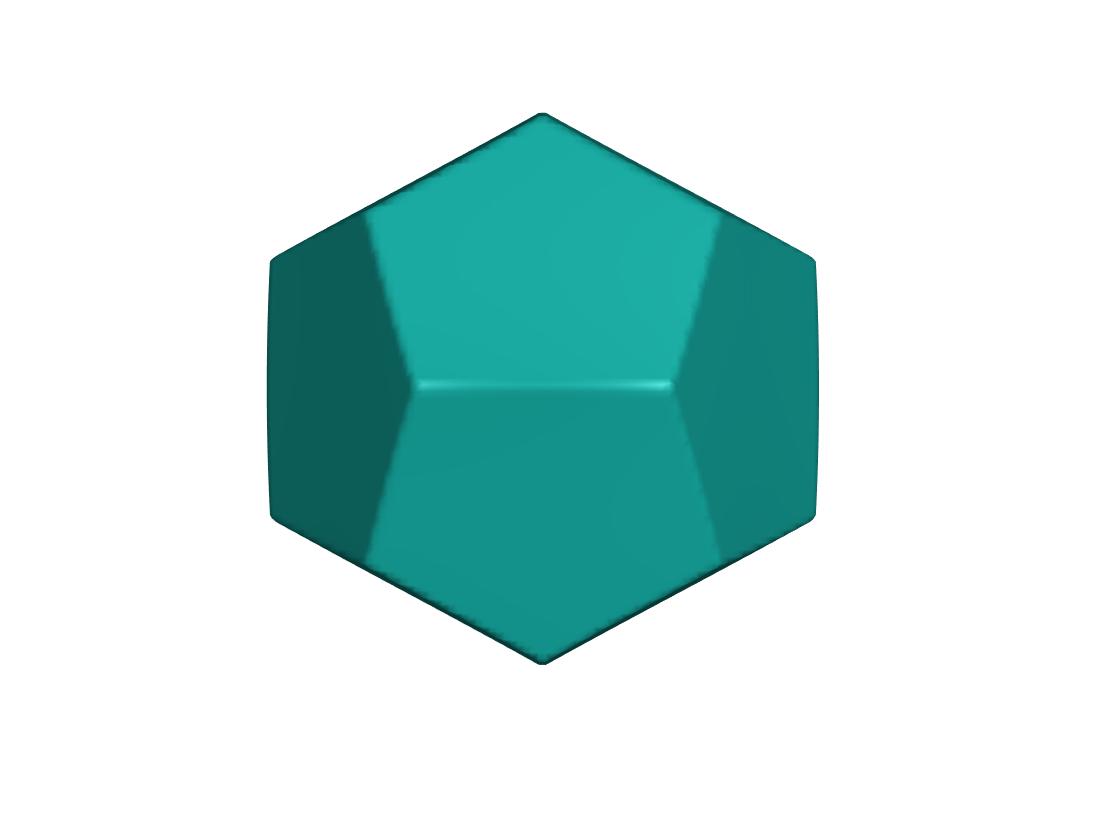}
\caption{ {\bf (left)} A type--two Weaire-Phelan structure, {\bf (center)} a vertical view, and {\bf (right)} a front view. The side view is same as the vertical view.} 
\label{fig:3d_3_2}
\end{figure}

For $k=8$ and  initialization using a random tessellation, we obtain a partition of the cube that is similar to the Weaire-Phelan structure. Figure~\ref{fig:3d_3_0} displays different views of a periodic extension of the partition.  Figures~\ref{fig:3d_3_1} and \ref{fig:3d_3_2} display different views of the first and second type Weaire-Phelan structures. In this experiment, the cube is discretized by $128^3$ uniform grid points and $\tau = 0.0625$. The CPU time for this experiment is $1200$ seconds. A rougher, but similar result can also be obtained by discretizing the cube with $64^3$ uniform grid points and using $\tau=0.0625$ in $81$ seconds. In the numerical experiments, our algorithm occasionally converged to other local minimizers. However, our experiments indicate that the algorithm usually converges to the Weaire-Phelan structure,  implying that the basin of attraction for this structure is larger.

\begin{figure}[ht]
\centering
\includegraphics[scale=0.18,clip,trim= 9cm 4cm 9cm 2cm]{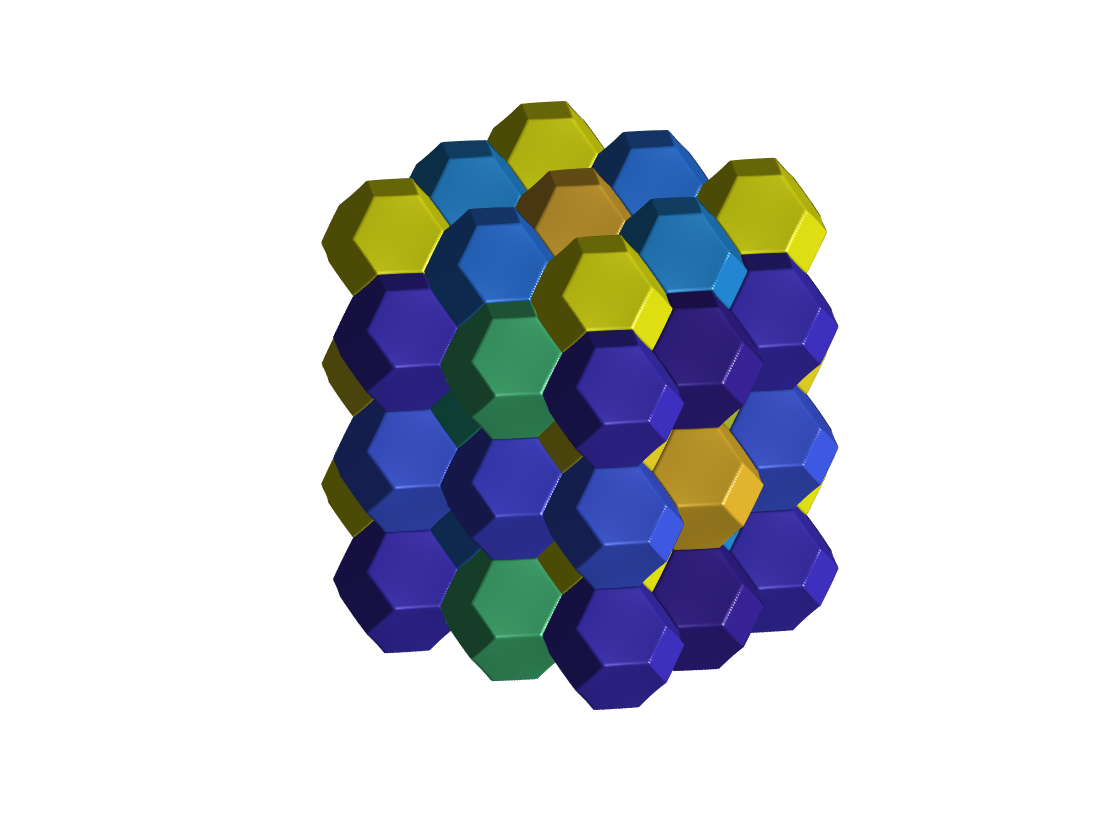}
\includegraphics[scale=0.14,clip,trim= 8cm 2cm 6cm 2cm]{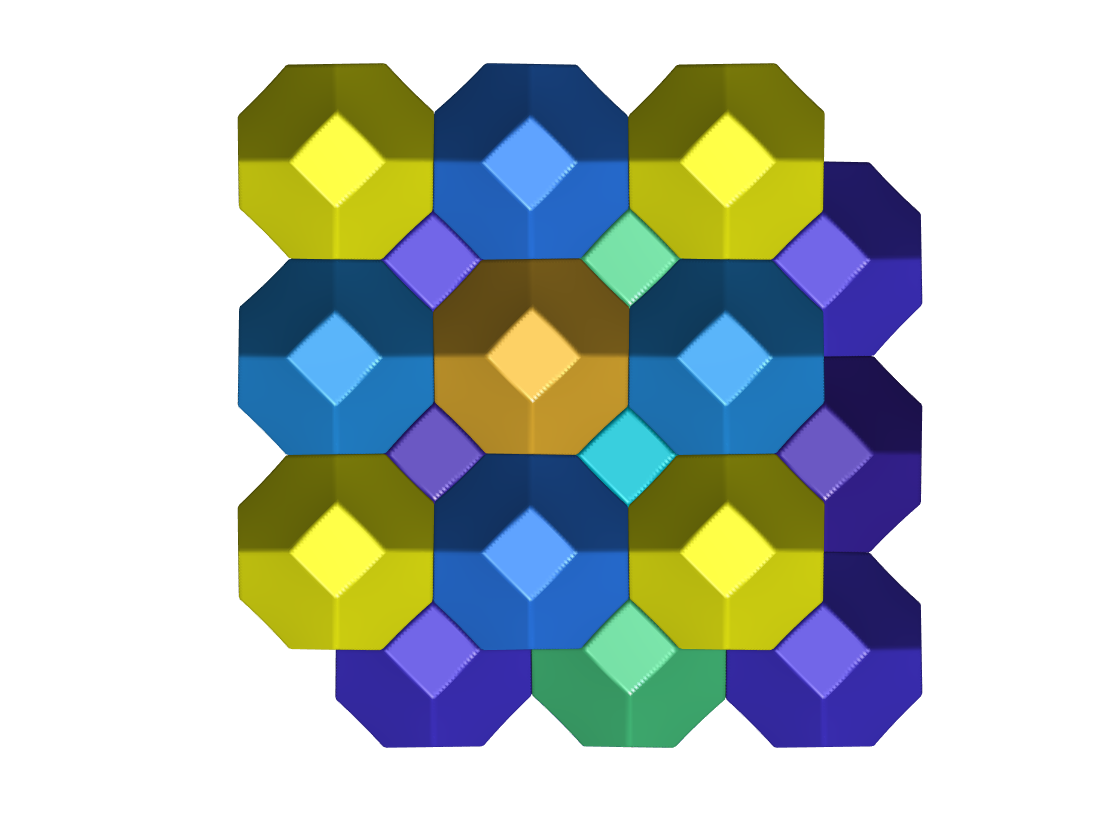}
\includegraphics[scale=0.14,clip,trim= 8cm 2cm 6cm 2cm]{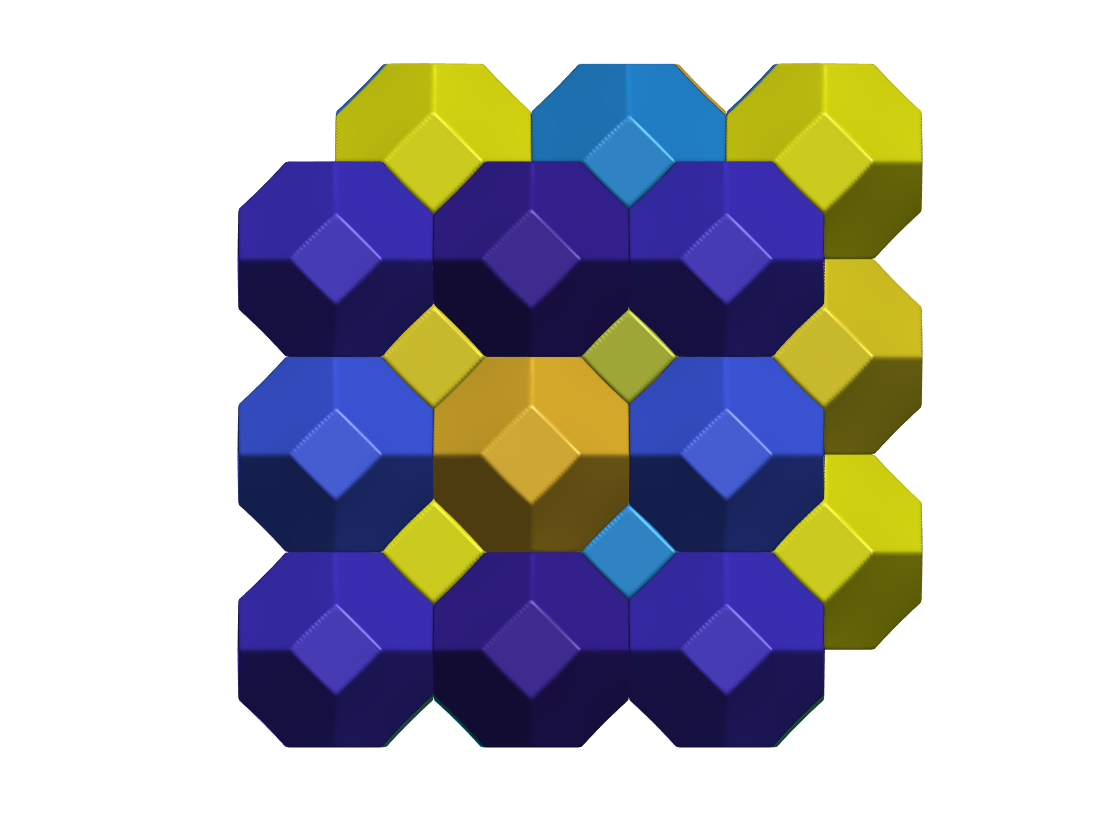}
\caption{ 
{\bf (left)}  A $k=12$ Dirichlet partition of the periodic cube, $[-1,1]^3$, by equal truncated octahedra, similar to Kelvin's structure. The partition has been periodically extended. 
{\bf (center)} A vertical view. 
{\bf (right)} A side view. 
The front view is same as the vertical view. 
In this experiment, the cube is discretized by $128^3$ uniform grid points and $\tau=0.0625$. 
The CPU time for this experiment is $3556$ seconds.}
\label{fig:3d_4_0}
\end{figure}

For $k=16$ and initialization using a random tessellation, we obtain a partition of the cube that is a packing of truncated octahedra, similar to the structure Lord Kelvin studied. Figure~\ref{fig:3d_4_0} displays different views of a periodic extension of this partition. 
In this experiment, the cube is discretized by $128^3$ uniform grid points and $\tau = 0.0625$. The CPU time for this experiment is $ 3556$ seconds.

In Table~\ref{tab:3d}, we also tabulate the values of  $\tilde E$ in \eqref{e:aEnergy}, the CPU time, and the $\tau$ used for different values of $k$.

\begin{table}[ht]
\centering
\caption{Values of  $\tilde E$ in \eqref{e:aEnergy}, the CPU time, and the $\tau$ used for different values of $k$.} \label{tab:3d}
\begin{tabular}{|c|c|c|c|c|c|c|c|c|c|c|c|c|}
\hline
$k$ & 2(left) & 2(right) & 4& 8& 16  \\
\hline 
$\tilde E$  & 3.43 & 3.61 & 3.07 & 2.68 & 2.47 \\
 \hline
 CPU time (s) & 26 & 3 & 112  & 1200 & 3556 \\
 \hline 
 $\tau$ & 0.25 & 0.25 & 0.125 & 0.0625 & 0.0625\\
 \hline    
\end{tabular}
\end{table}
 
\subsection{4d flat torus} \label{s:4d}
To our knowledge, neither partitions that minimize the total surface area or Dirichlet partitions in four dimensional space have been studied. In this section, we compute Dirichlet partitions using Algorithm~\ref{a:MBO} for the tesseract, $[-1,1]^4$, with periodic boundary conditions  and $k=4,8$.

For $k=4$ and initialization using a random tessellation, we obtain a constant extension of a rhombic dodecahedral honeycomb along the fourth direction. A rhombic dodecahedral honeycomb is plotted in Figure~\ref{fig:3d_2}; we do not include a figure of this extension. 

\begin{figure}
\centering
\includegraphics[scale=0.18,clip,trim= 3cm 1.5cm 3cm 1.2cm]{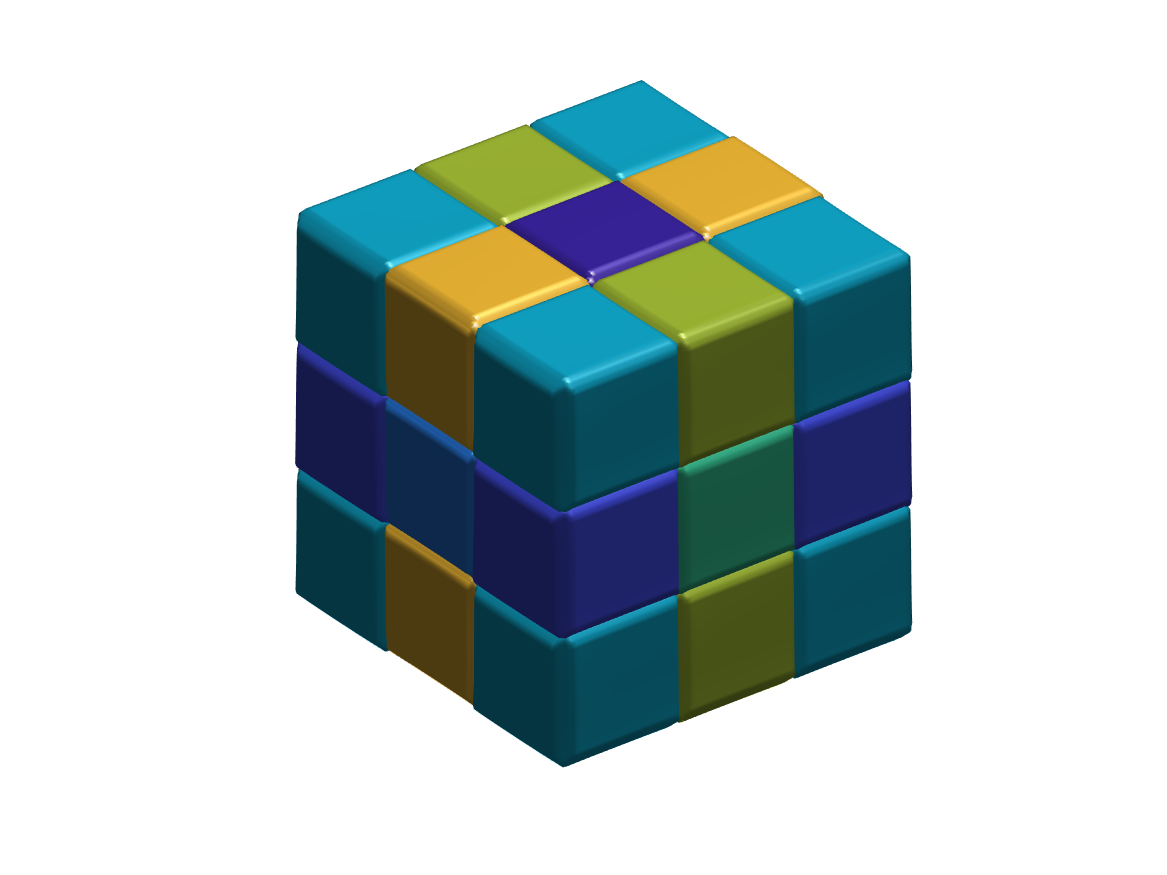}
\includegraphics[scale=0.18,clip,trim= 3cm 1.5cm 3cm 1.2cm]{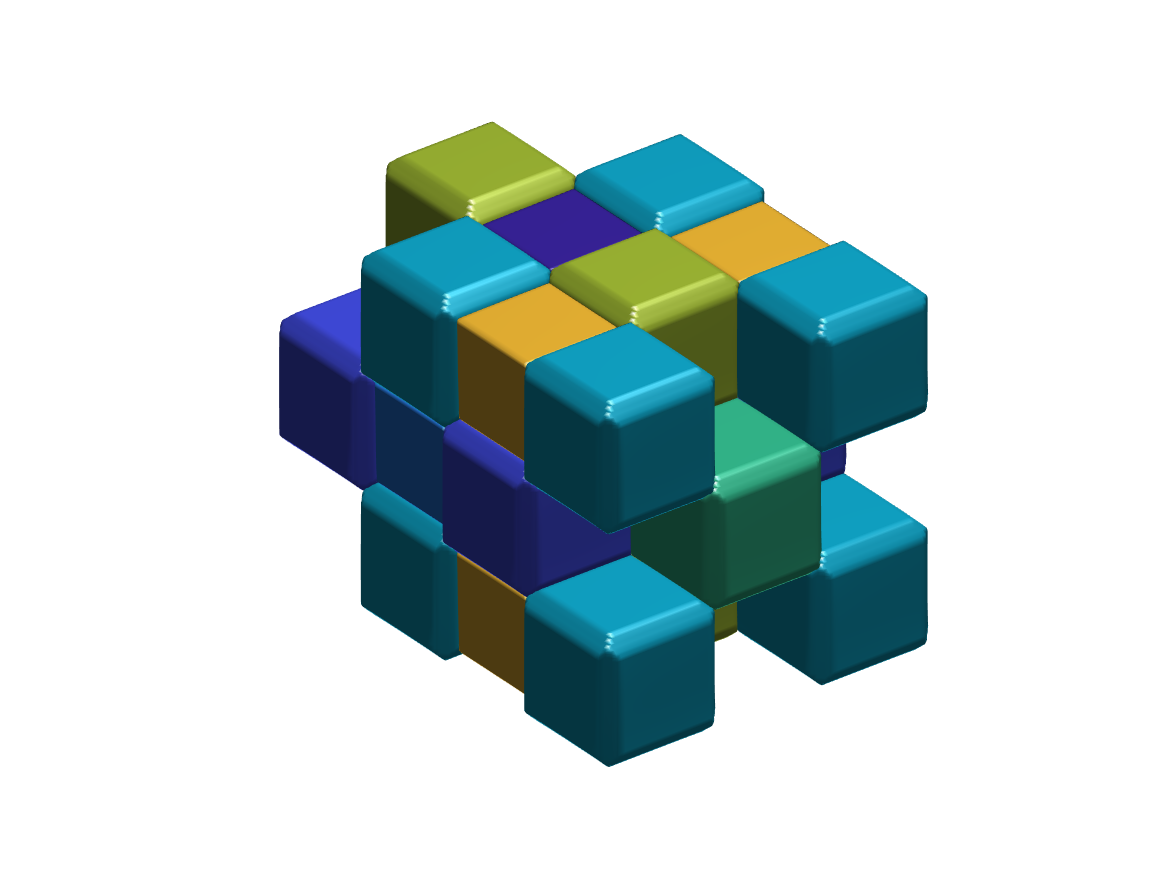}
\includegraphics[scale=0.18,clip,trim= 3cm 1.5cm 3cm 1.2cm]{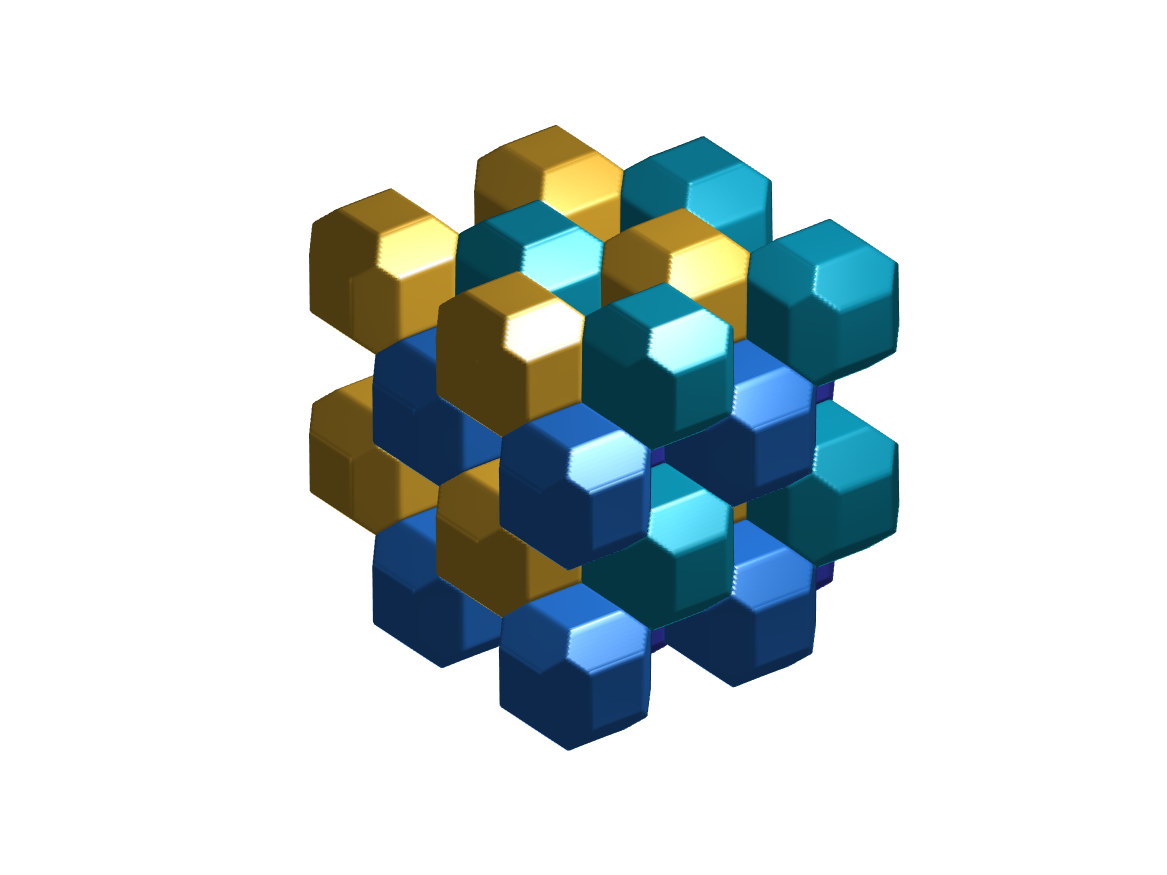}
\includegraphics[scale=0.18,clip,trim= 3cm 1.5cm 3cm 1.2cm]{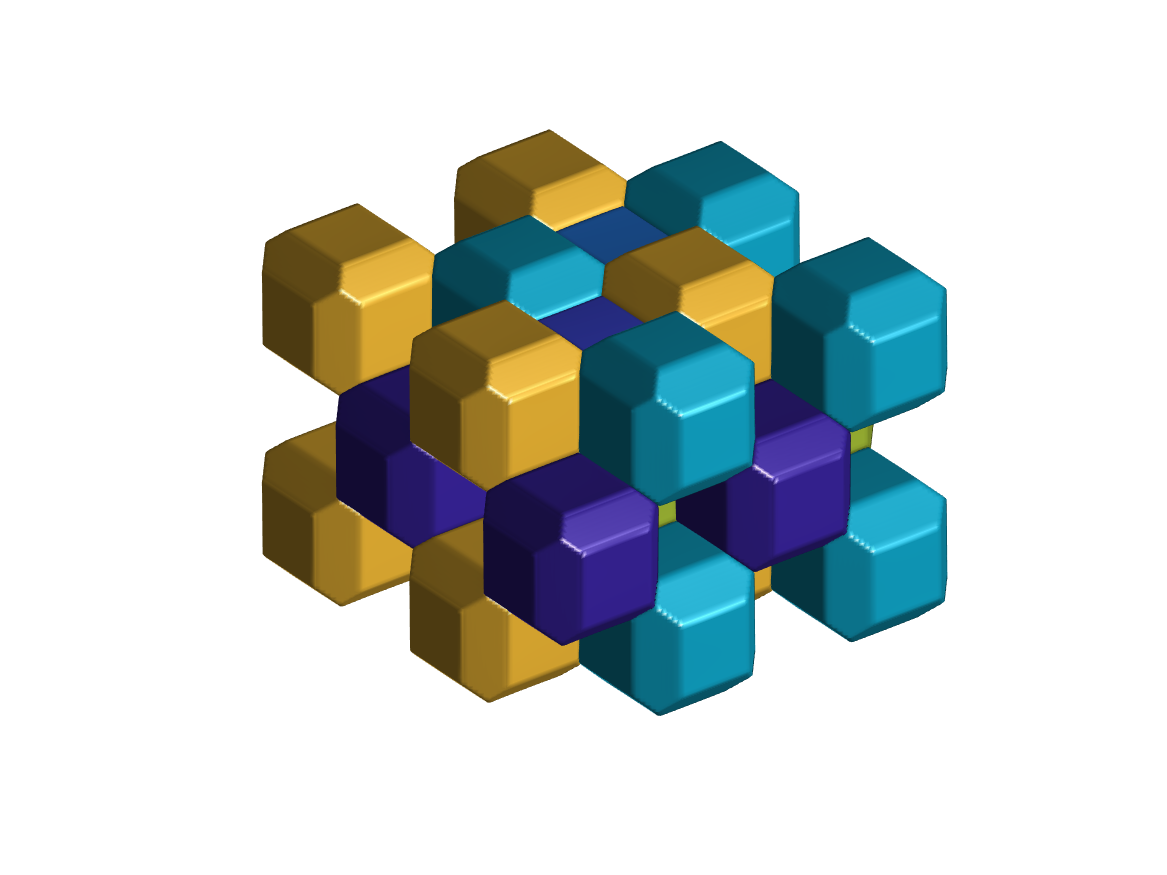}\\
\includegraphics[scale=0.18,clip,trim=3cm 1.5cm 3cm 1.2cm]{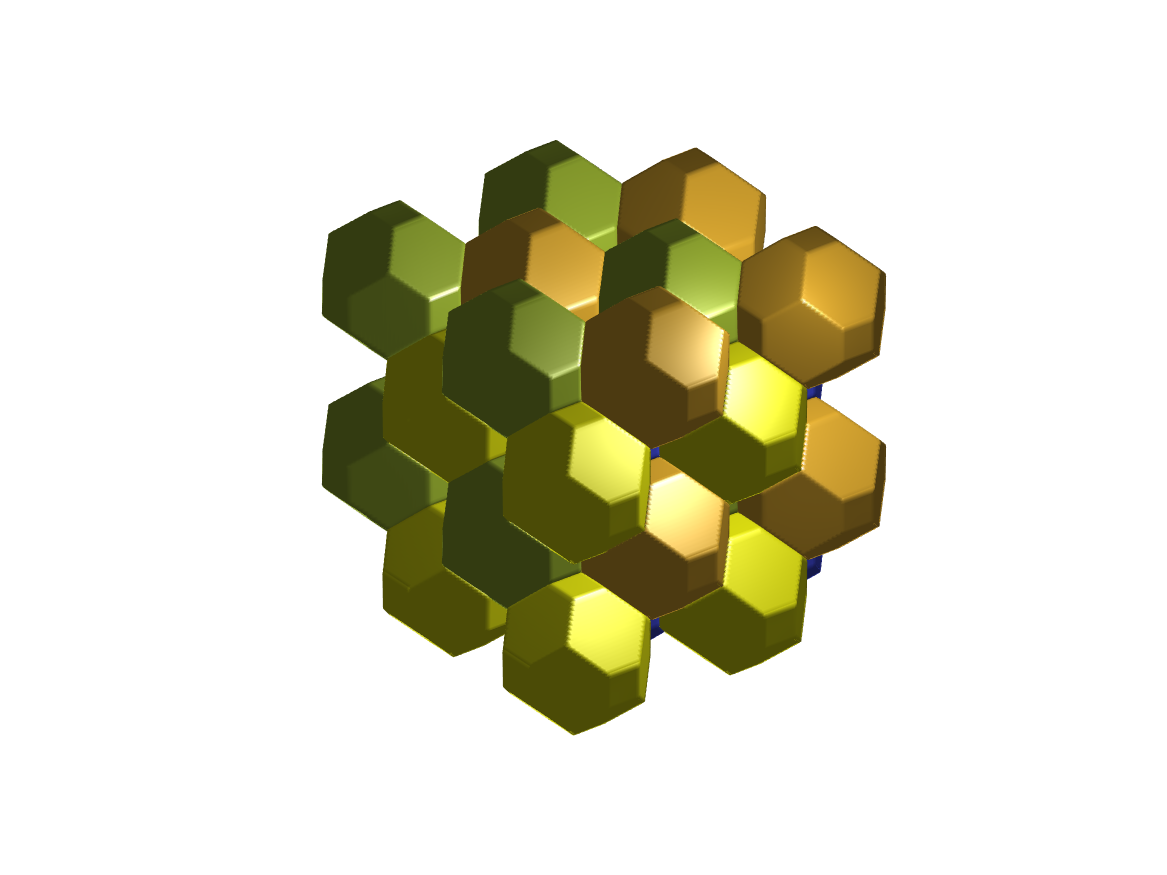}
\includegraphics[scale=0.18,clip,trim= 3cm 1.5cm 3cm 1.2cm]{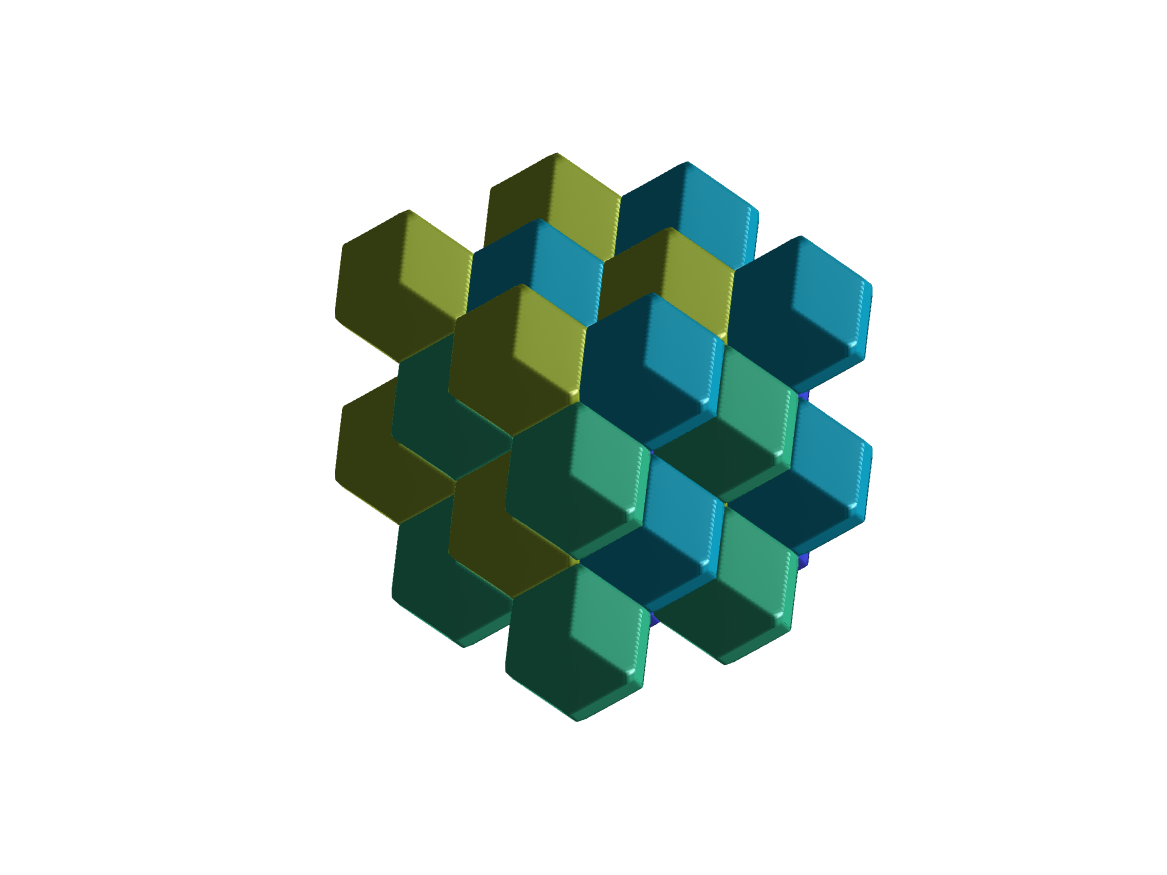}
\includegraphics[scale=0.18,clip,trim= 3cm 1.5cm 3cm 1.2cm]{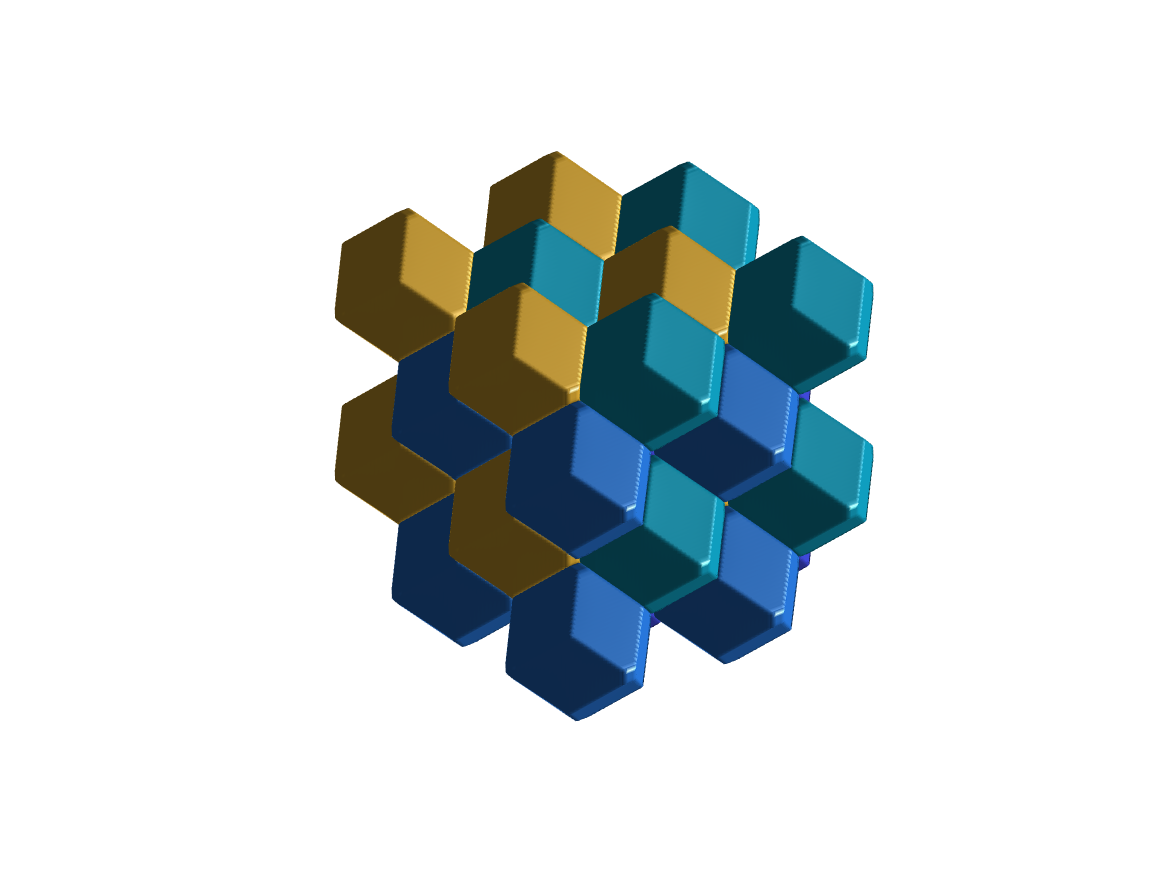}
\includegraphics[scale=0.18,clip,trim= 3cm 1.5cm 3cm 1.2cm]{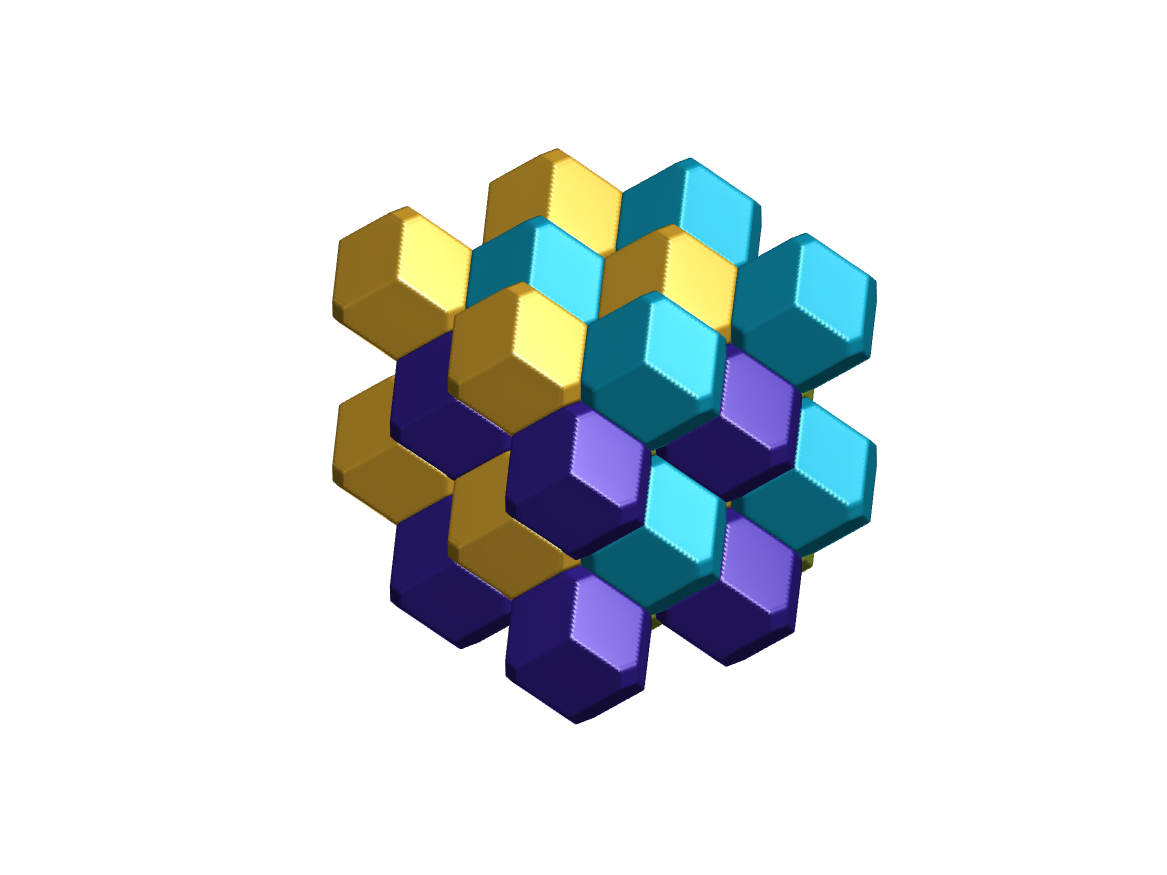}\\
\includegraphics[scale=0.18,clip,trim= 3cm 1.5cm 3cm 1.2cm]{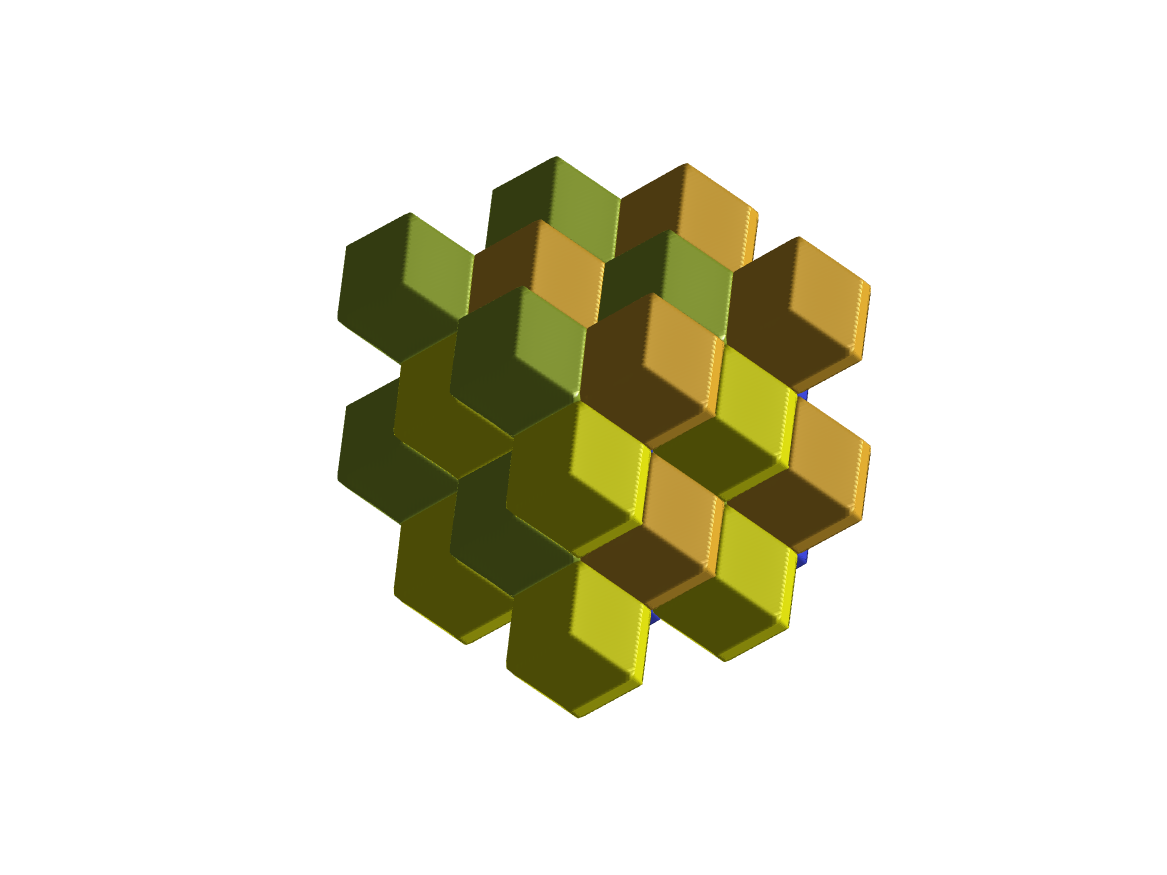}
\includegraphics[scale=0.18,clip,trim= 3cm 1.5cm 3cm 1.2cm]{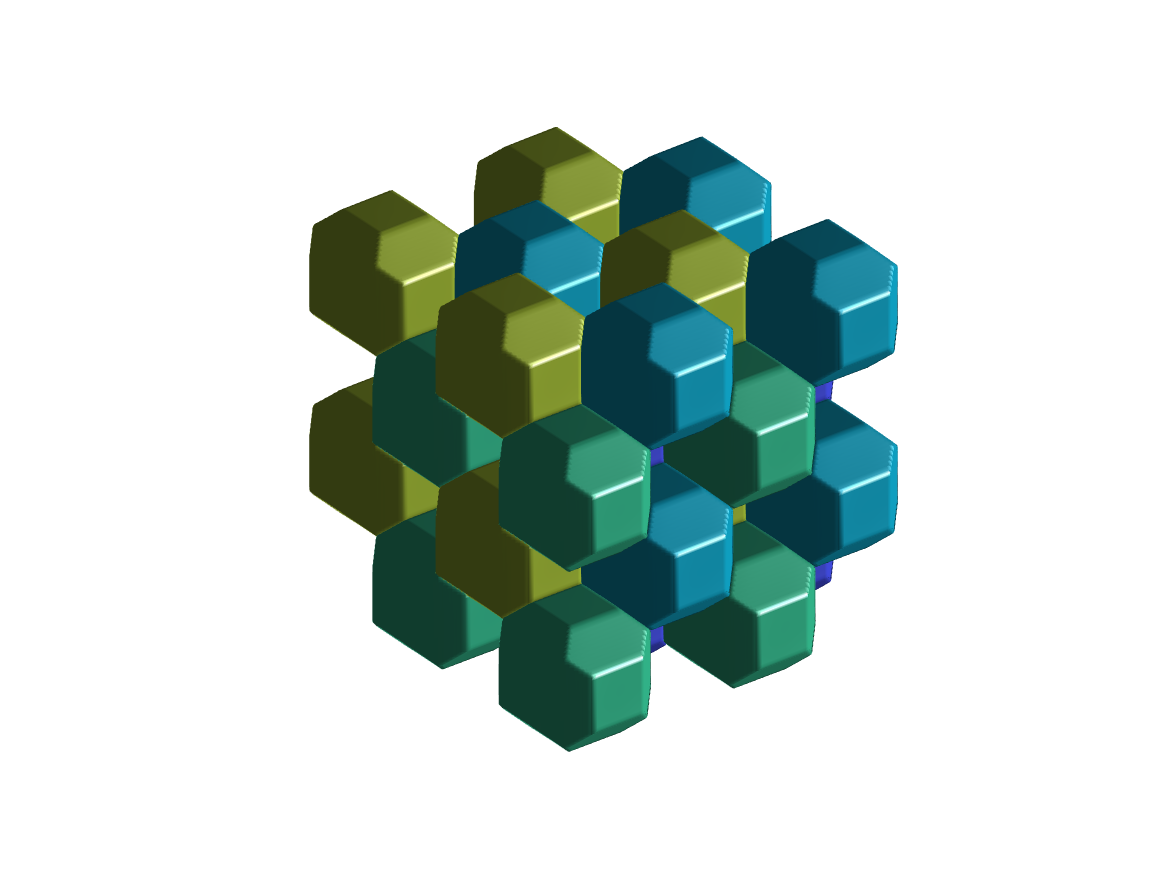}
\includegraphics[scale=0.18,clip,trim= 3cm 1.5cm 3cm 1.2cm]{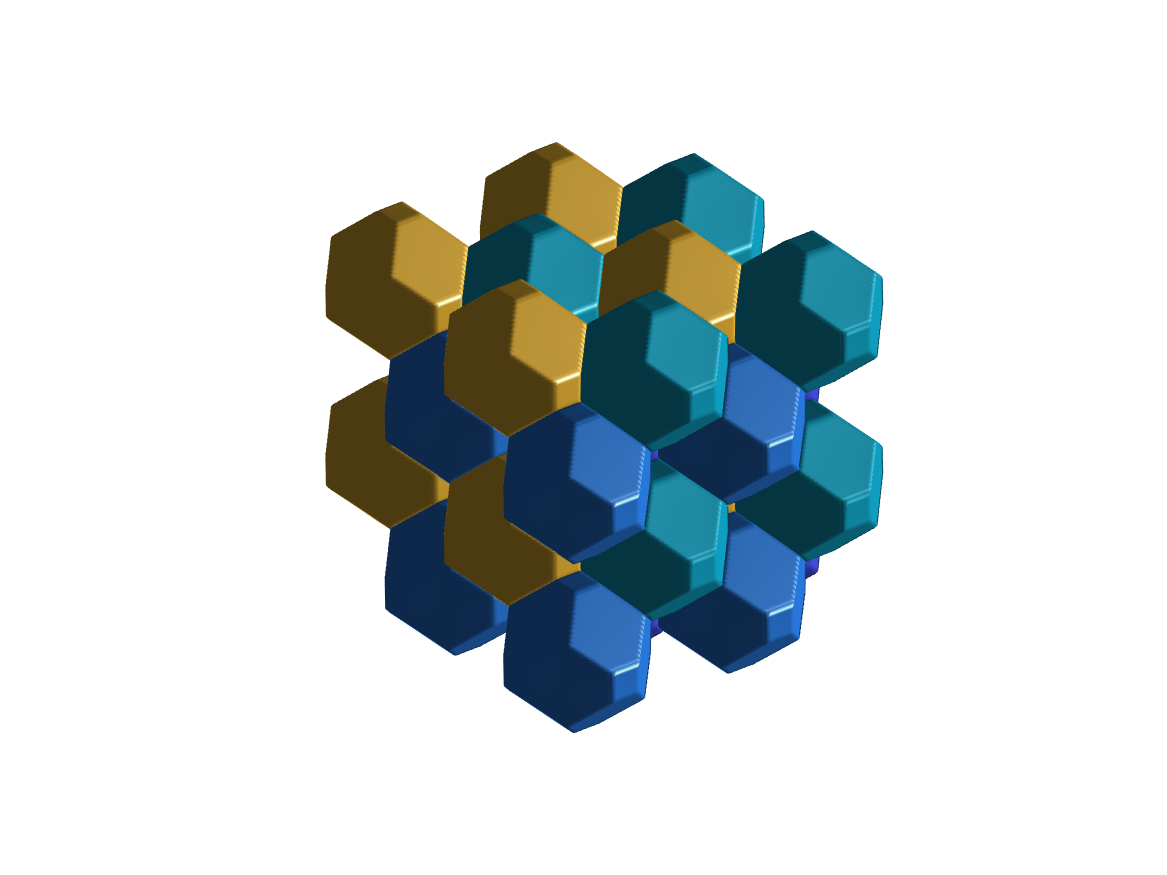}
\includegraphics[scale=0.18,clip,trim= 3cm 1.5cm 3cm 1.2cm]{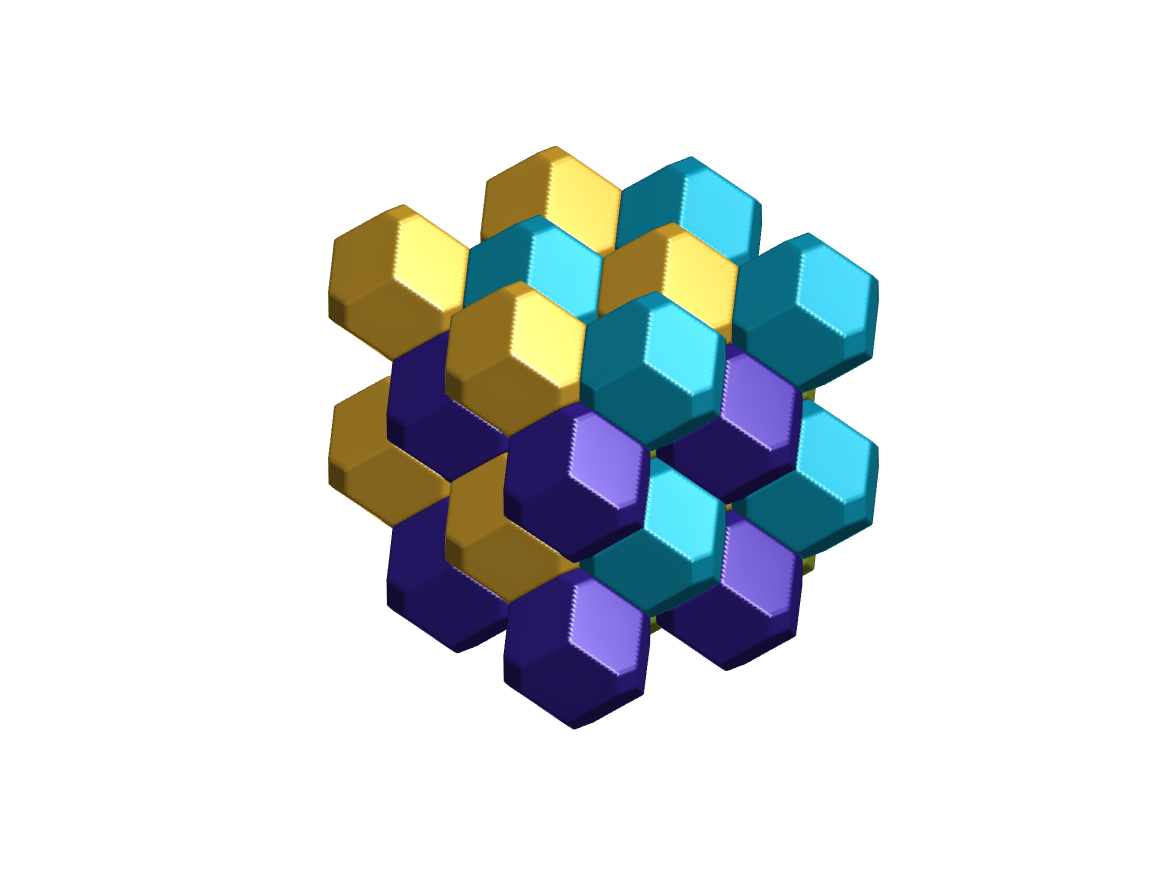}\\
\includegraphics[scale=0.18,clip,trim= 3cm 1.5cm 3cm 1.2cm]{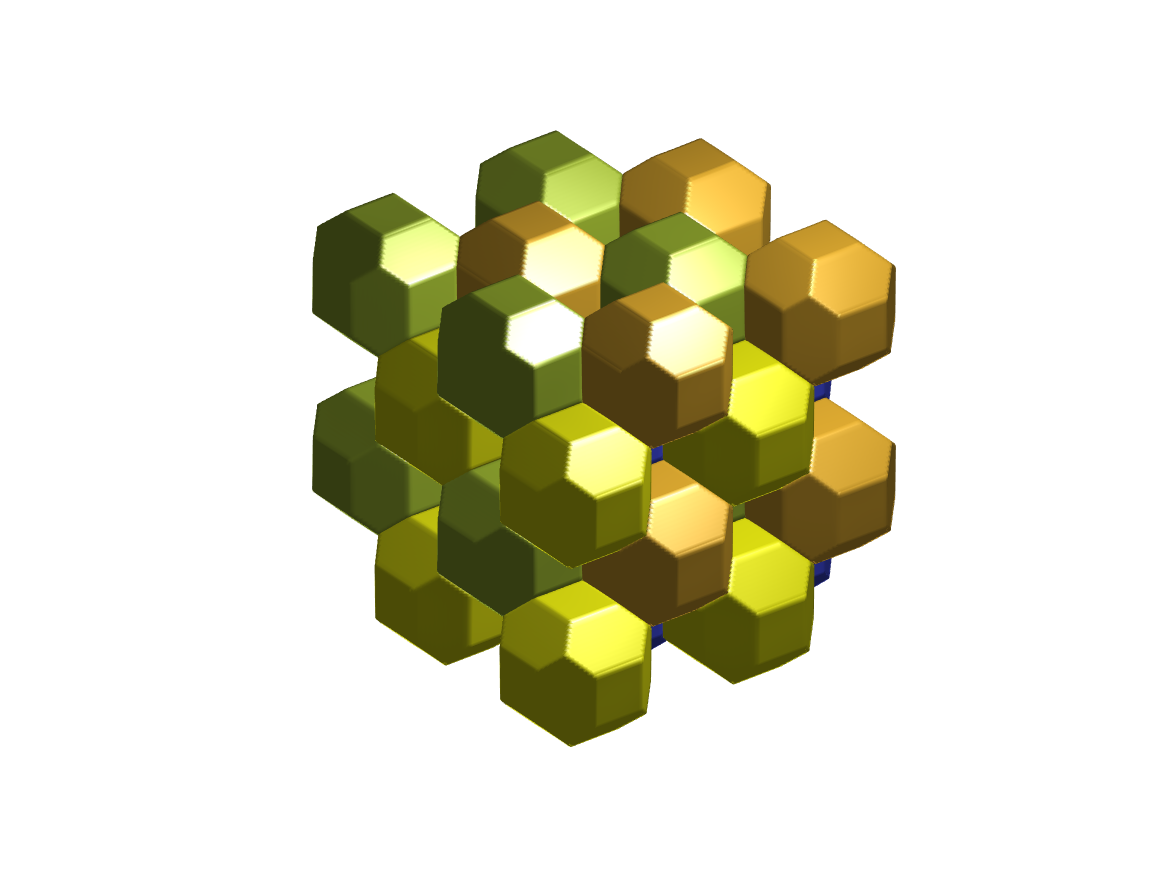}
\includegraphics[scale=0.18,clip,trim= 3cm 1.5cm 3cm 1.2cm]{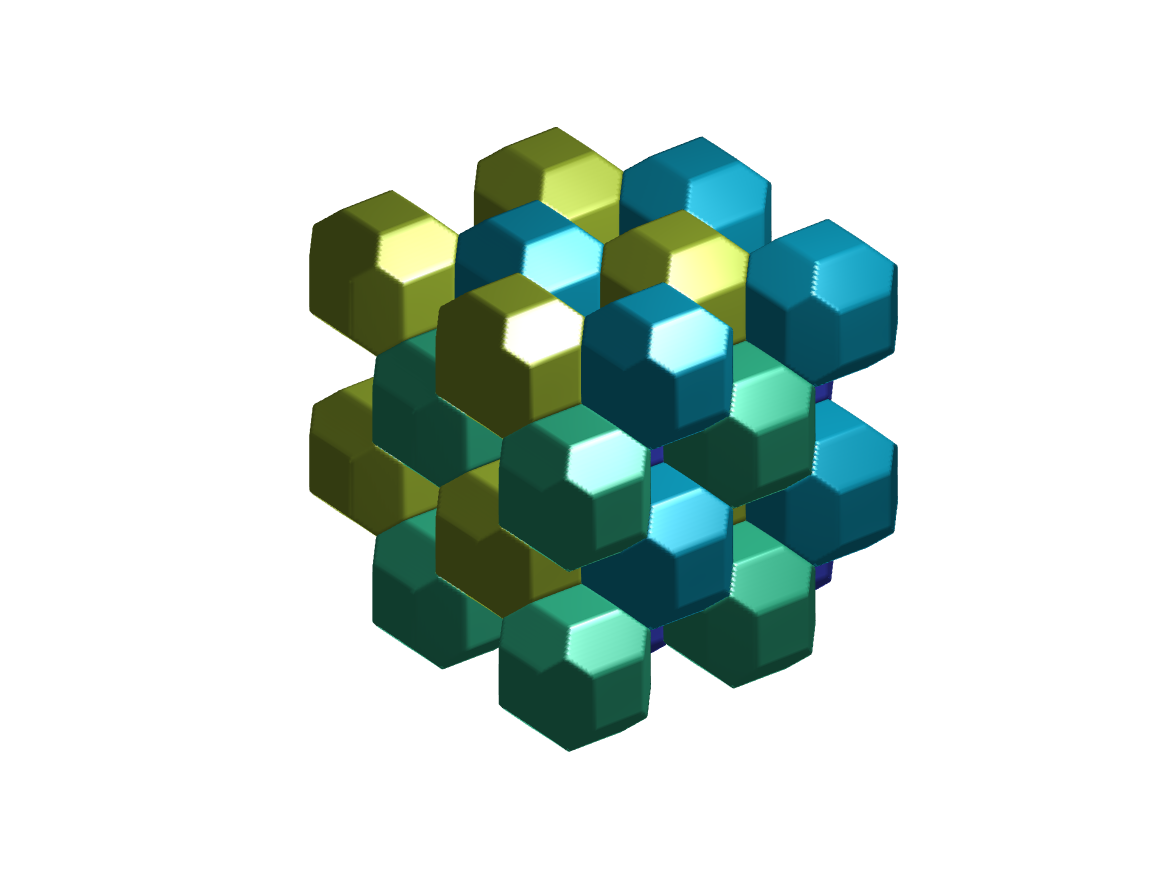}
\includegraphics[scale=0.18,clip,trim= 3cm 1.5cm 3cm 1.2cm]{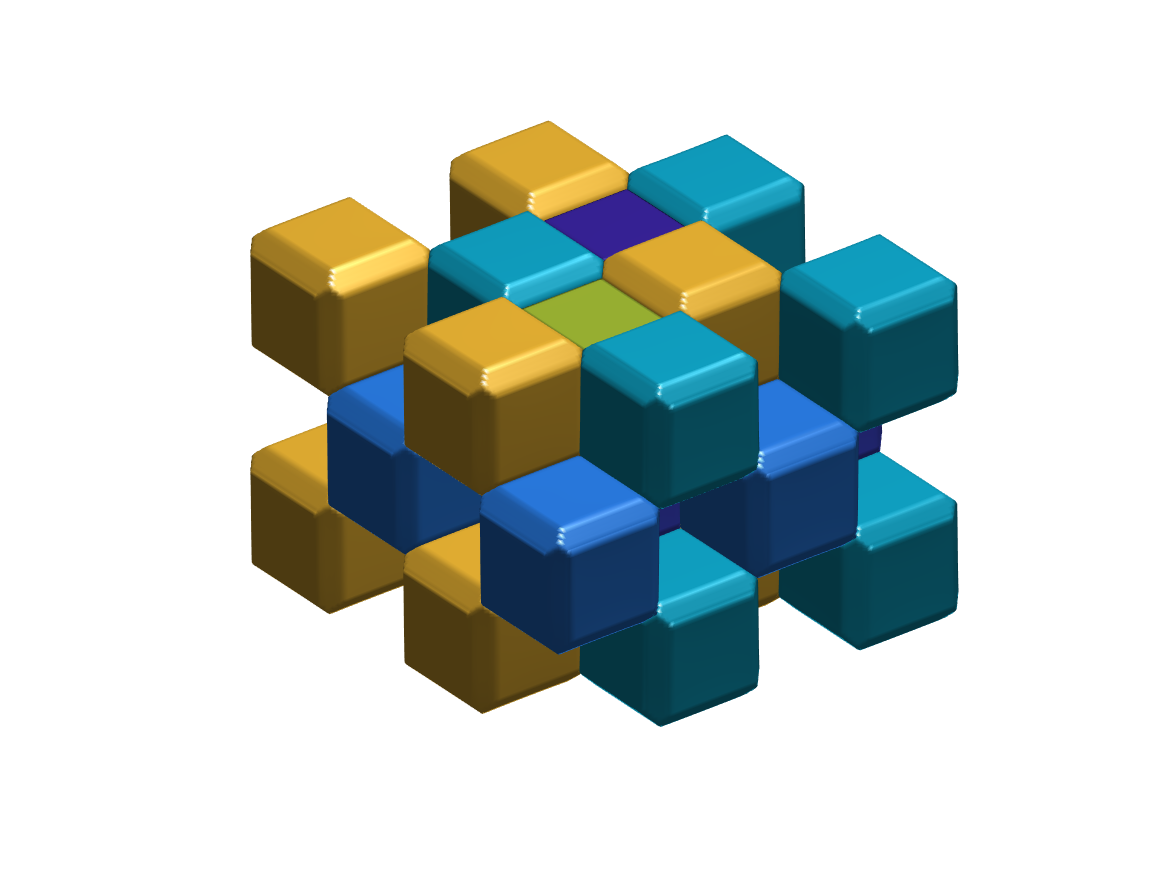}\includegraphics[scale=0.18,clip,trim= 3cm 1.5cm 3cm 1.2cm]{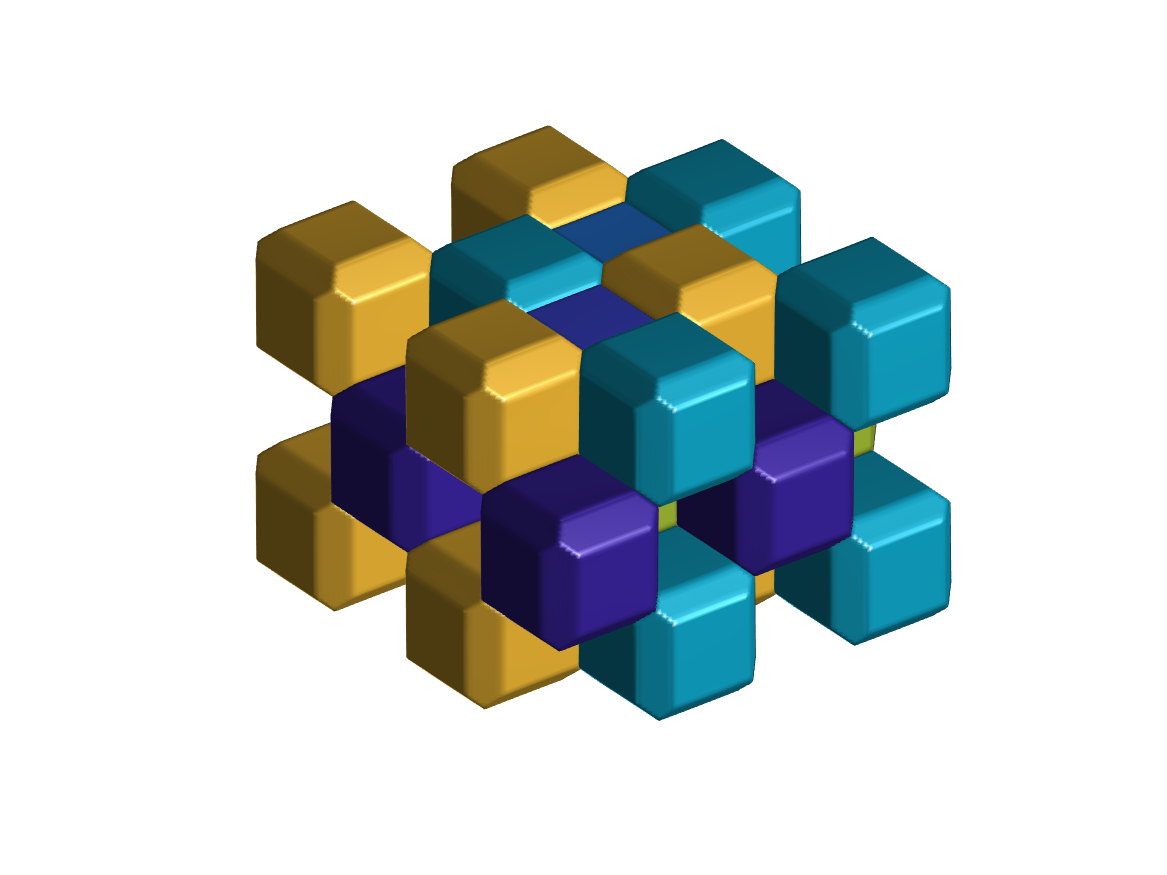}\\
\includegraphics[scale=0.18,clip,trim= 3cm 1.5cm 3cm 1.2cm]{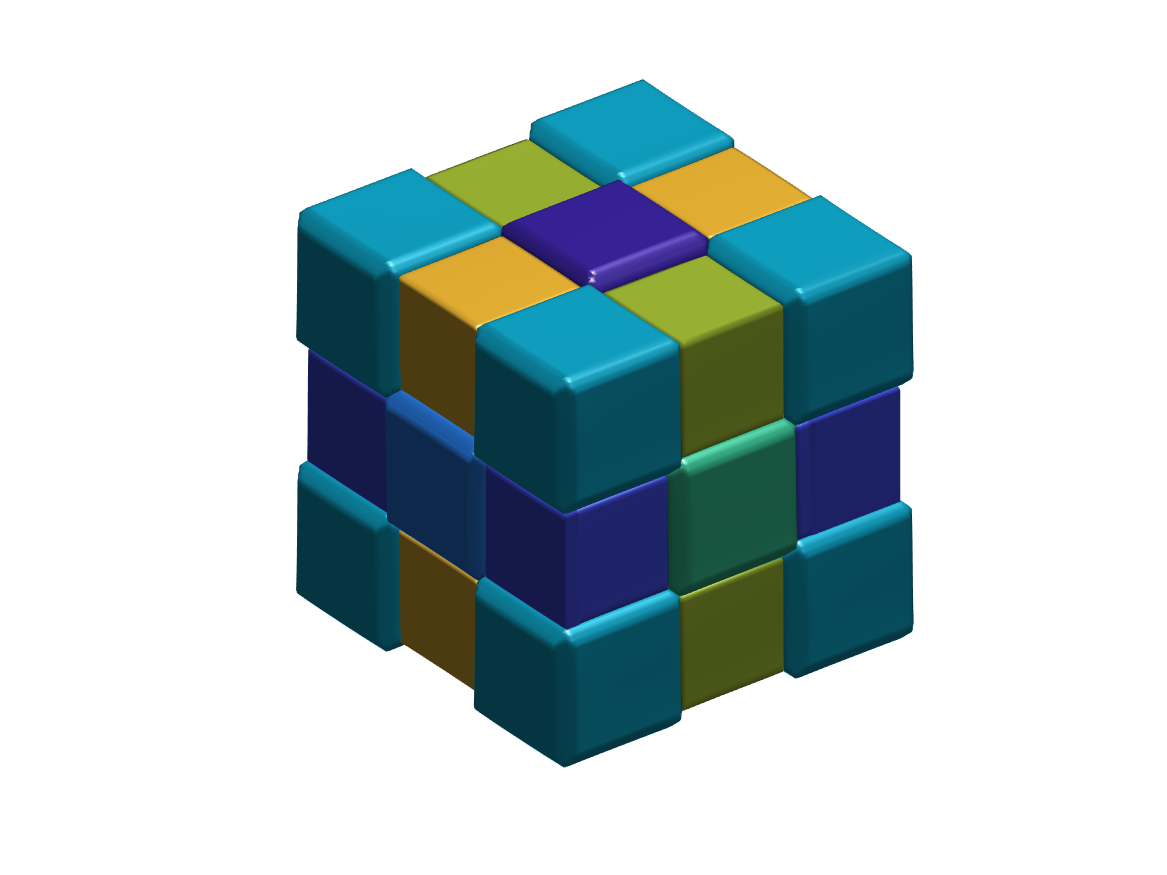}
\includegraphics[scale=0.18,clip,trim= 3cm 1.5cm 3cm 1.2cm]{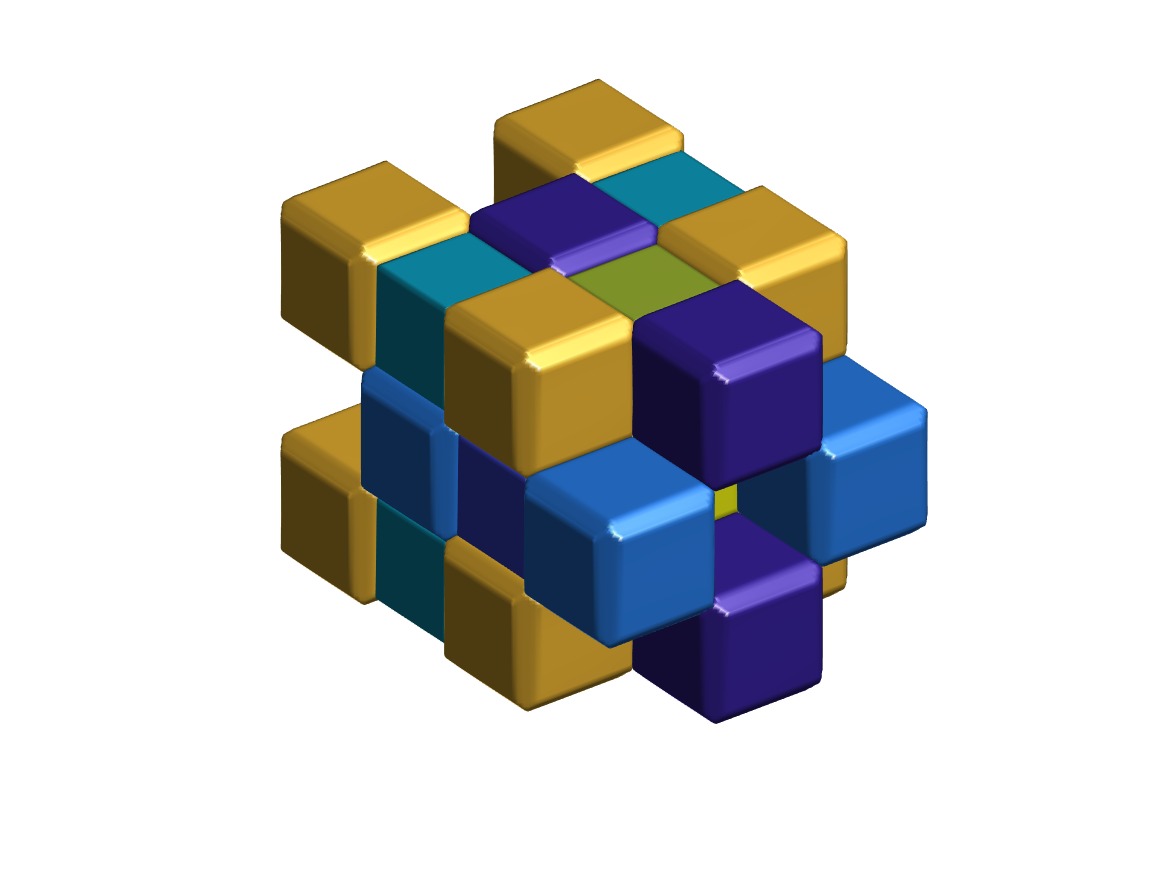}
\includegraphics[scale=0.18,clip,trim= 3cm 1.5cm 3cm 1.2cm]{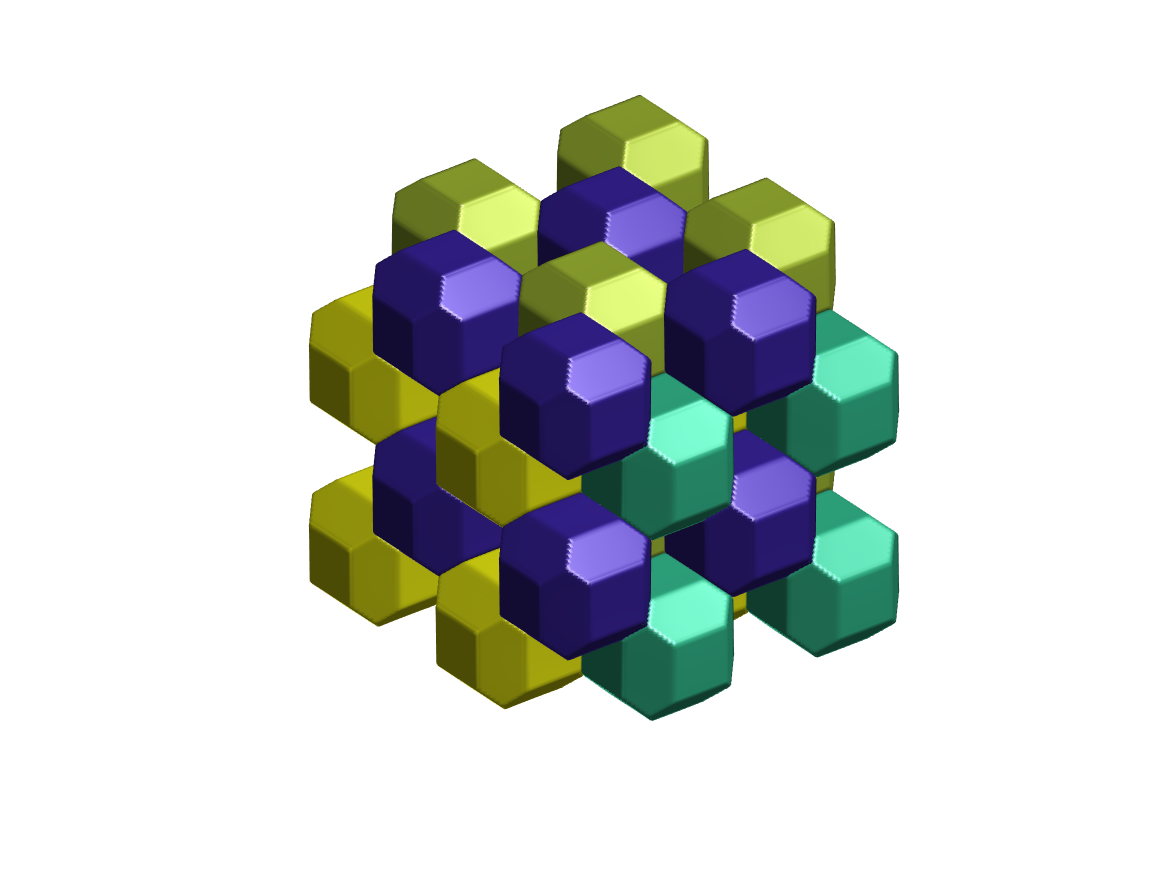}
\includegraphics[scale=0.18,clip,trim= 3cm 1.5cm 3cm 1.2cm]{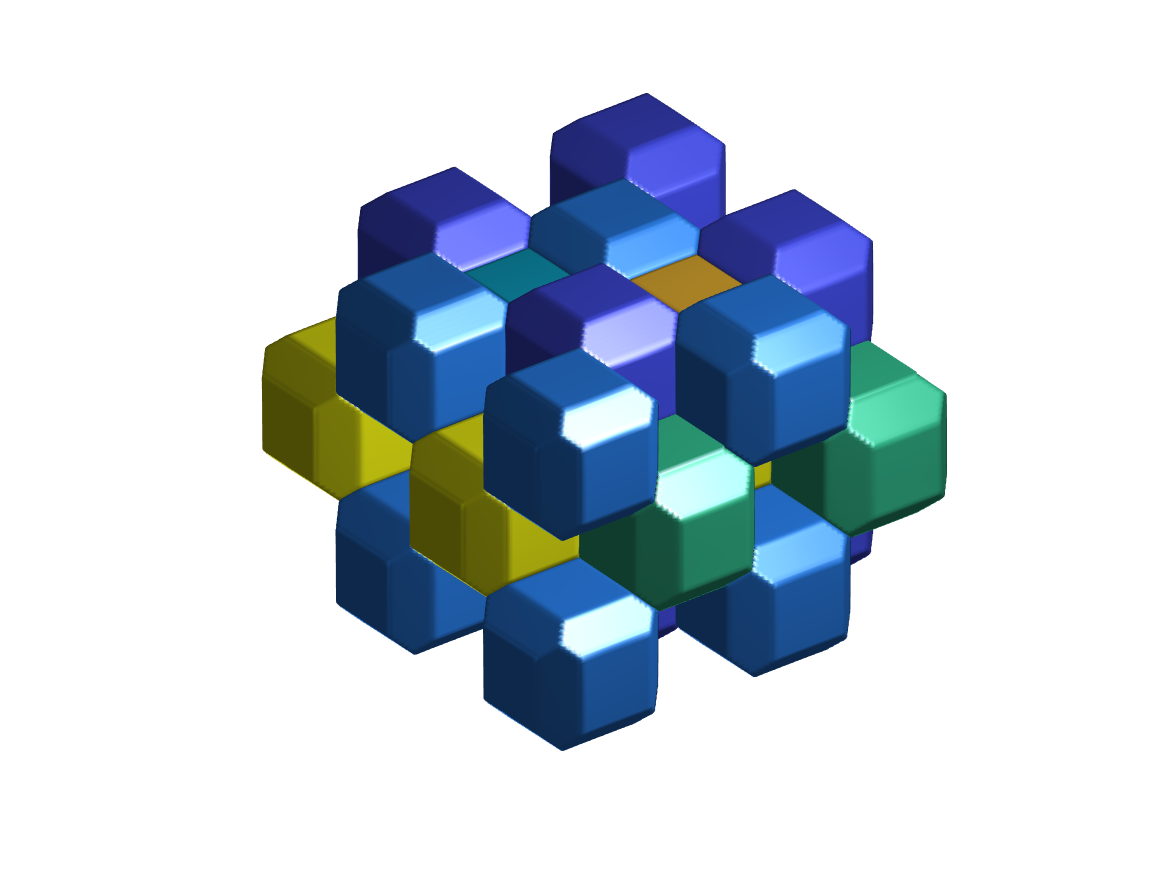}\\
\includegraphics[scale=0.18,clip,trim= 3cm 1.5cm 3cm 1.2cm]{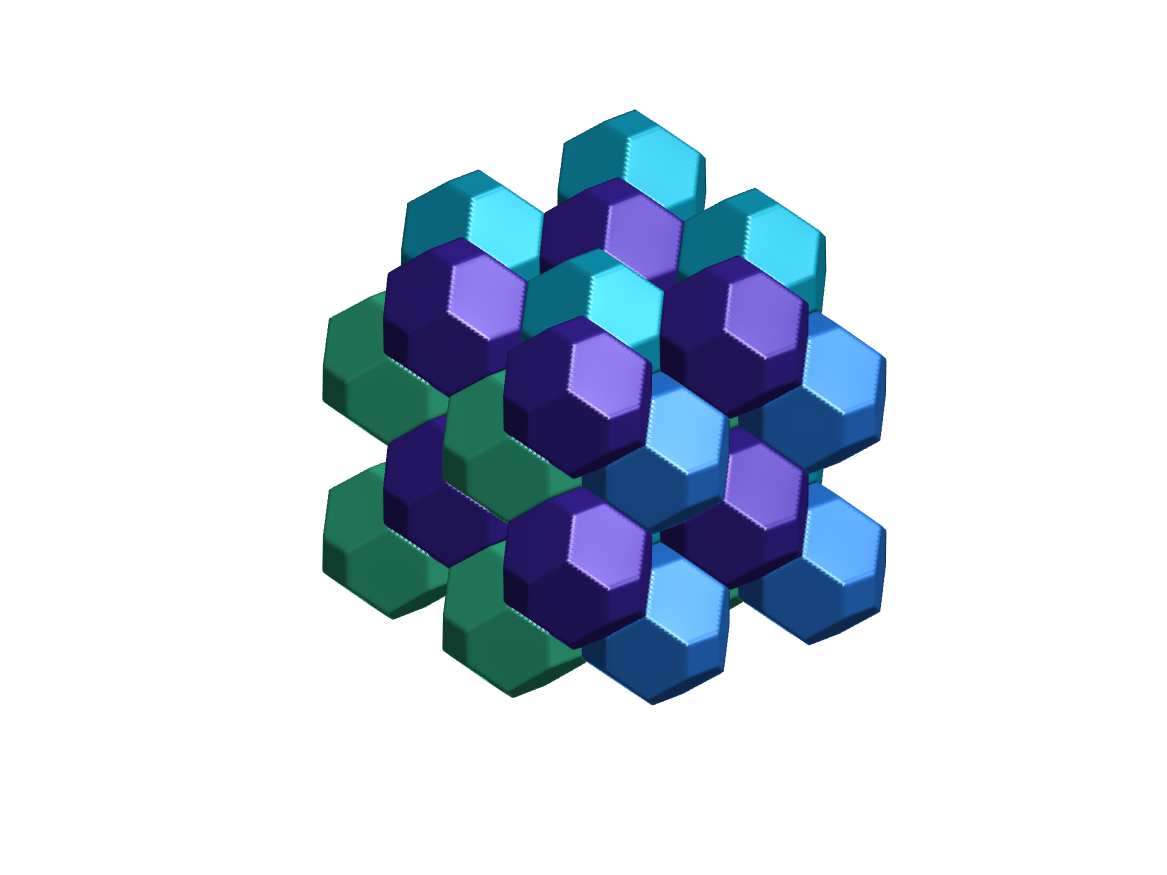}
\includegraphics[scale=0.18,clip,trim= 3cm 1.5cm 3cm 1.2cm]{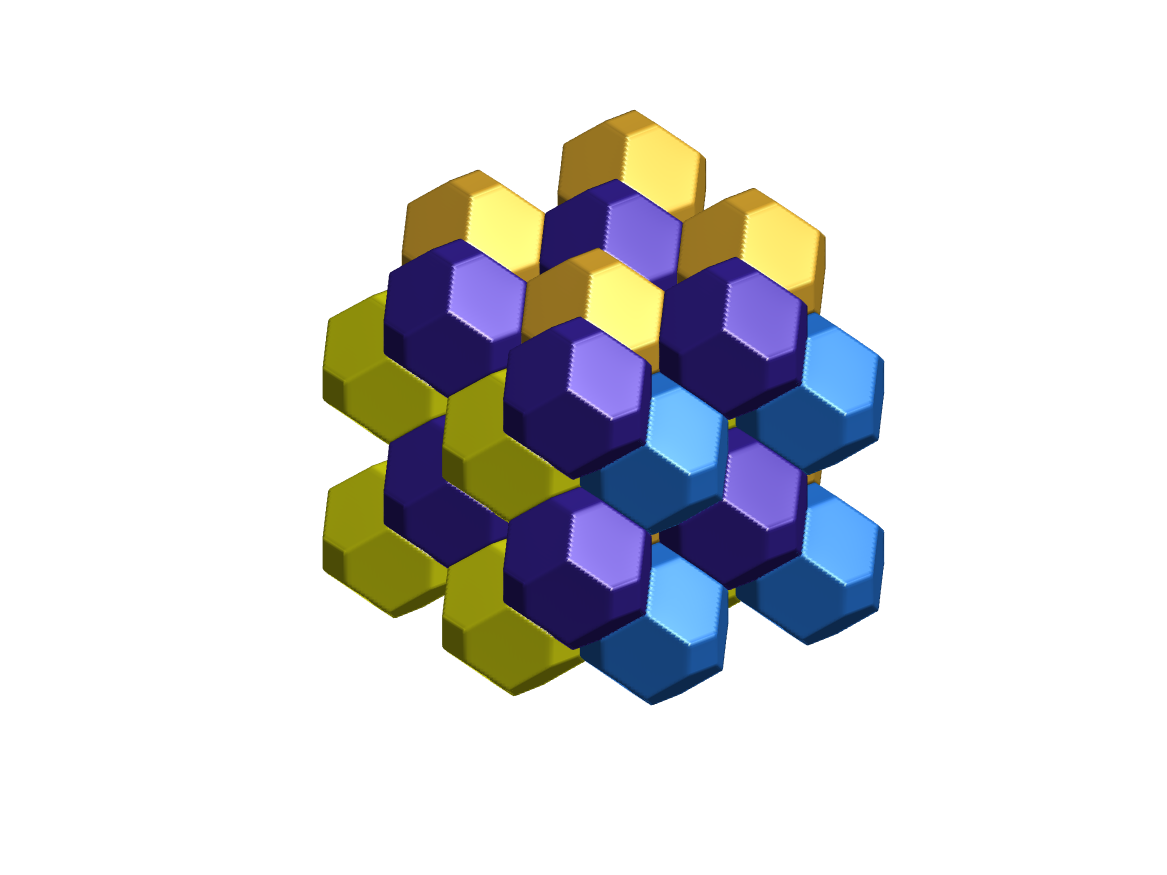}
\includegraphics[scale=0.18,clip,trim= 3cm 1.5cm 3cm 1.2cm]{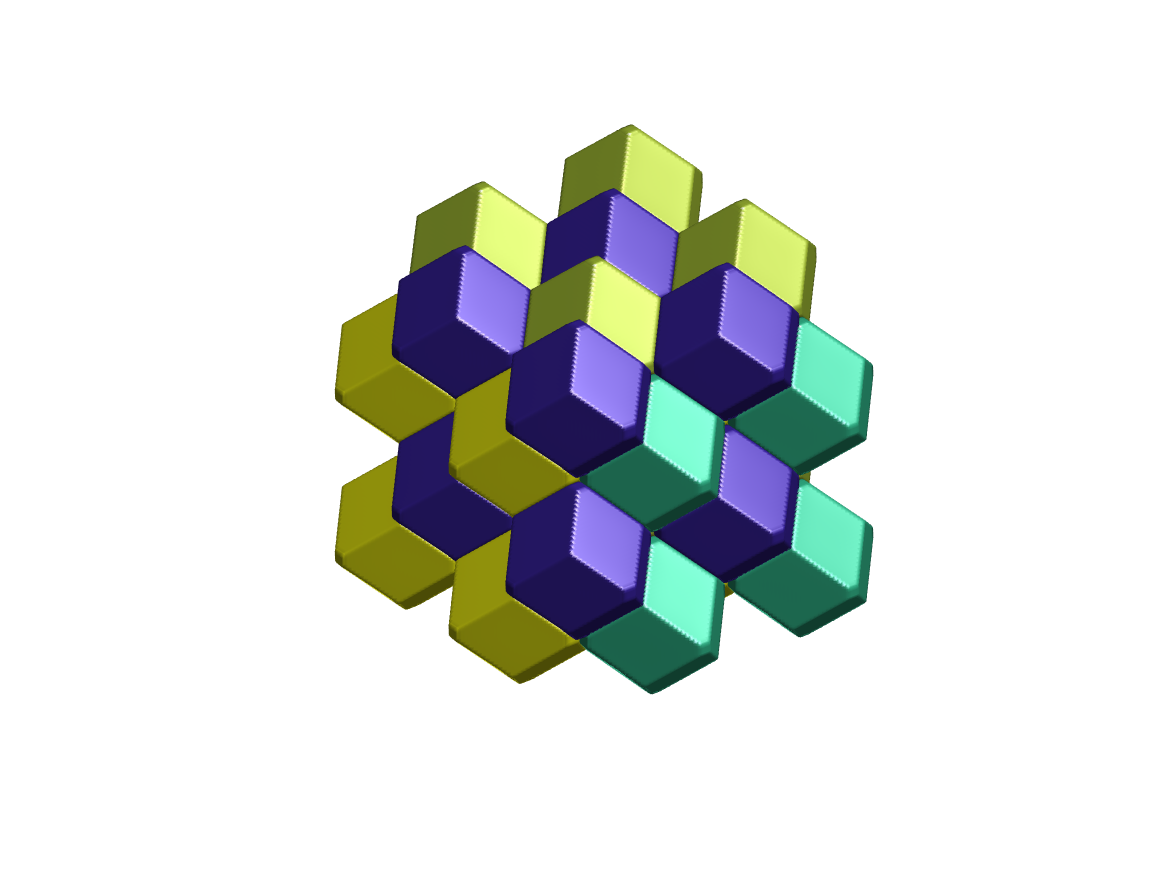}
\includegraphics[scale=0.18,clip,trim= 3cm 1.5cm 3cm 1.2cm]{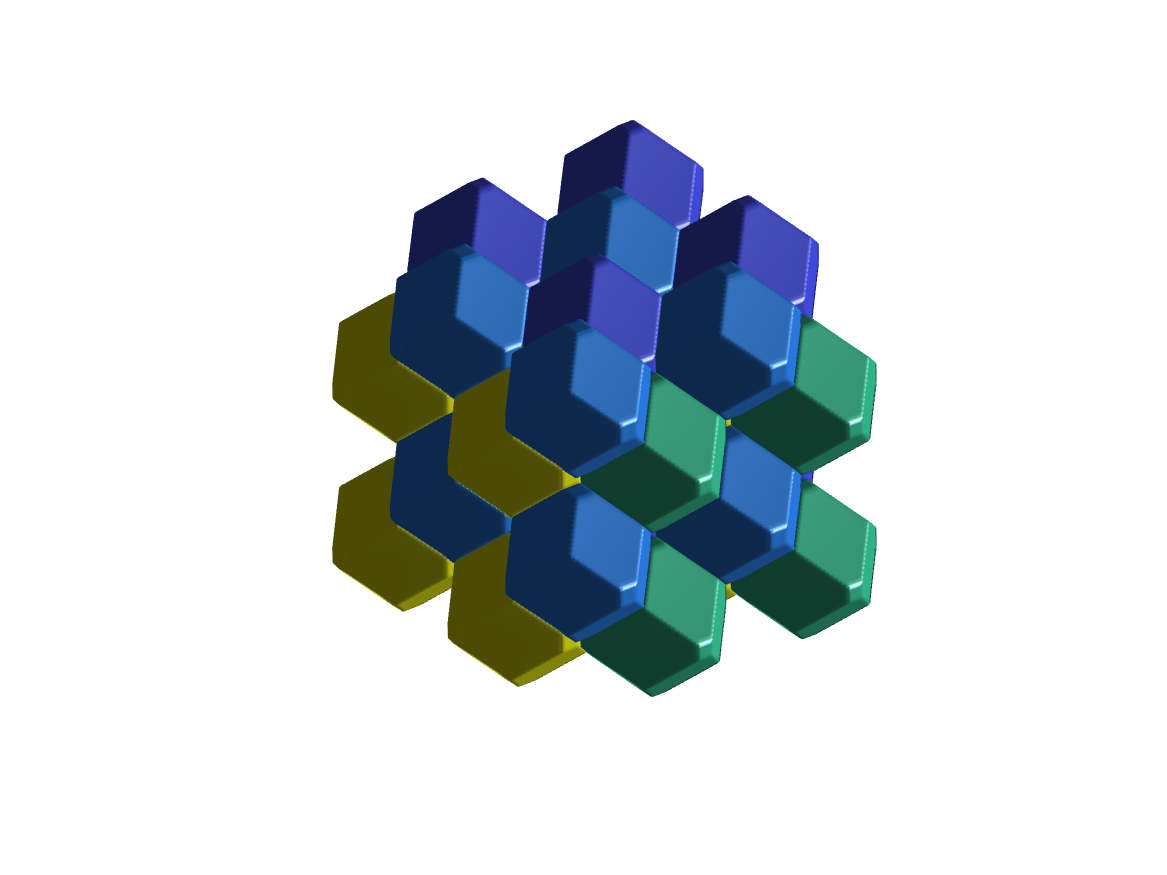}\\
\includegraphics[scale=0.18,clip,trim= 3cm 1.5cm 3cm 1.2cm]{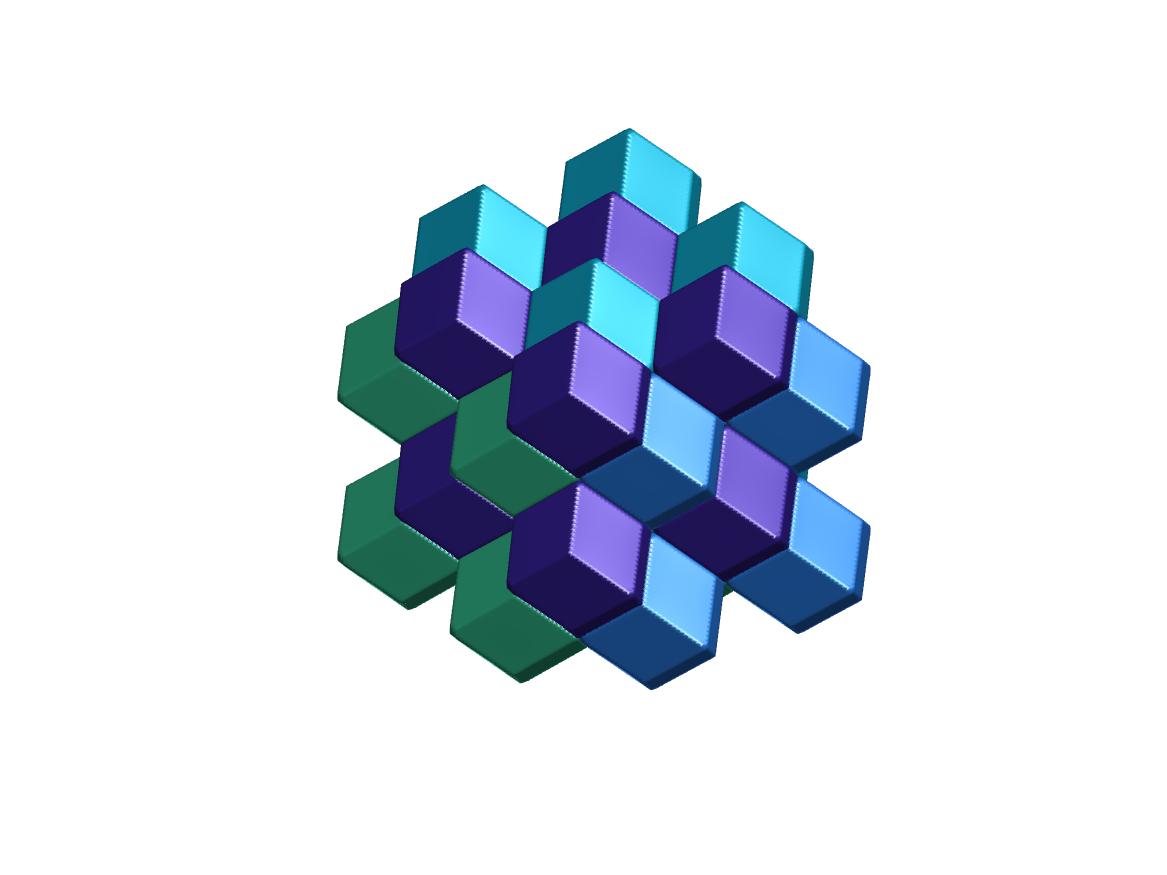}
\includegraphics[scale=0.18,clip,trim= 3cm 1.5cm 3cm 1.2cm]{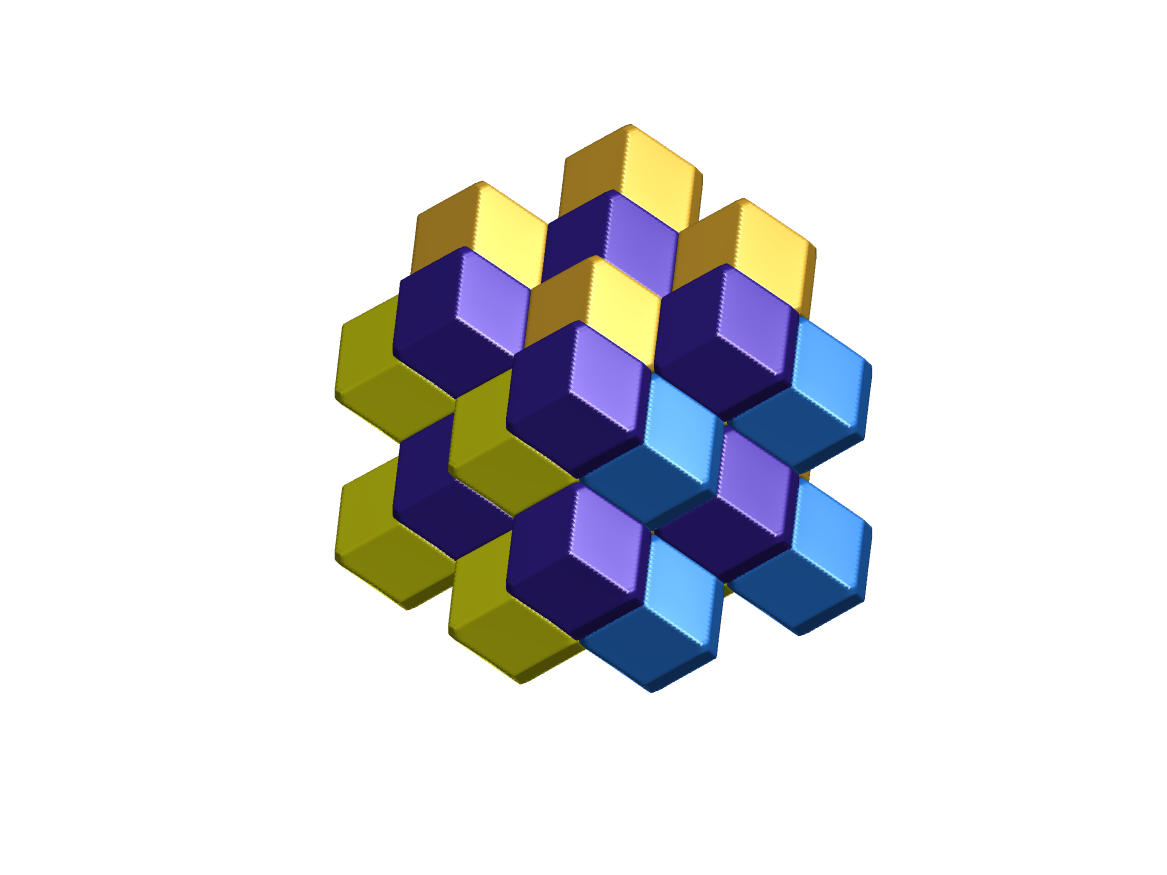}
\includegraphics[scale=0.18,clip,trim= 3cm 1.5cm 3cm 1.2cm]{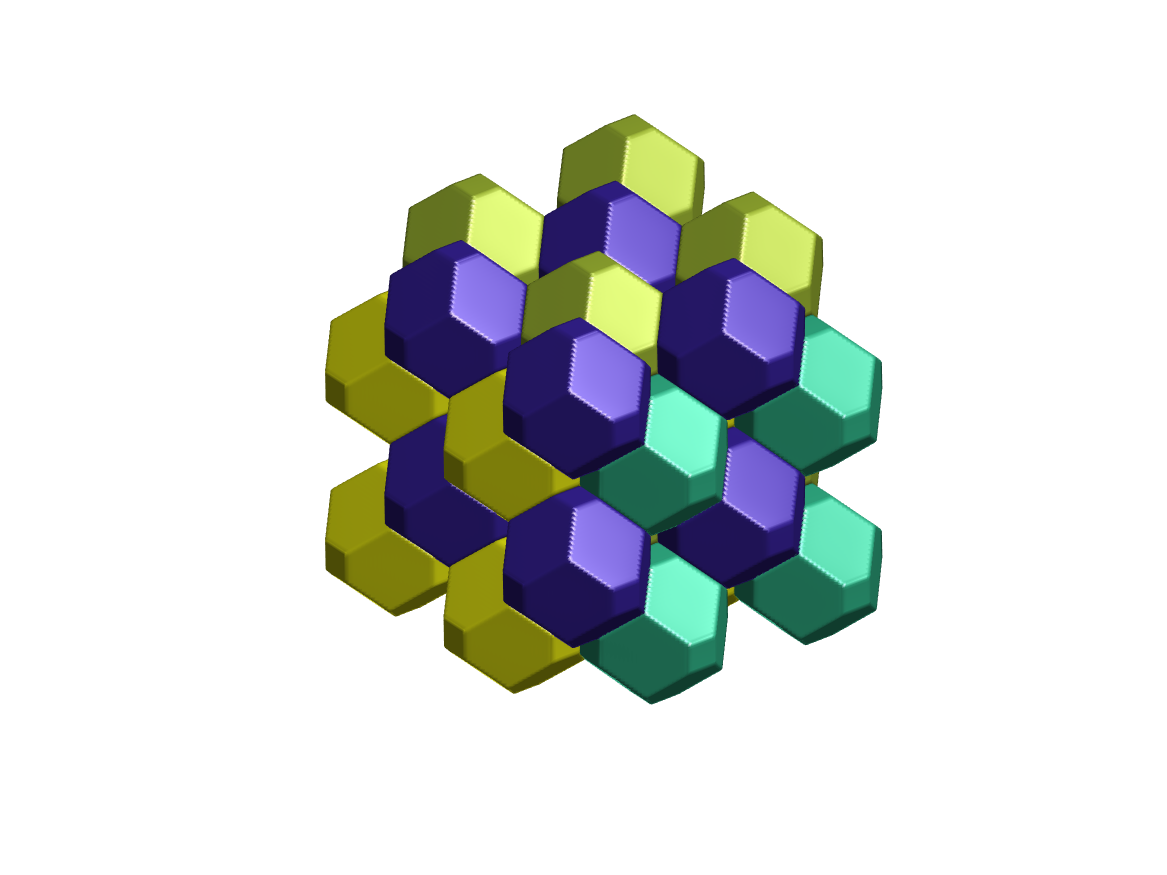}
\includegraphics[scale=0.18,clip,trim= 3cm 1.5cm 3cm 1.2cm]{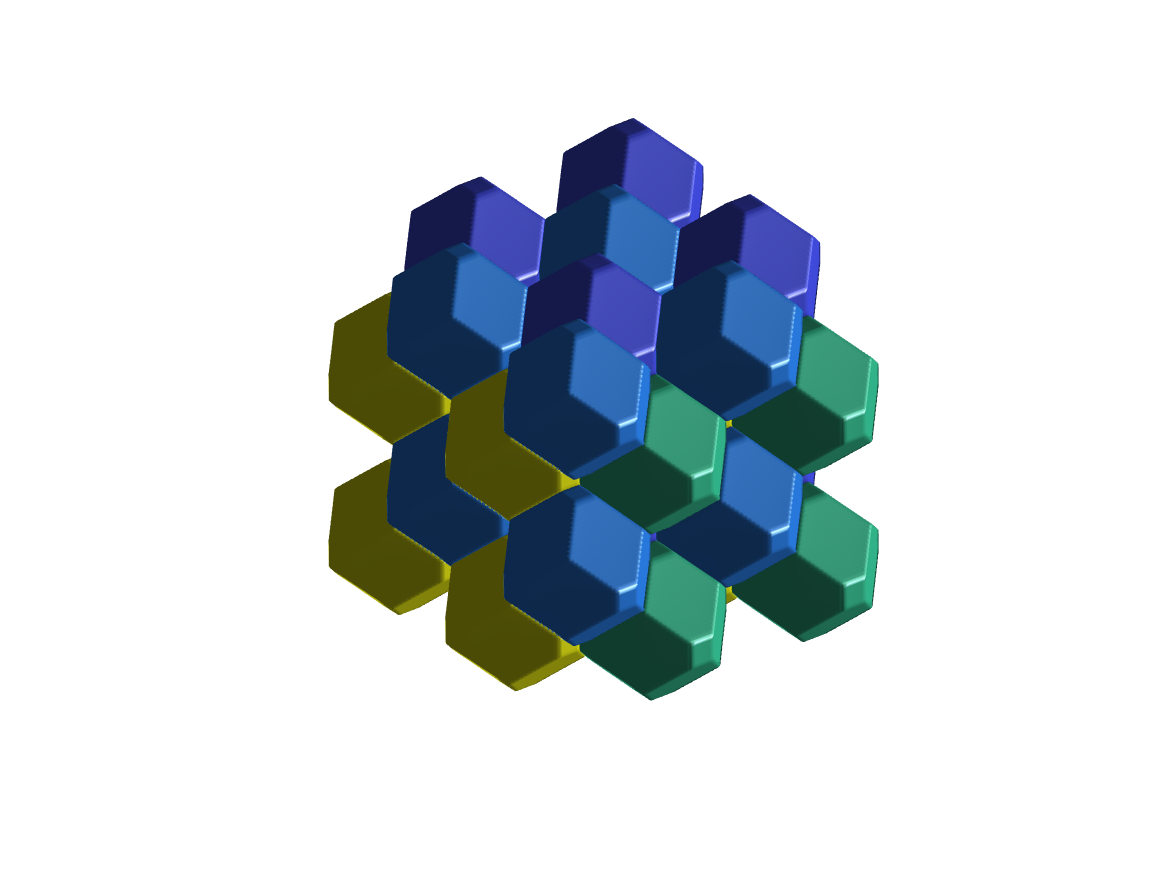}\\
\includegraphics[scale=0.18,clip,trim= 3cm 1.5cm 3cm 1.2cm]{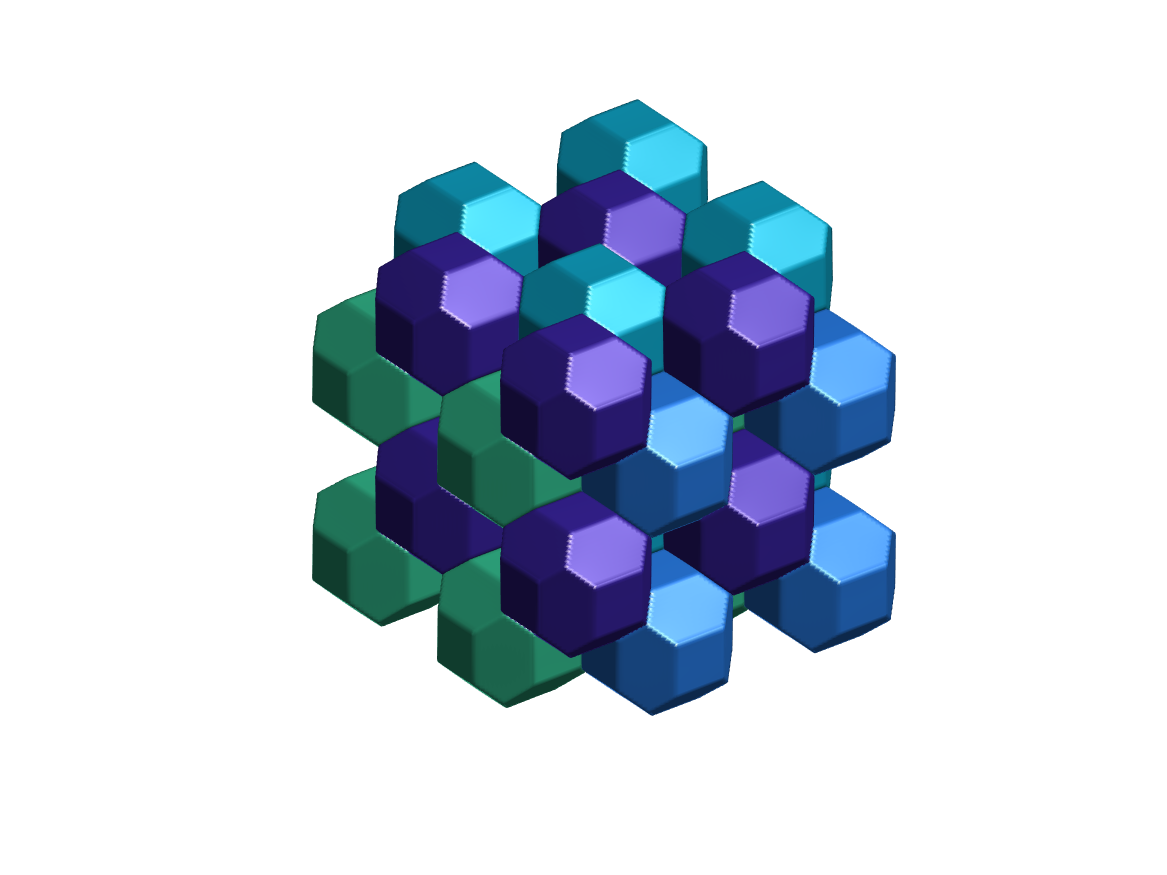}
\includegraphics[scale=0.18,clip,trim= 3cm 1.5cm 3cm 1.2cm]{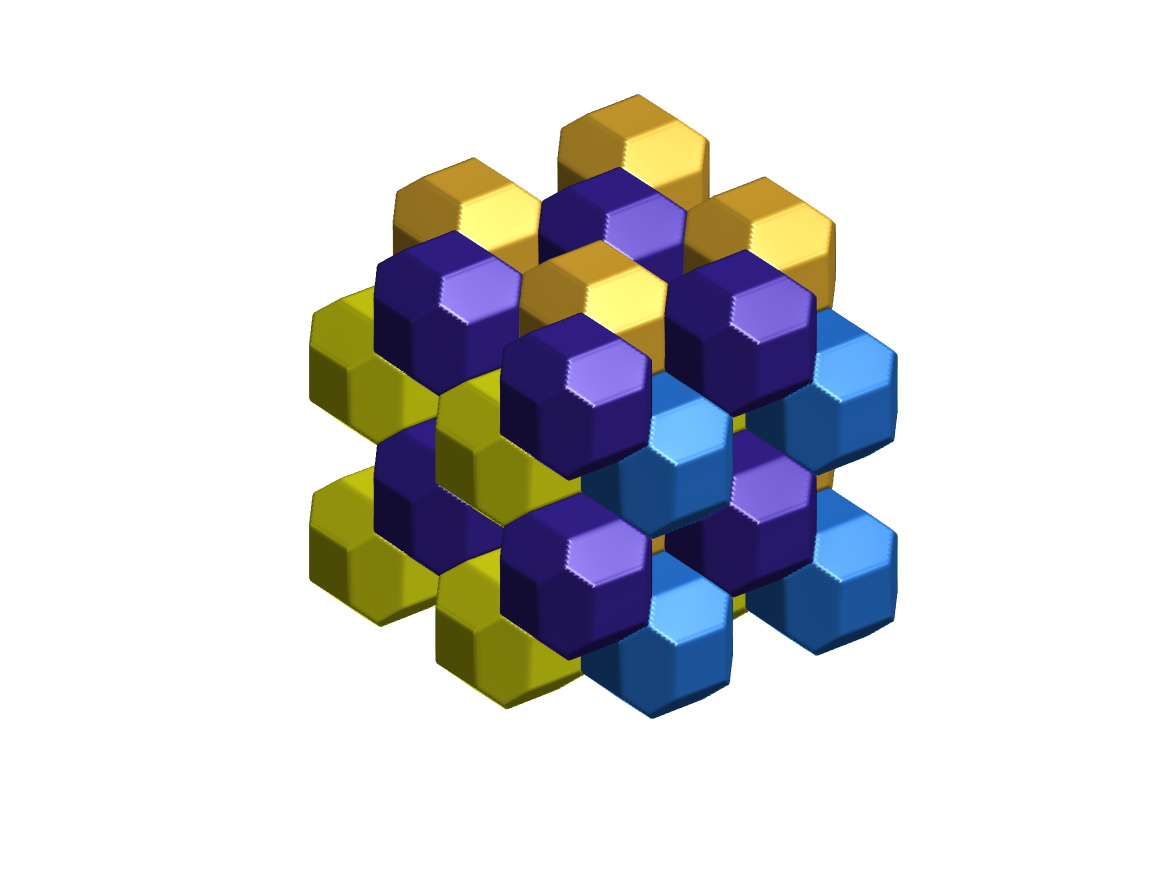}
\includegraphics[scale=0.18,clip,trim= 3cm 1.5cm 3cm 1.2cm]{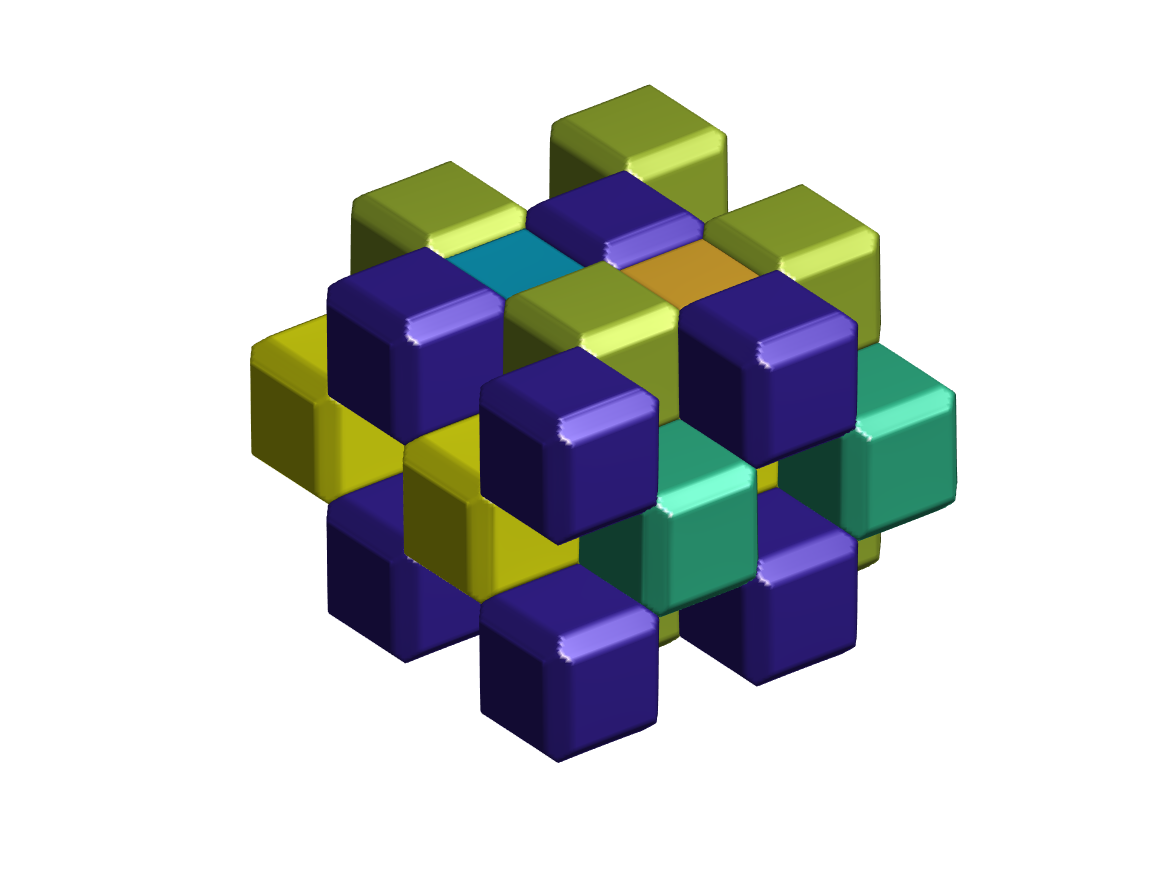}
\includegraphics[scale=0.18,clip,trim= 3cm 1.5cm 3cm 1.2cm]{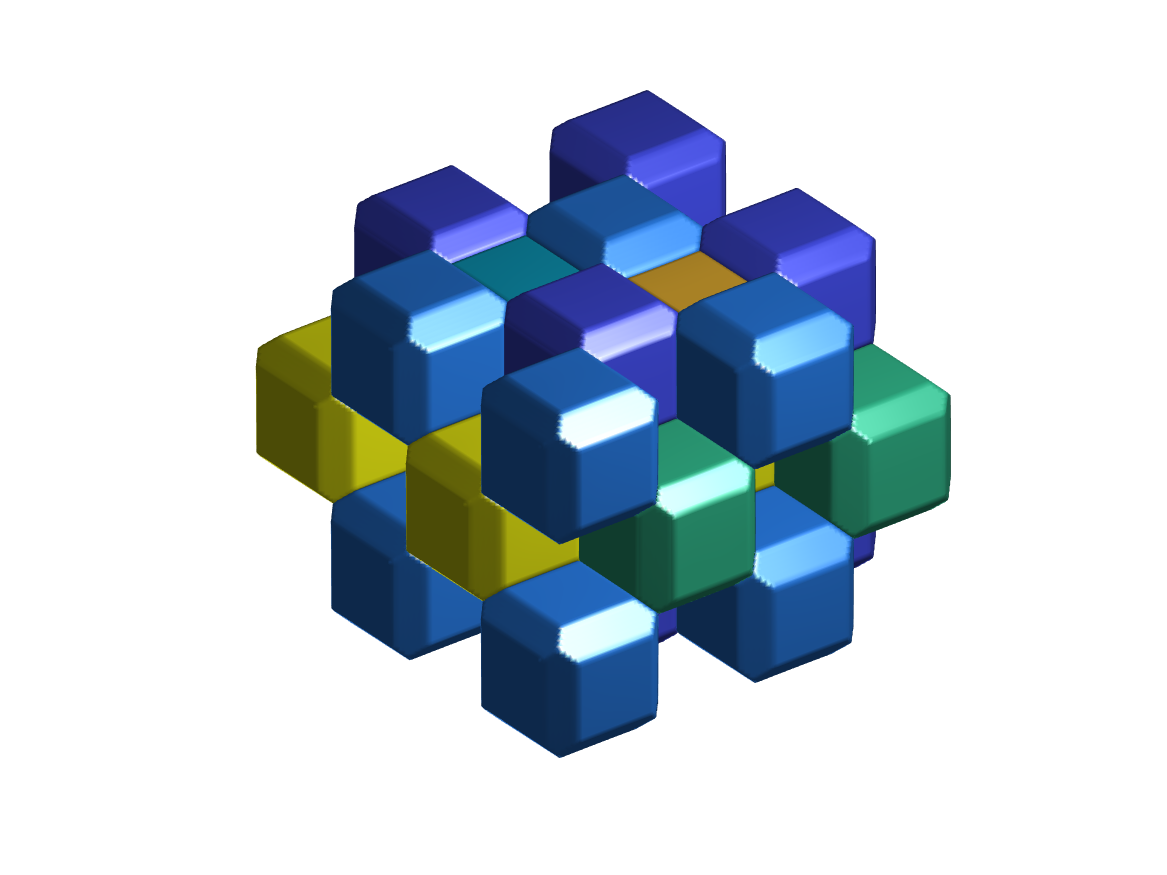}
\caption{A $k=8$ Dirichlet partition of the periodic tesseract, $[-1,1]^4$, by 24-cells. 
The four columns correspond to the slides perpendicular to the $x_1-$, $x_2-$, $x_3-$, and $x_4-$axis respectively. 
The eight rows correspond to the slices at 
$x_j=$-1, -0.75, -0.5, -0.25, 0, 0.25, 0.5, 0.75, respectively. 
The CPU time was $9803$ seconds. } 
\label{fig:4d}
\end{figure}

For $k=8$ and initialization using a random tessellation, we obtain a partition of the tesseract as shown in Figure~\ref{fig:4d}. 
The four columns of this plot correspond to slices perpendicular to the  $x_1-$, $x_2-$, $x_3-$, and $x_4-$axes, respectively. 
The eight rows  correspond to the slices at $x_j=$-1, -0.75, -0.5, -0.25, 0, 0.25, 0.5, and 0.75, respectively. 
The partition obtained is known as a \emph{24-cell honeycomb}, which is a tessellation by 24-cells. 
In the experiment, the tesseract was discretized  by $64^4$ grid points and $\tau = 0.0625$. The CPU time was $9803$ seconds.

\subsection{Results for Sphere} \label{s:sphere}
Finally, we consider Dirichlet partitions for a sphere. It has been conjectured that the  3 Dirichlet partition of the sphere is the ``Y-partition'' \cite{Helffer2010}. Dirichlet partitions have been computed on the sphere for several values of $k$, see  \cite{osting2013minimal,elliott2015computational,bogosel2017efficient}. In this section, we compute Dirichlet partitions using Algorithm~\ref{a:MBO} for the sphere. Our results are consistent with previous results. 

\begin{figure}[ht]
\includegraphics[scale=0.15,clip,trim= 12cm 5cm 10cm 5cm]{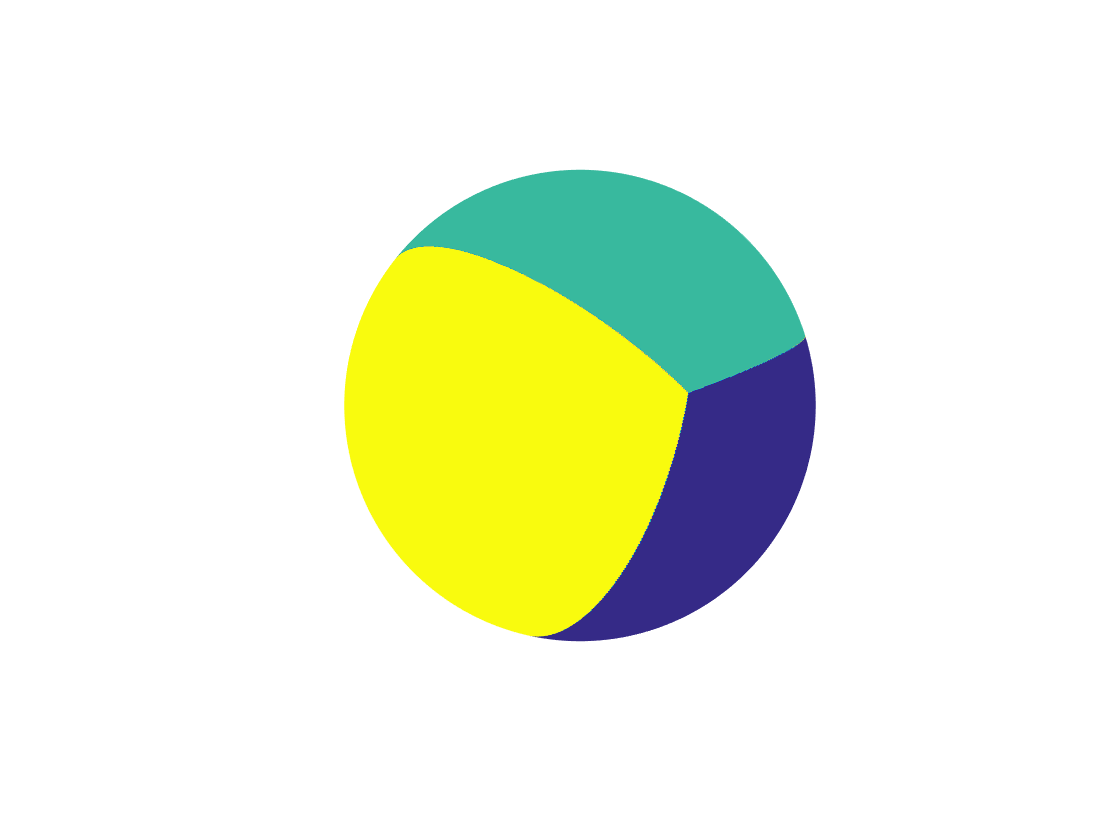} \  \  \  \ 
\includegraphics[scale=0.12,clip,trim= 8cm 2cm 7cm 2cm]{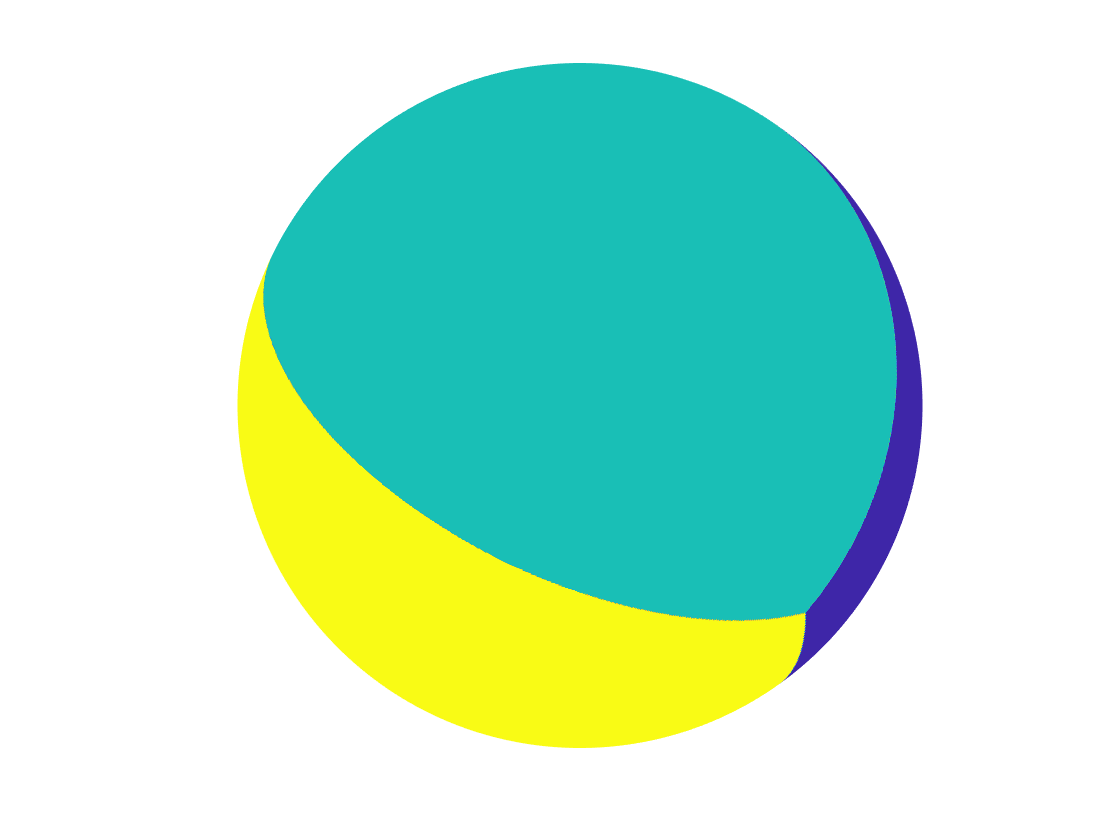}
\includegraphics[scale=0.12,clip,trim= 8cm 2cm 7cm 2cm]{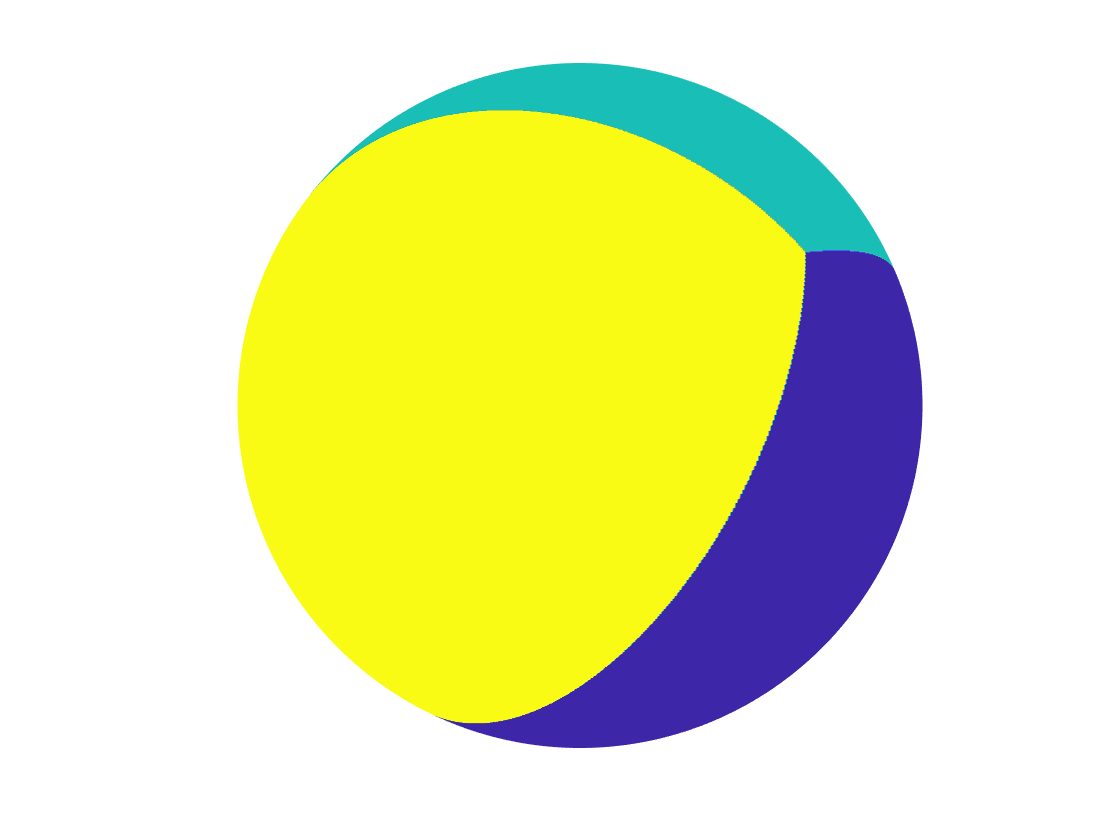}
\includegraphics[scale=0.12,clip,trim= 8cm 2cm 7cm 2cm]{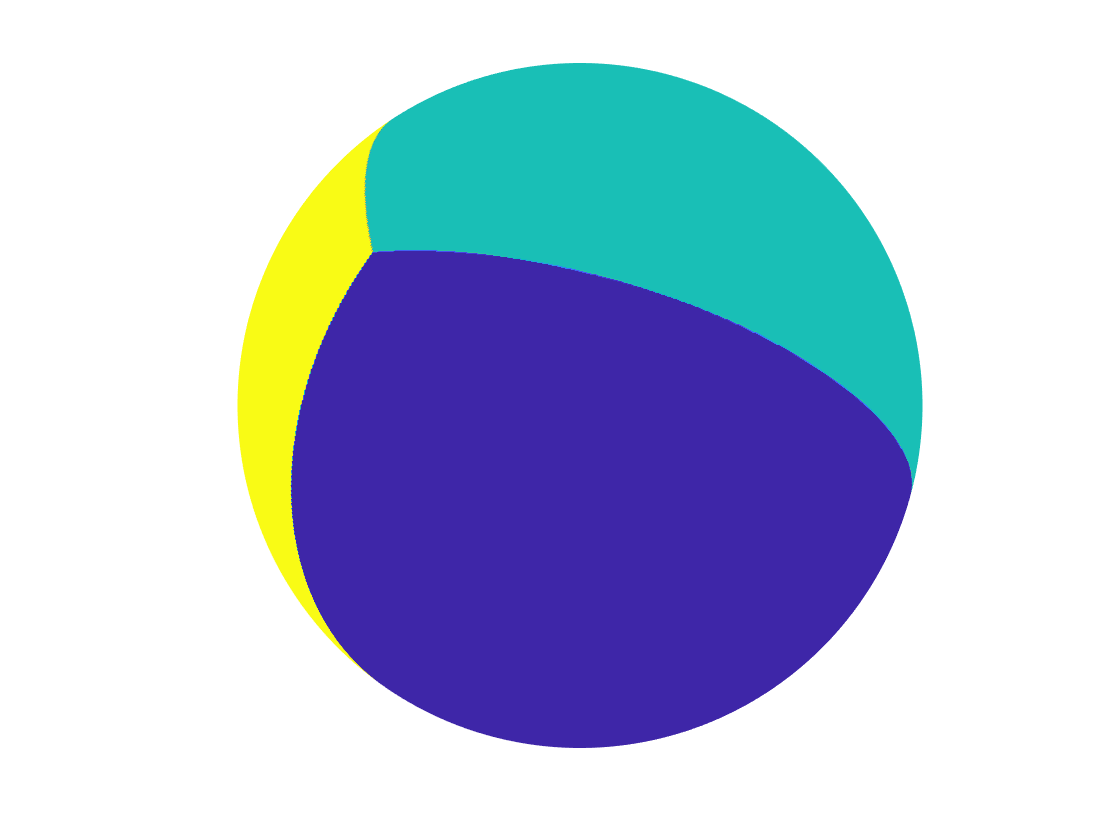}

\includegraphics[scale=0.15,clip,trim= 12cm 5cm 10cm 5cm]{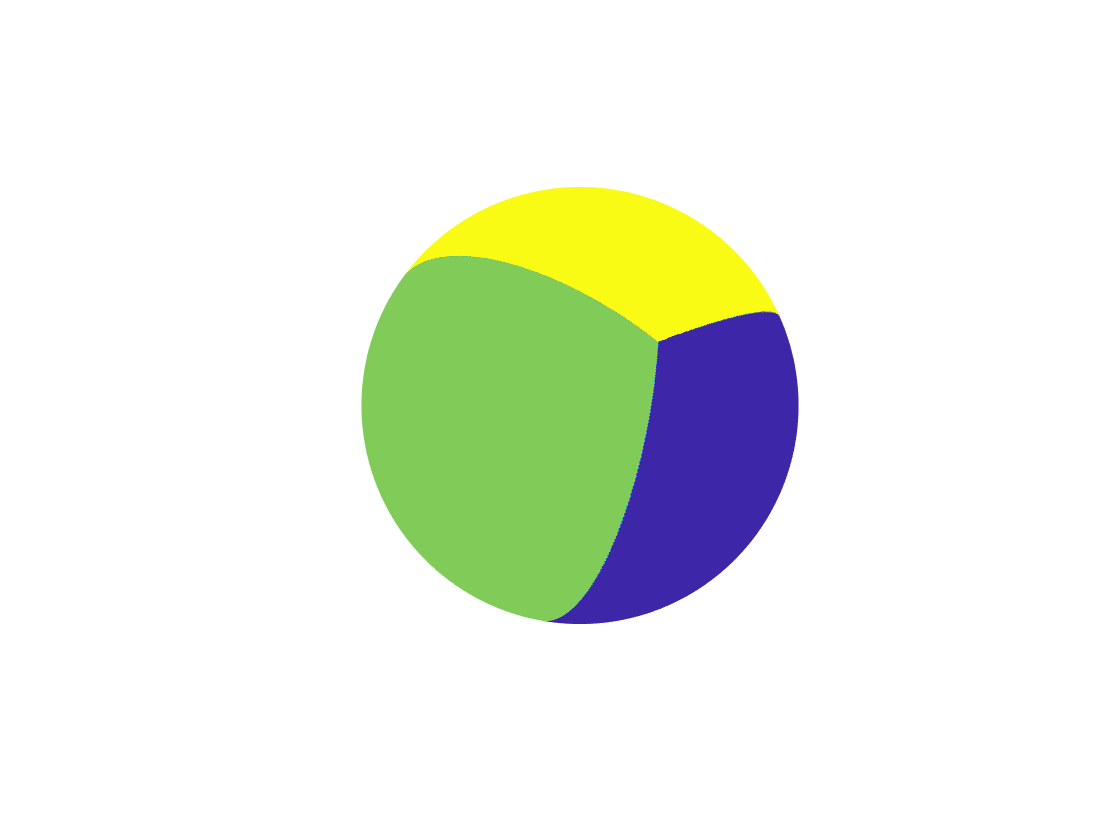} \  \  \ \ 
\includegraphics[scale=0.12,clip,trim= 8cm 2cm 7cm 2cm]{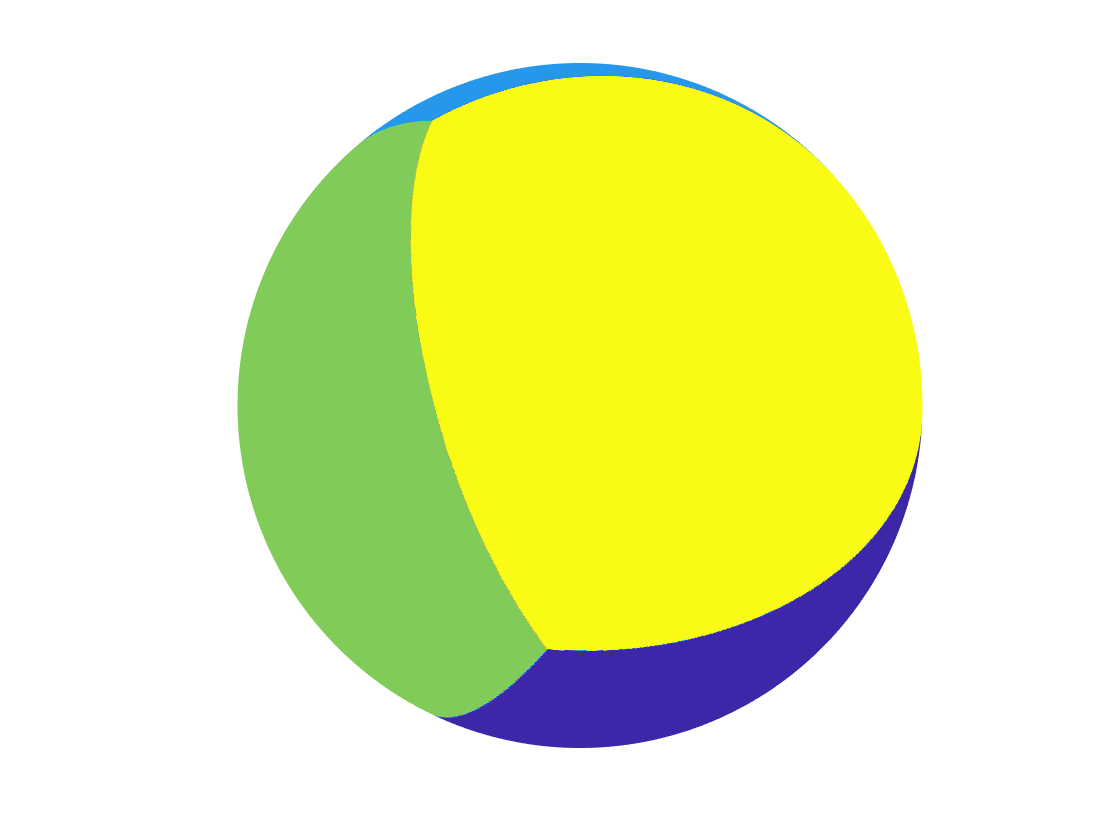}
\includegraphics[scale=0.12,clip,trim= 8cm 2cm 7cm 2cm]{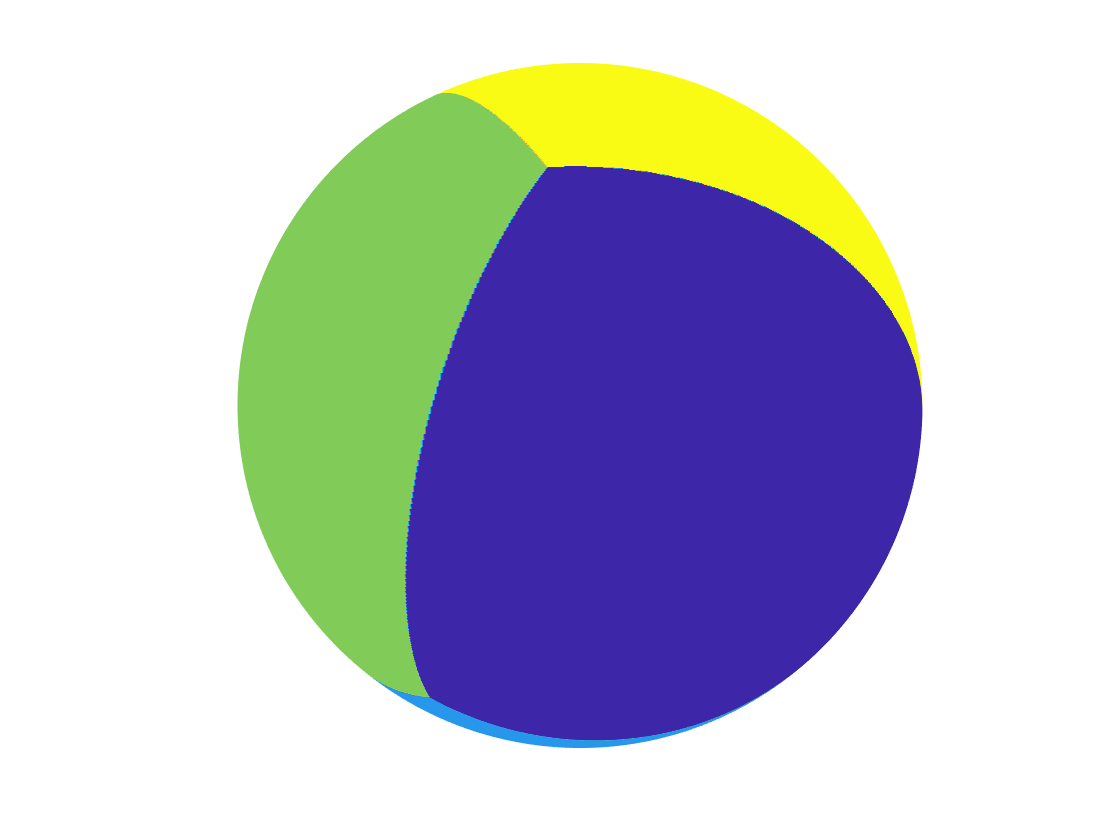}
\includegraphics[scale=0.12,clip,trim= 8cm 2cm 7cm 2cm]{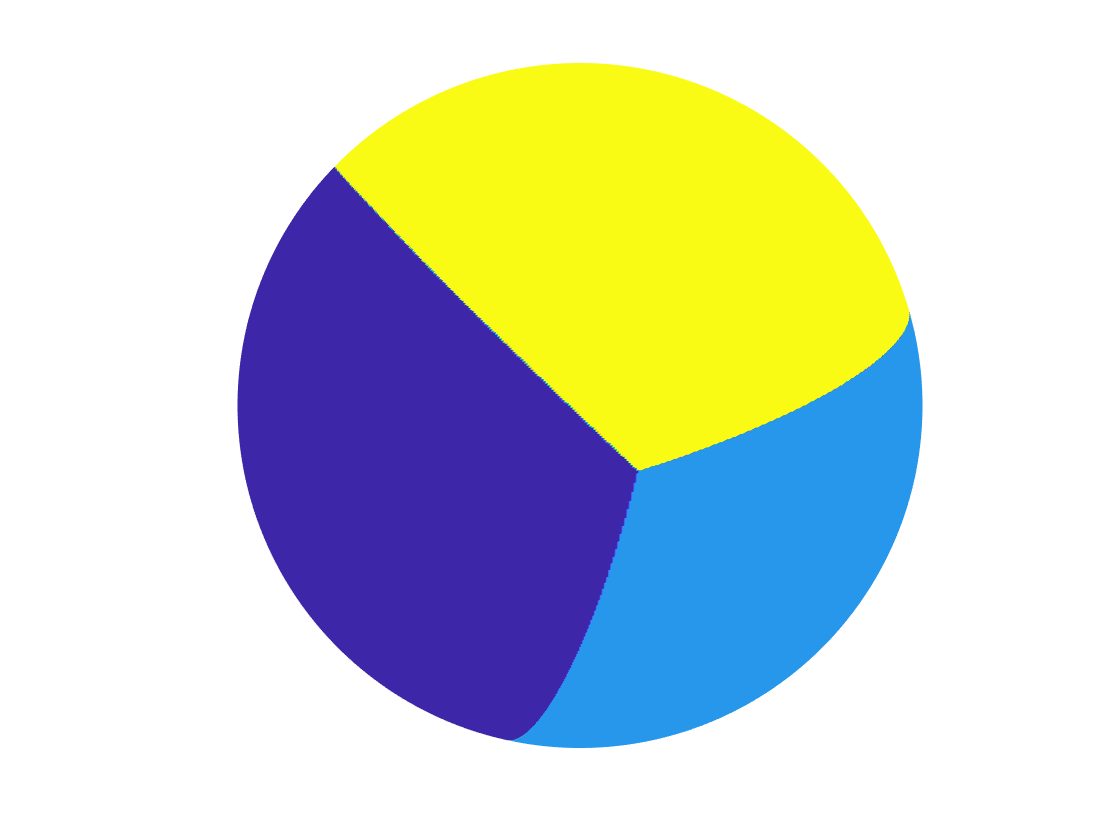}

\includegraphics[scale=0.15,clip,trim= 12cm 5cm 10cm 5cm]{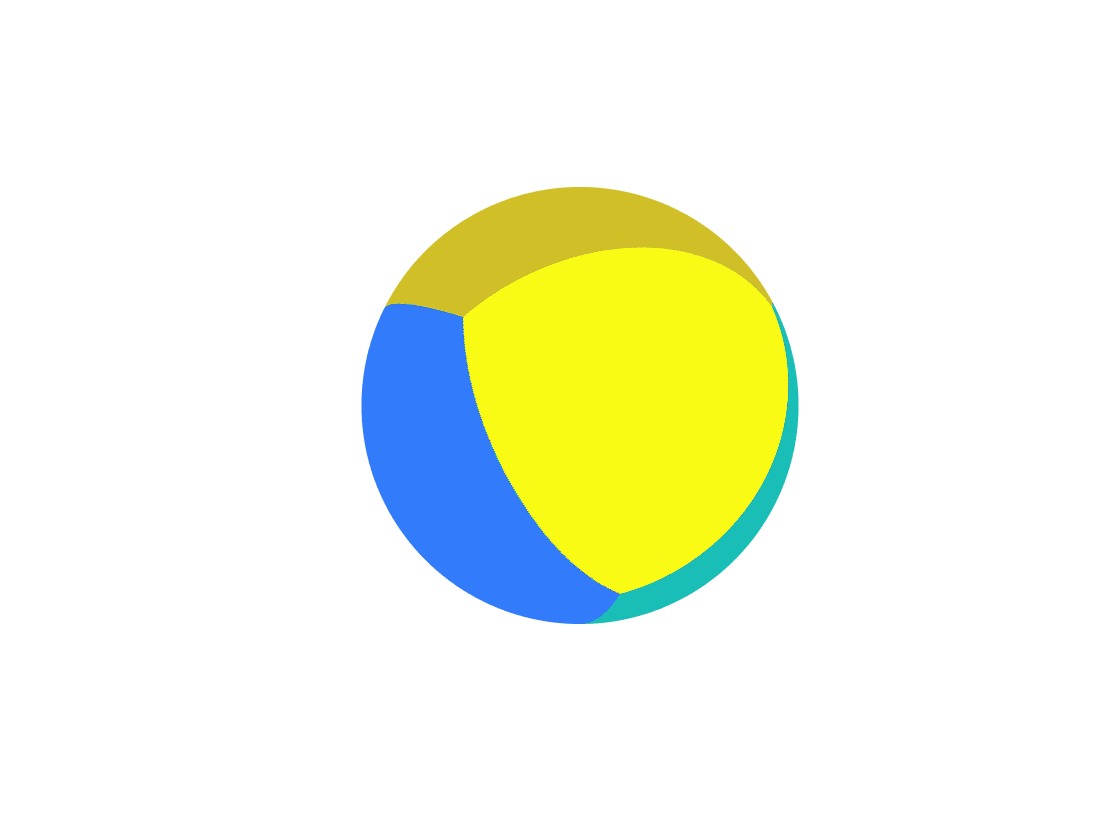} \  \  \ \ 
\includegraphics[scale=0.12,clip,trim= 8cm 2cm 7cm 2cm]{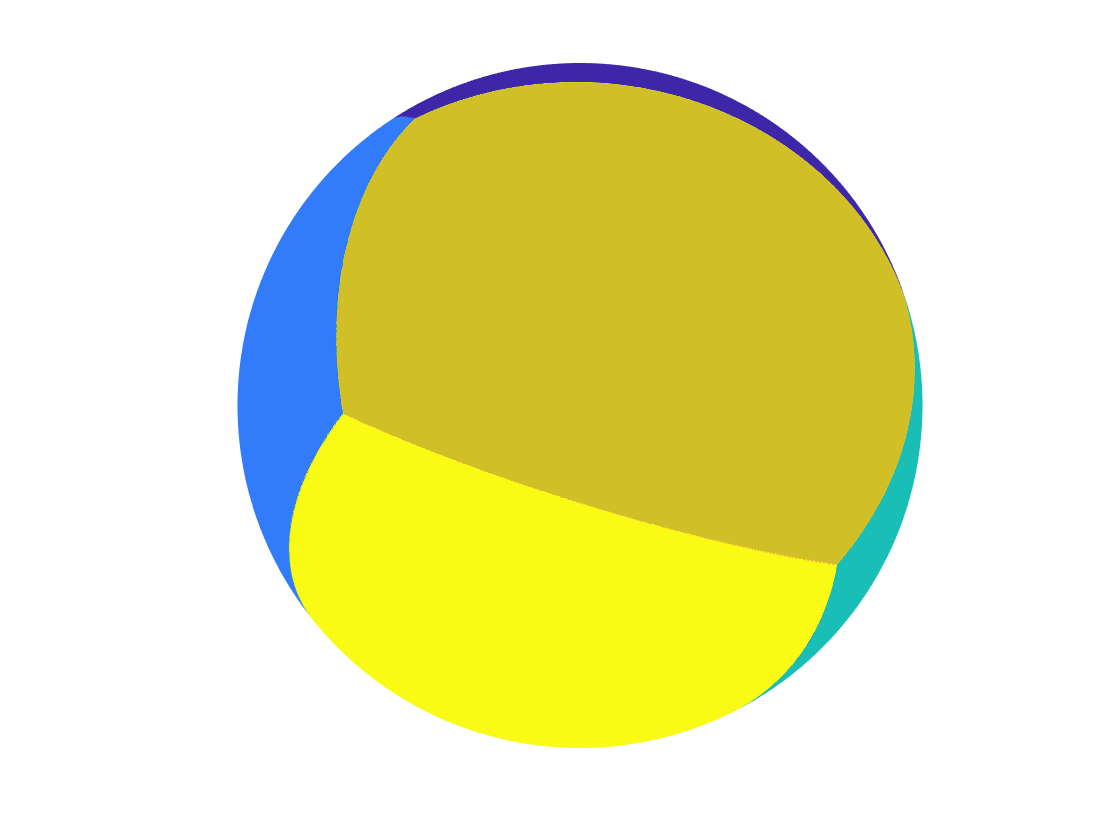}
\includegraphics[scale=0.12,clip,trim= 8cm 2cm 7cm 2cm]{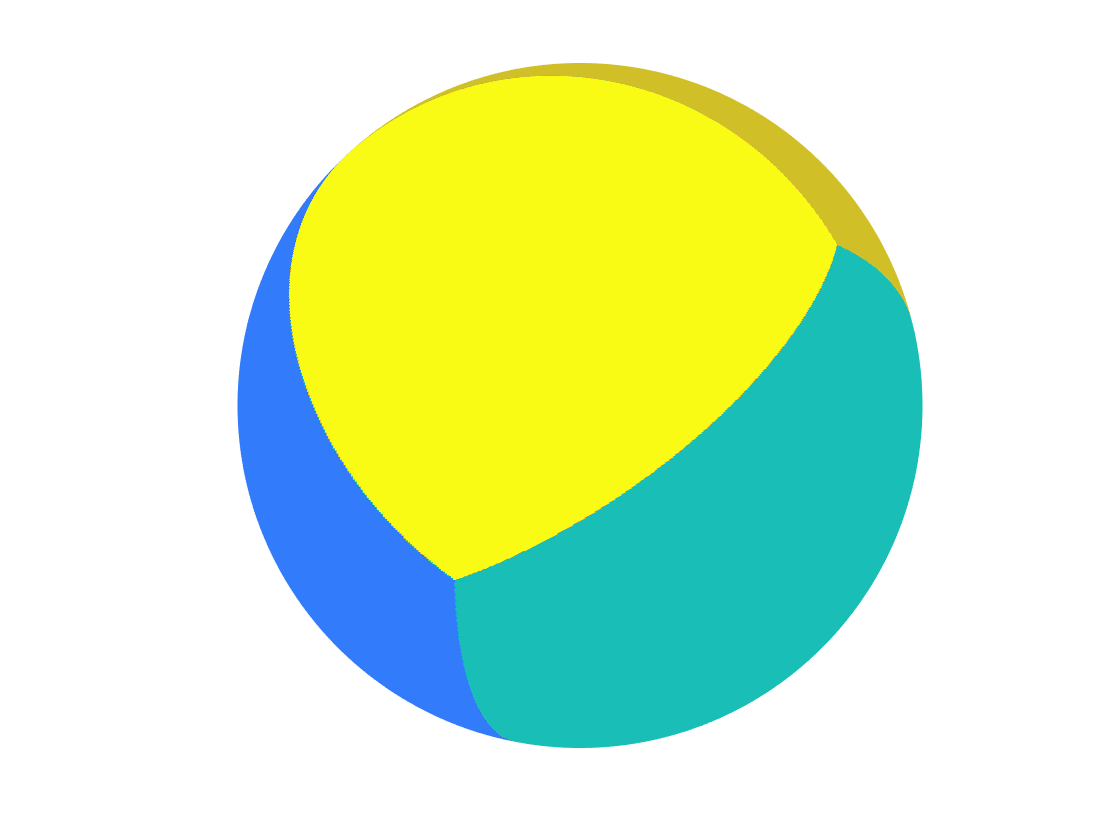}
\includegraphics[scale=0.12,clip,trim= 8cm 2cm 7cm 2cm]{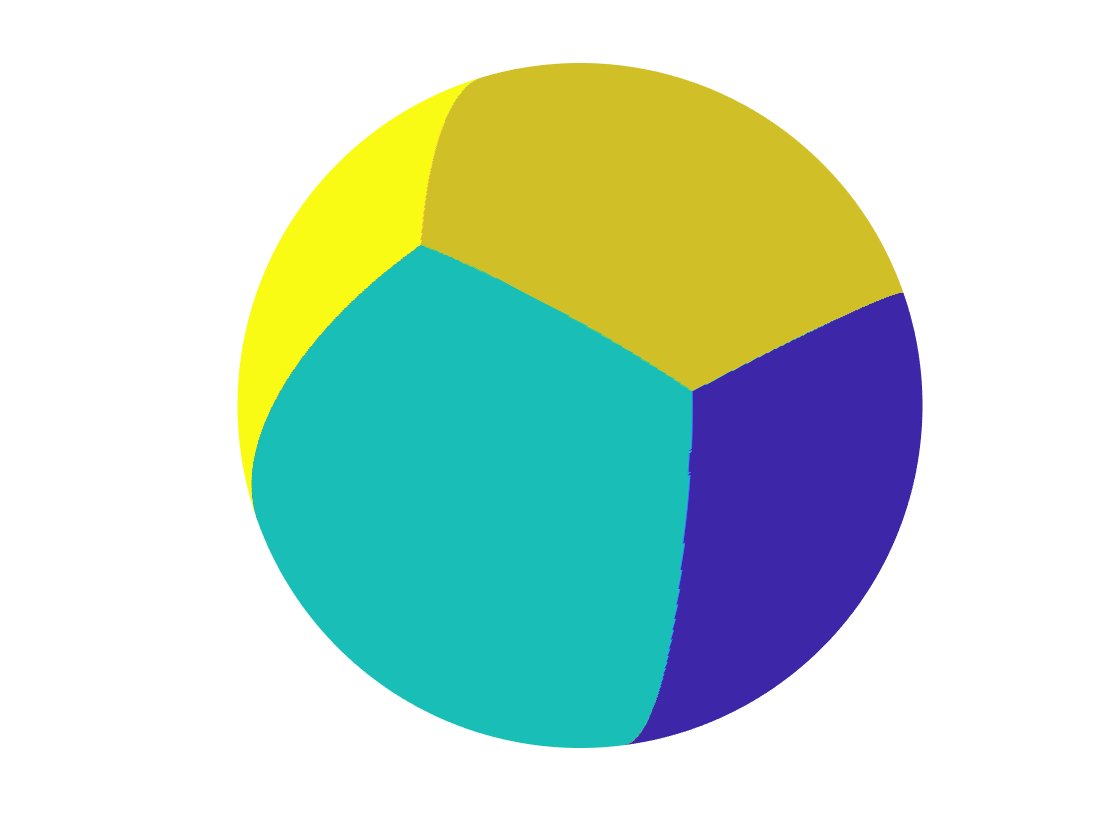}

\includegraphics[scale=0.15,clip,trim= 12cm 5cm 10cm 5cm]{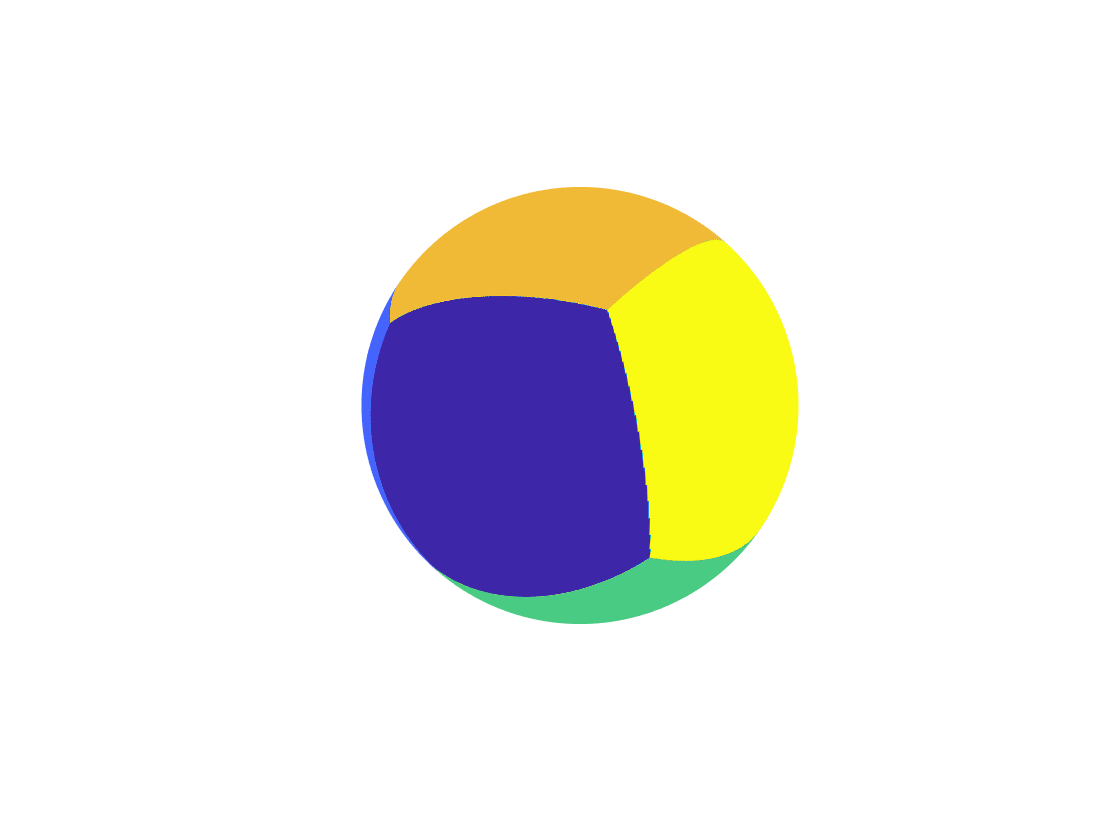} \  \  \ \ 
\includegraphics[scale=0.12,clip,trim= 8cm 2cm 7cm 2cm]{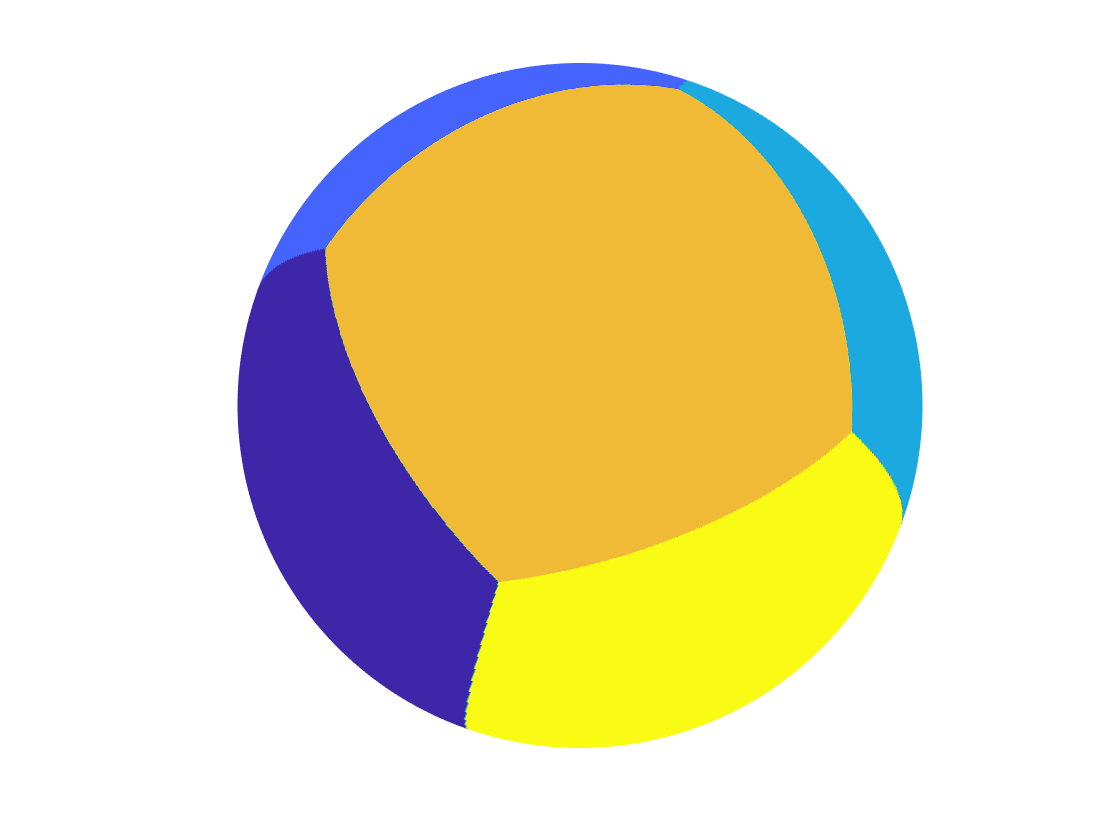}
\includegraphics[scale=0.12,clip,trim= 8cm 2cm 7cm 2cm]{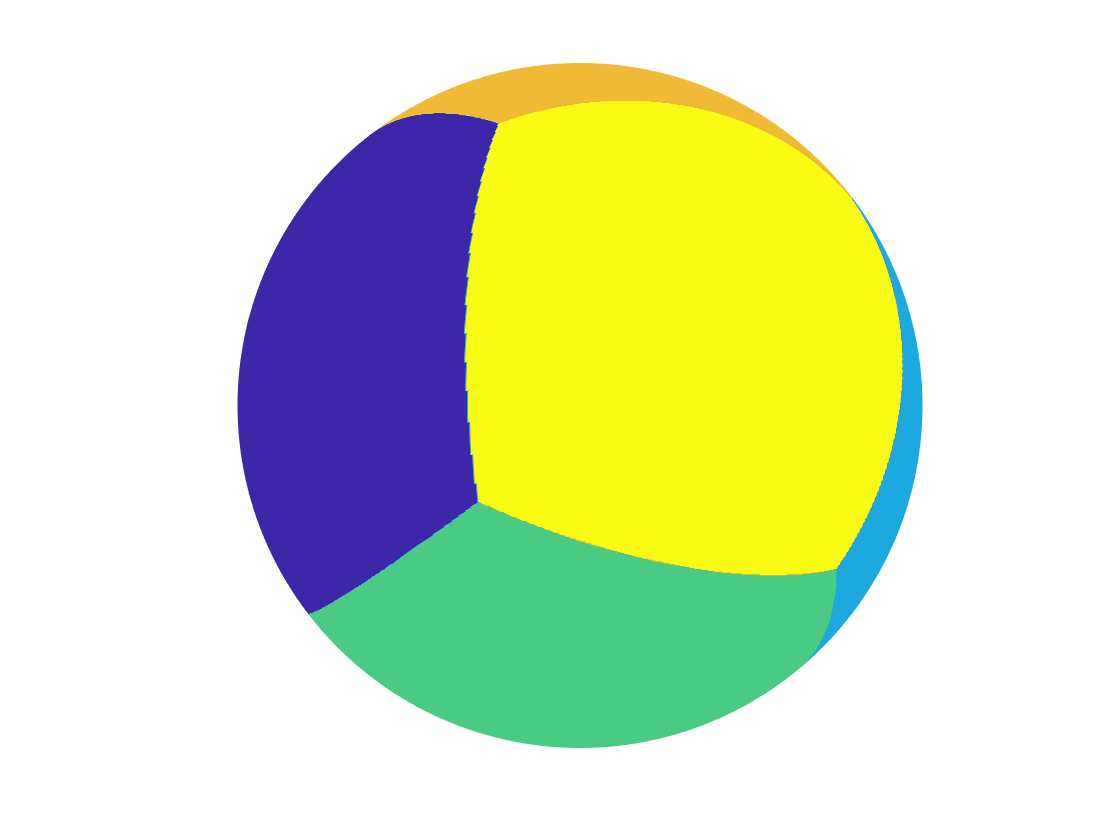}
\includegraphics[scale=0.12,clip,trim= 8cm 2cm 7cm 2cm]{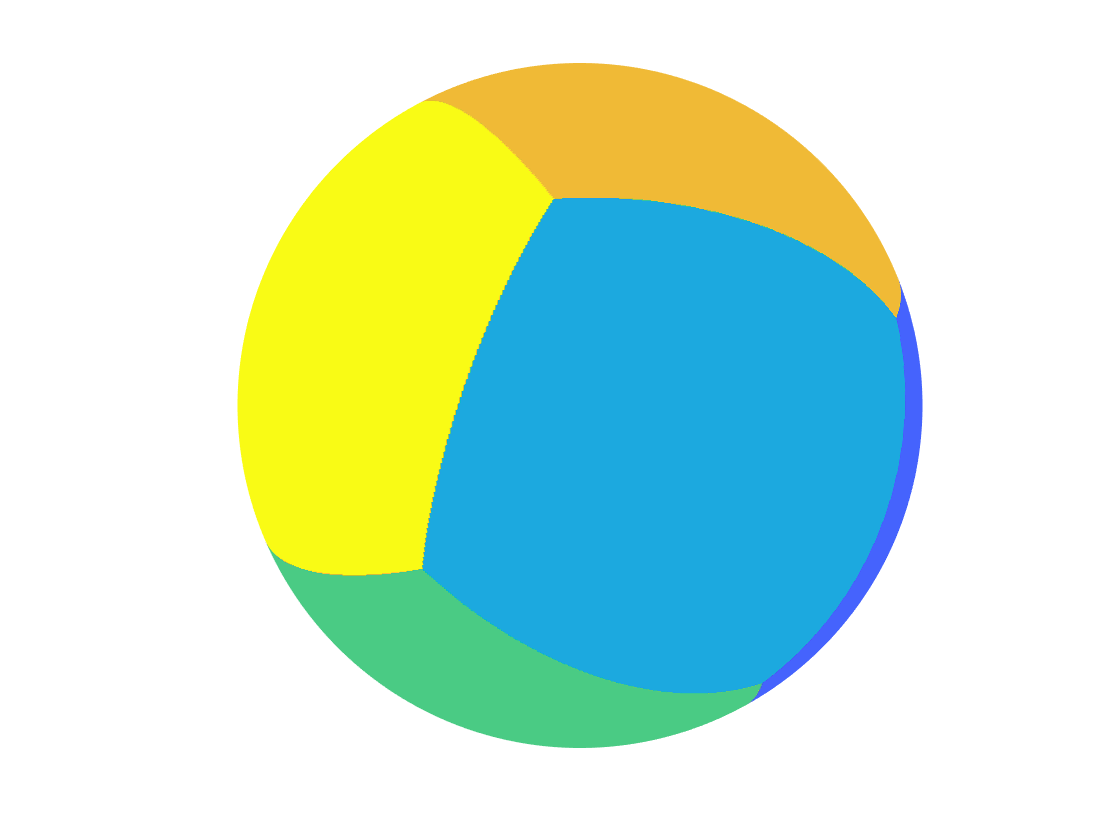}

\includegraphics[scale=0.15,clip,trim= 12cm 5cm 10cm 5cm]{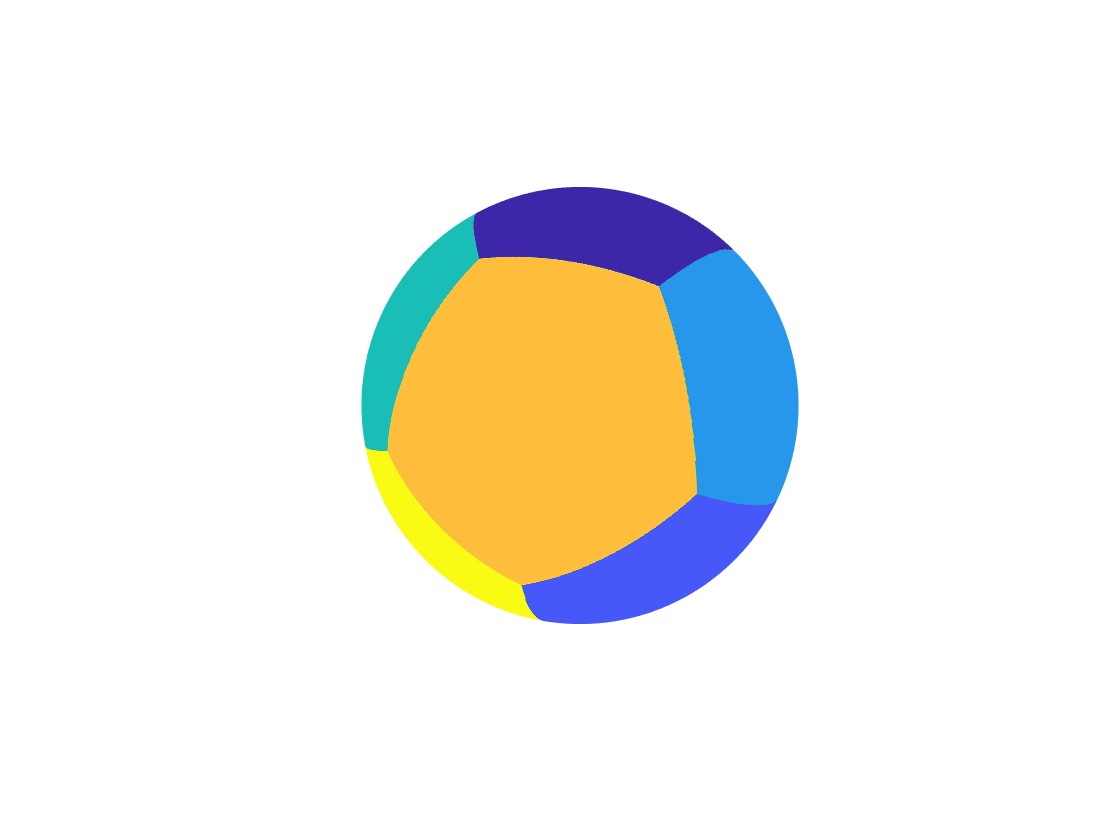} \  \  \ \ 
\includegraphics[scale=0.12,clip,trim= 8cm 2cm 7cm 2cm]{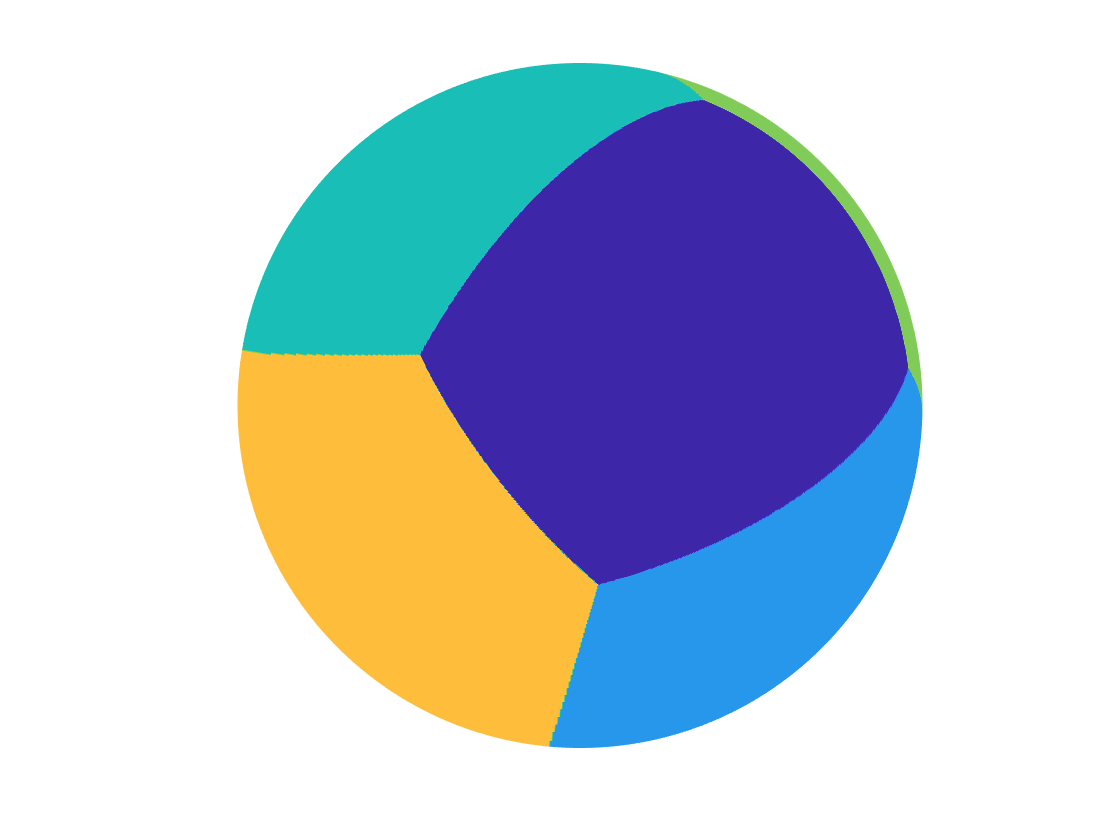}
\includegraphics[scale=0.12,clip,trim= 8cm 2cm 7cm 2cm]{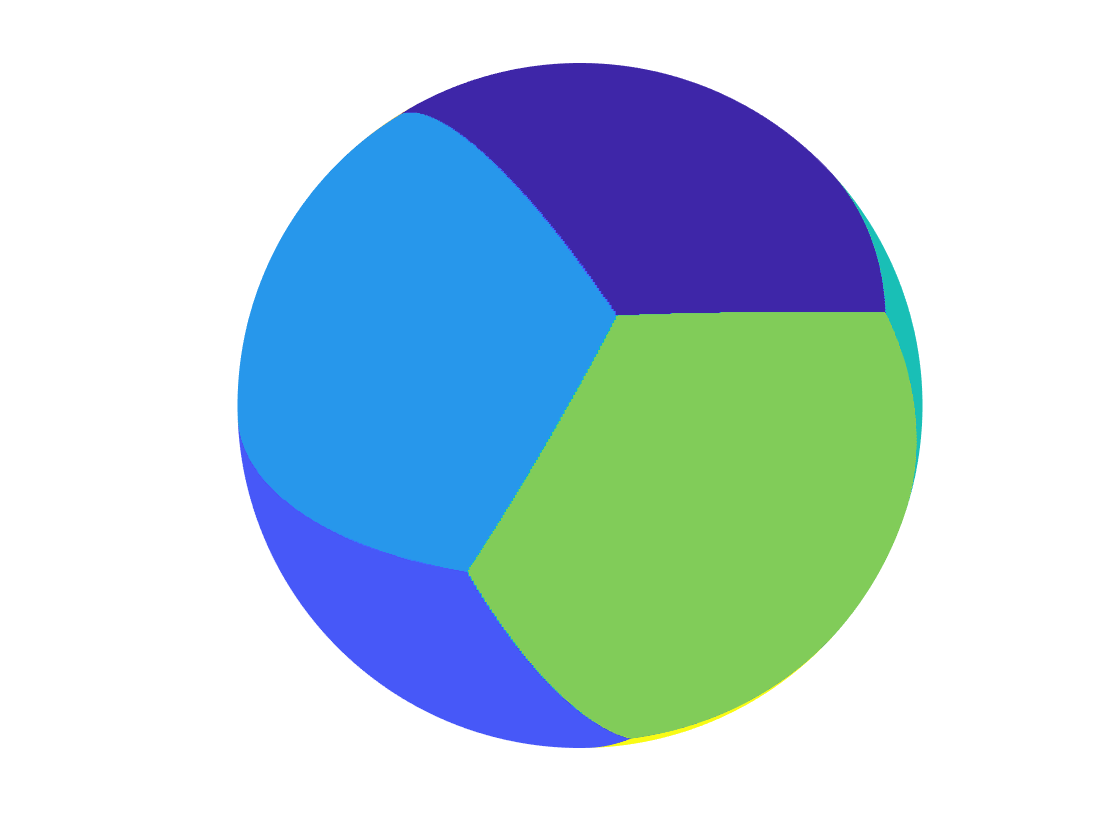}
\includegraphics[scale=0.12,clip,trim= 8cm 2cm 7cm 2cm]{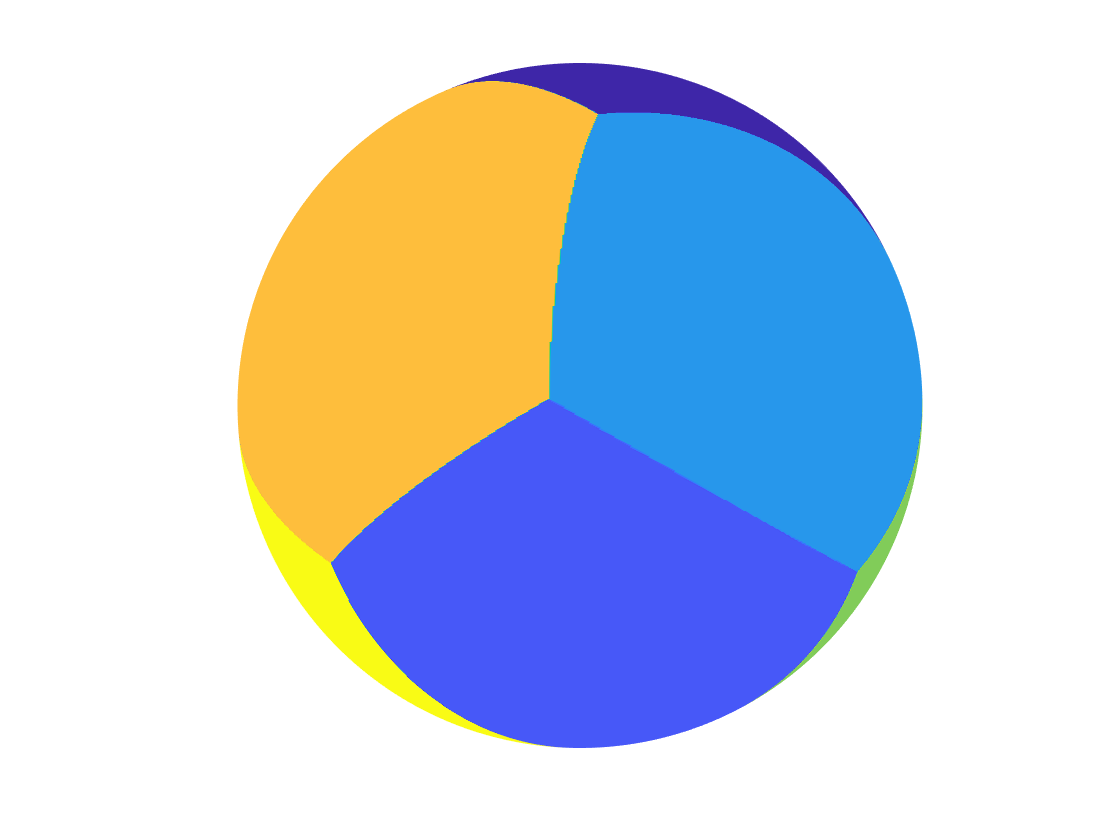}
\caption{Column 1: $k$-Dirichlet partitions of a sphere. Column 2: Vertical view. Column 3: Front view.  Column 4: Side view. From top to bottom:  $k$-Dirichlet partitions of a sphere with $k = 3$--7. The CPU time for each case was $180$, $485$, $727$, $901$, and $1231$ seconds respectively.} 
\label{fig:sphere1}
\end{figure}

\begin{figure}
\includegraphics[scale=0.15,clip,trim= 12cm 5cm 10cm 5cm]{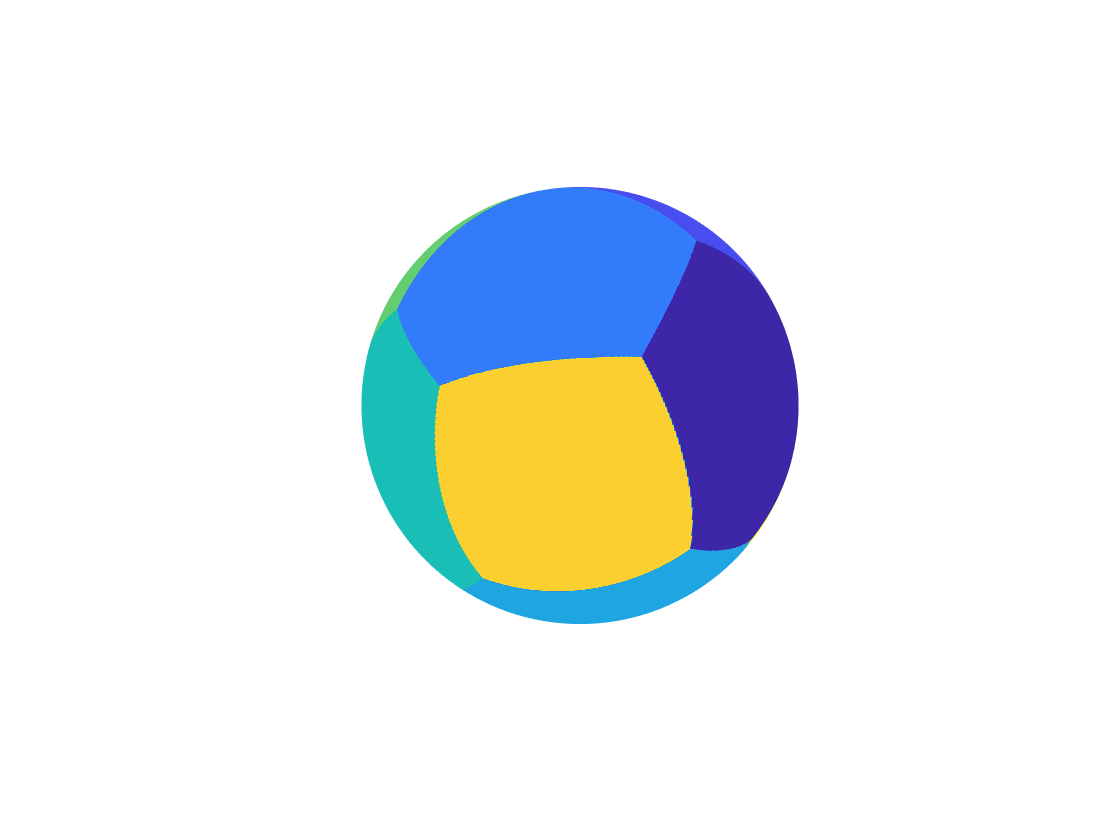} \  \  \ \ 
\includegraphics[scale=0.12,clip,trim= 8cm 2cm 7cm 2cm]{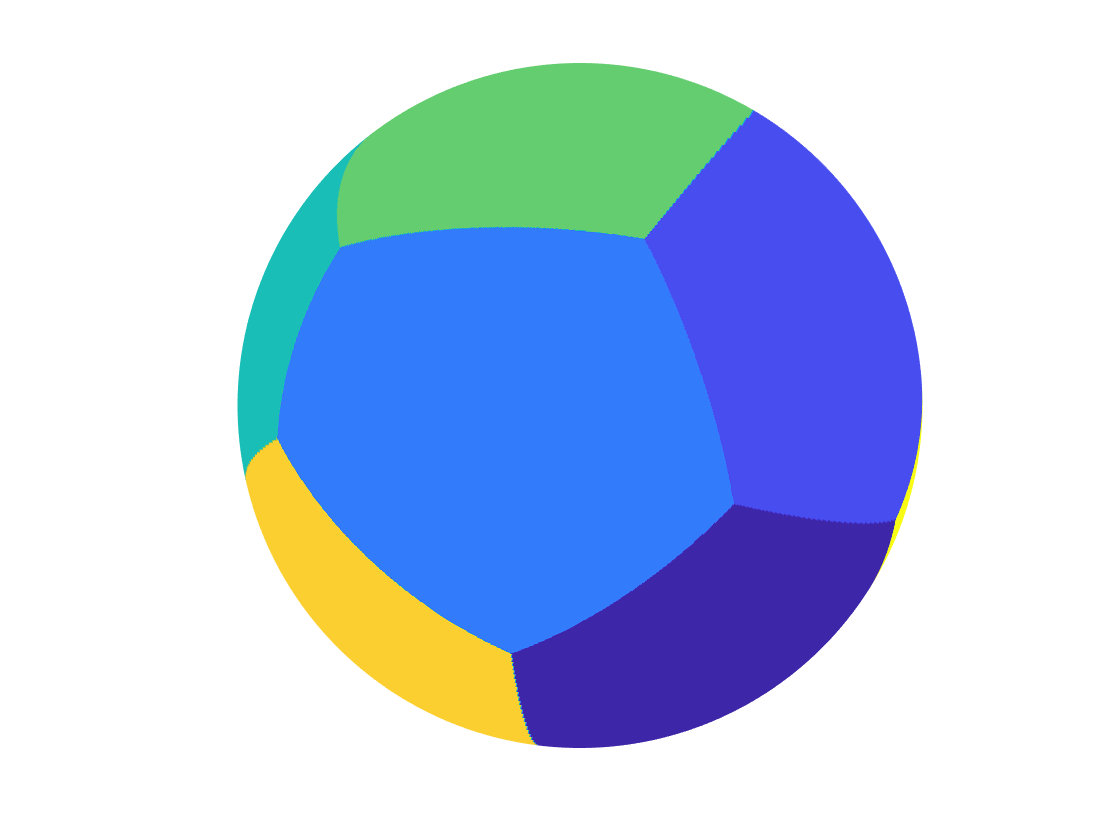}
\includegraphics[scale=0.12,clip,trim= 8cm 2cm 7cm 2cm]{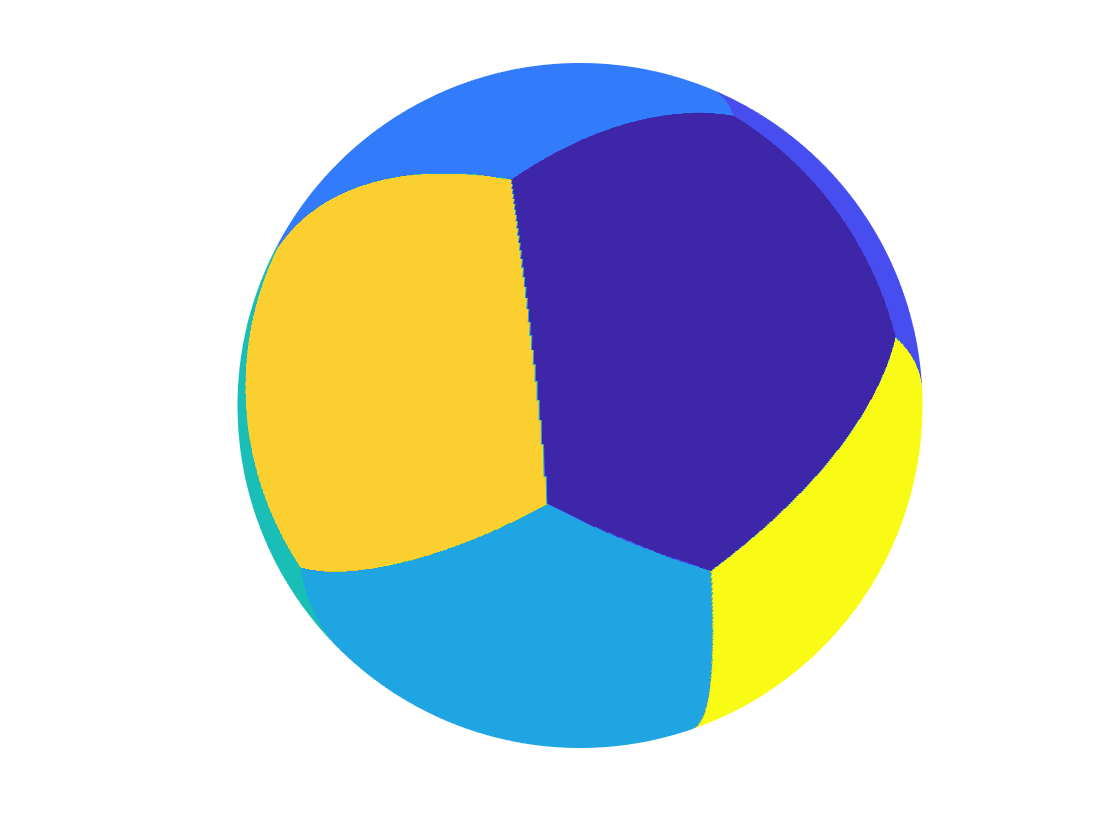}
\includegraphics[scale=0.12,clip,trim= 8cm 2cm 7cm 2cm]{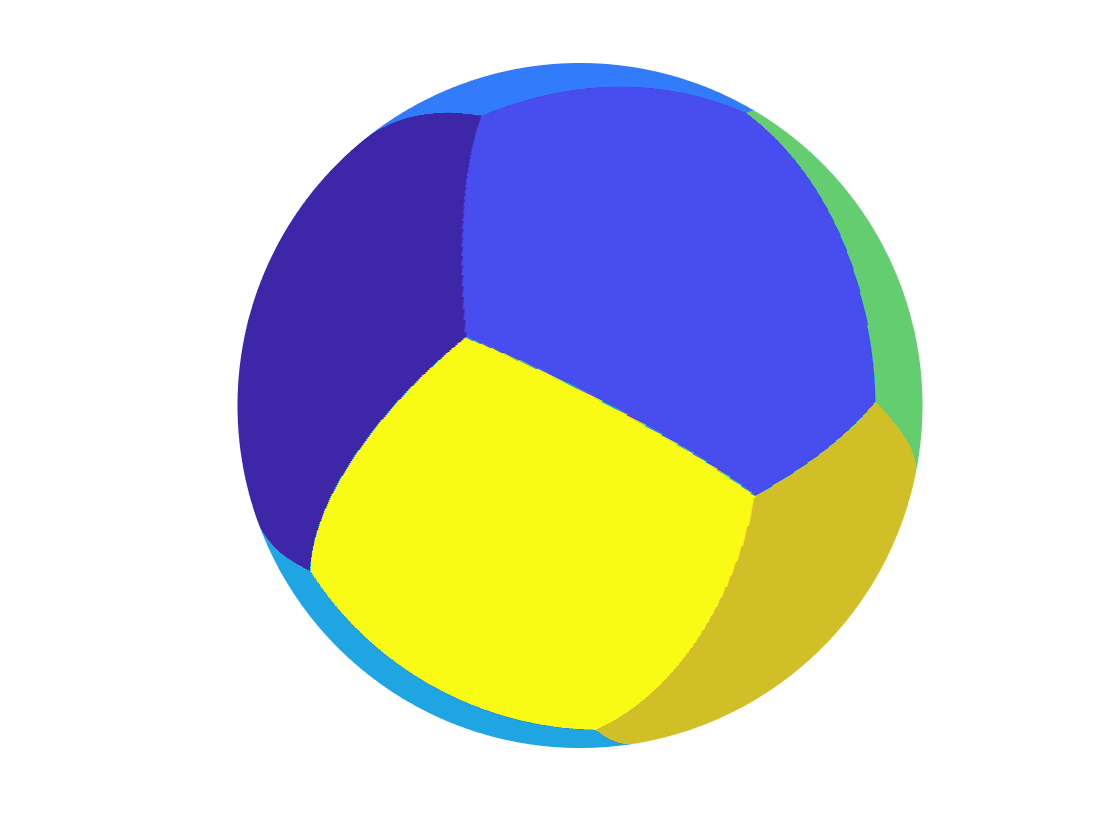}

\includegraphics[scale=0.15,clip,trim= 12cm 5cm 10cm 5cm]{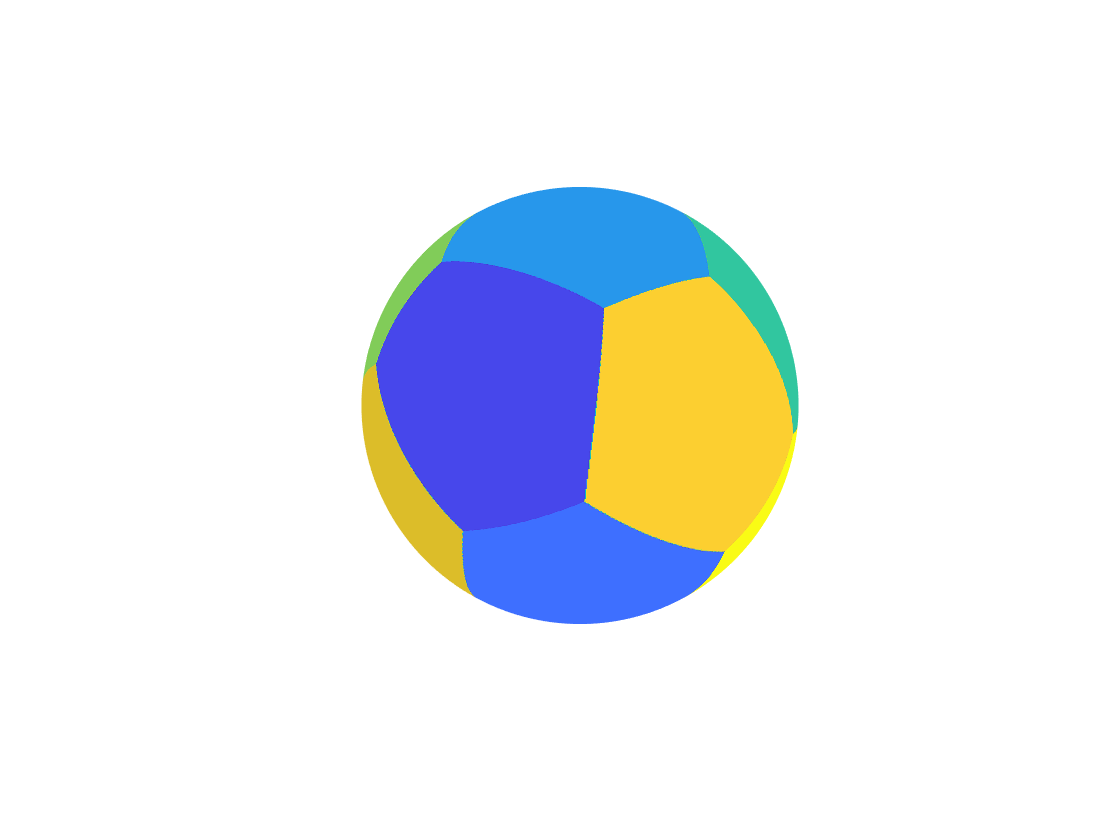} \  \  \ \ 
\includegraphics[scale=0.12,clip,trim= 8cm 2cm 7cm 2cm]{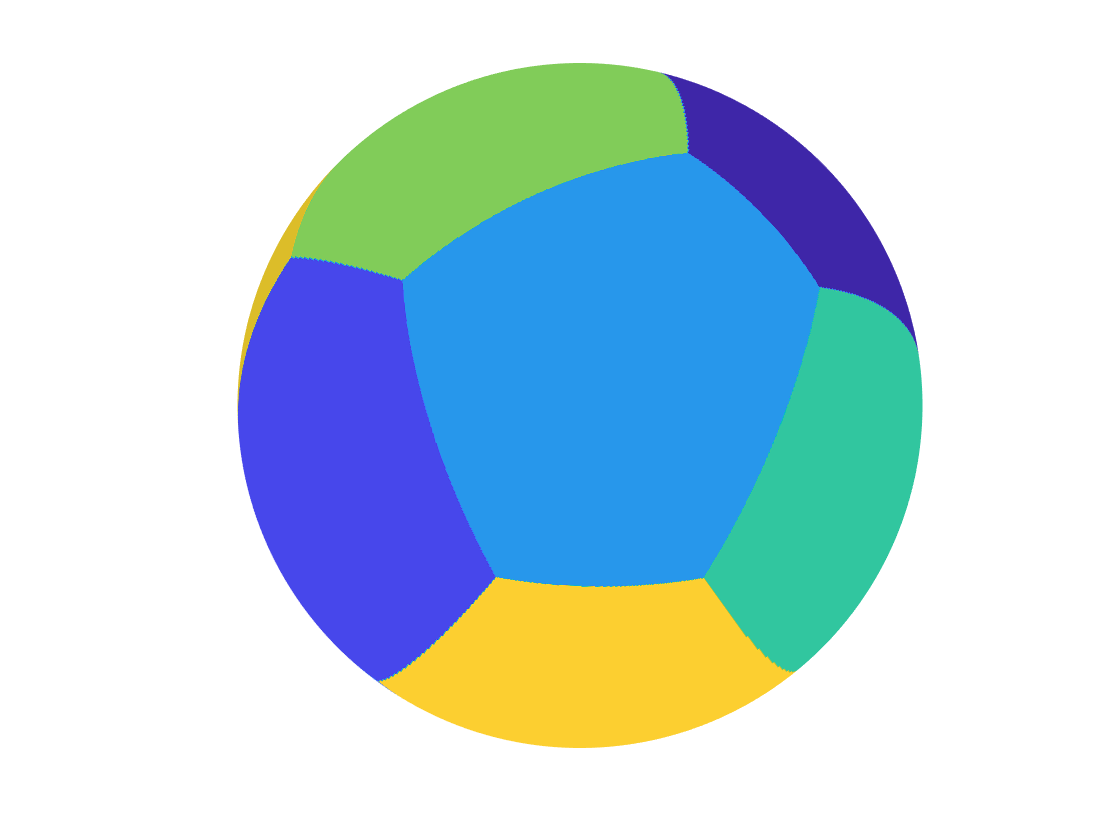}
\includegraphics[scale=0.12,clip,trim= 8cm 2cm 7cm 2cm]{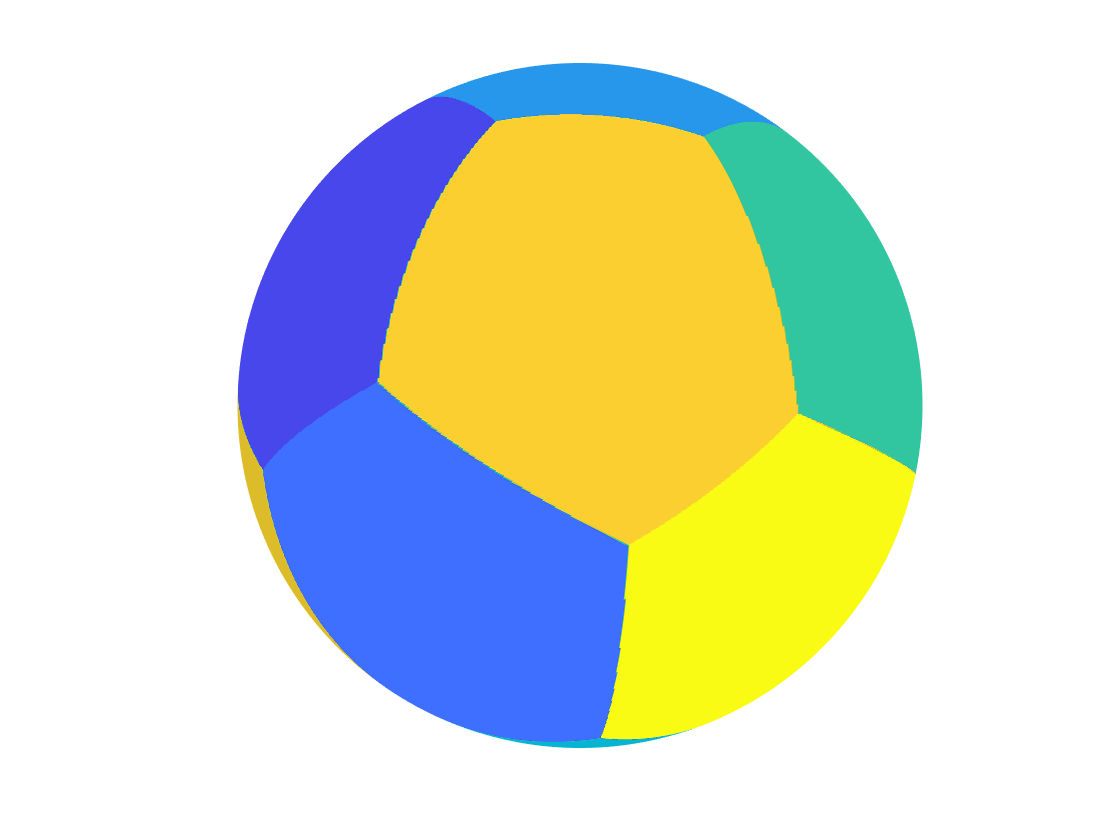}
\includegraphics[scale=0.12,clip,trim= 8cm 2cm 7cm 2cm]{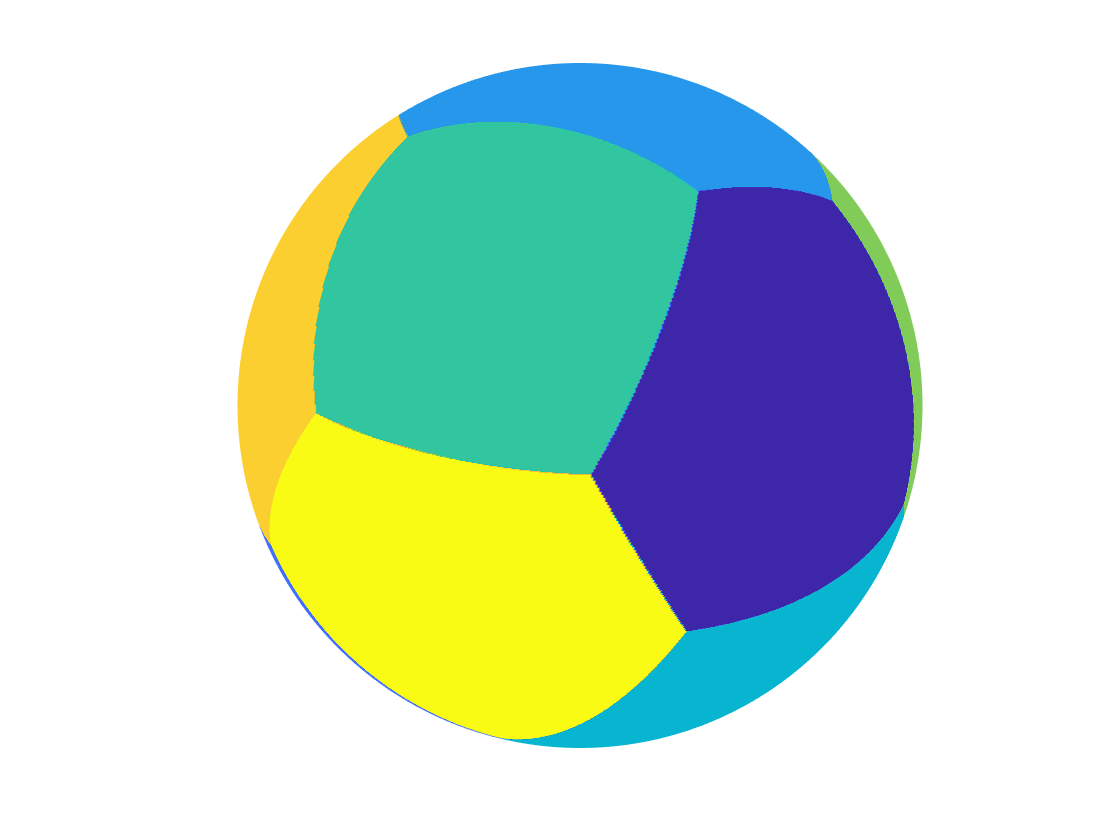}

\includegraphics[scale=0.15,clip,trim= 12cm 5cm 10cm 5cm]{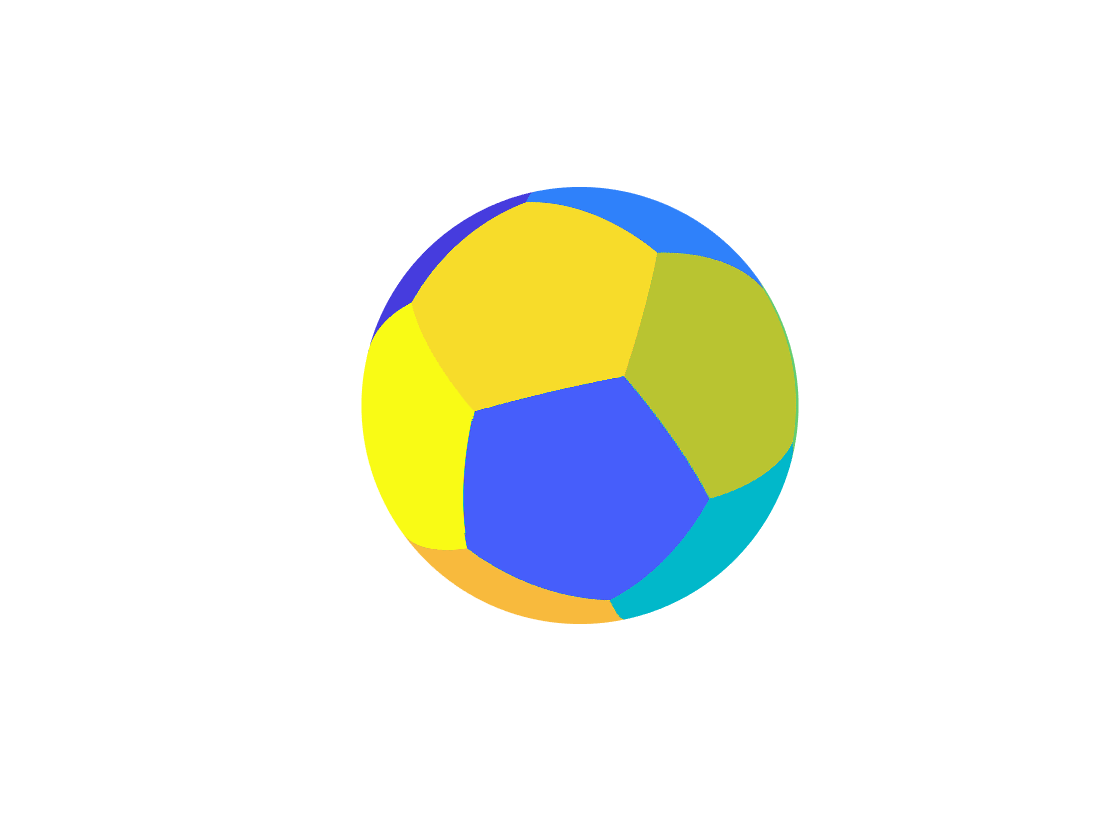} \  \  \ \ 
\includegraphics[scale=0.12,clip,trim= 8cm 2cm 7cm 2cm]{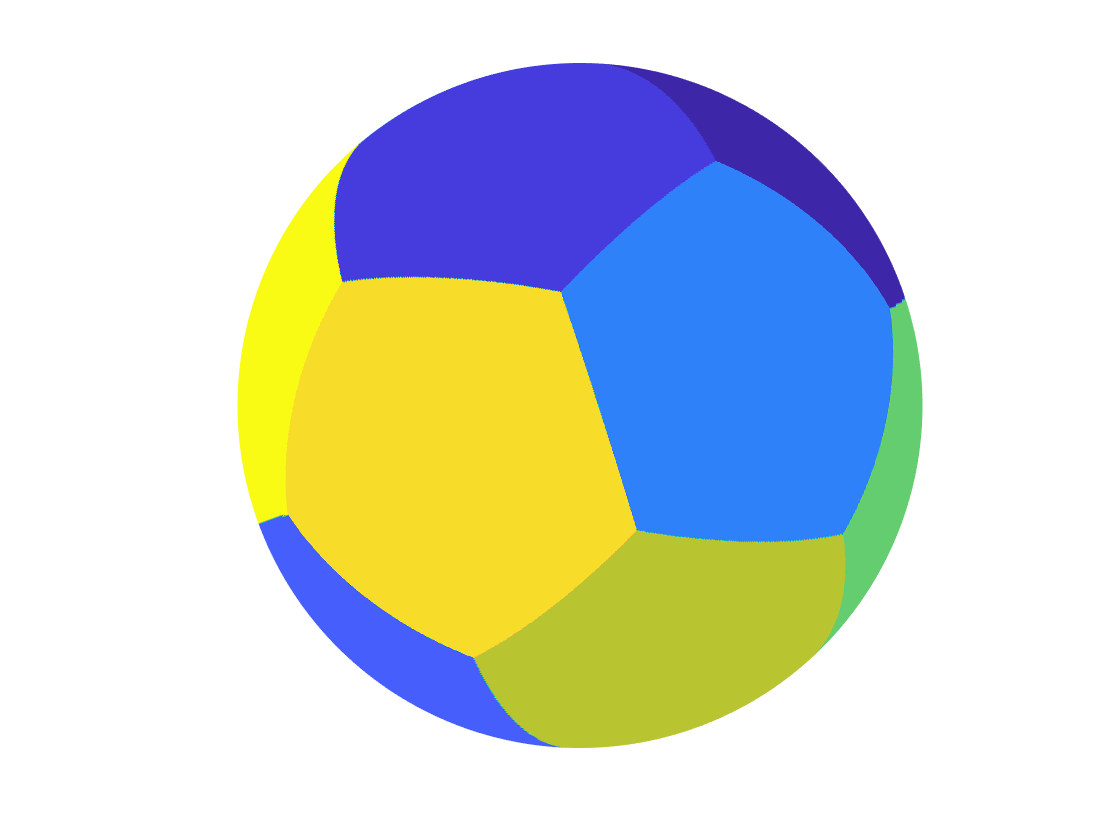}
\includegraphics[scale=0.12,clip,trim= 8cm 2cm 7cm 2cm]{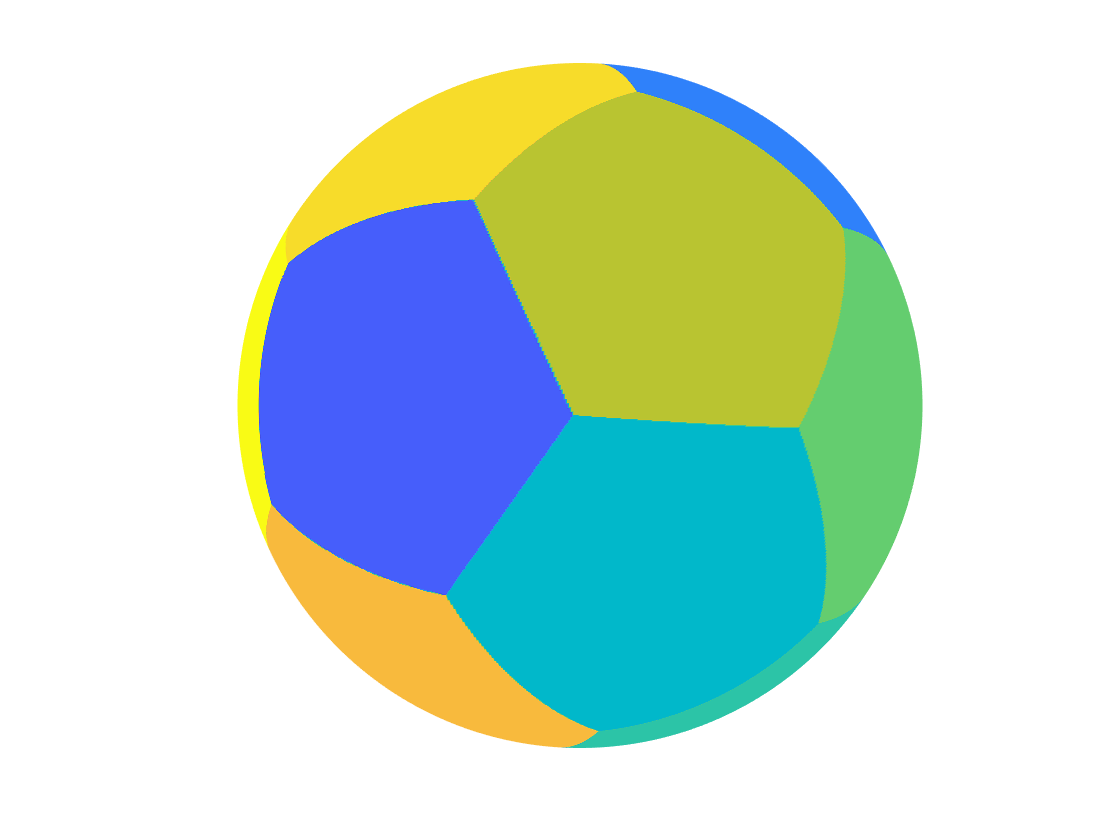}
\includegraphics[scale=0.12,clip,trim= 8cm 2cm 7cm 2cm]{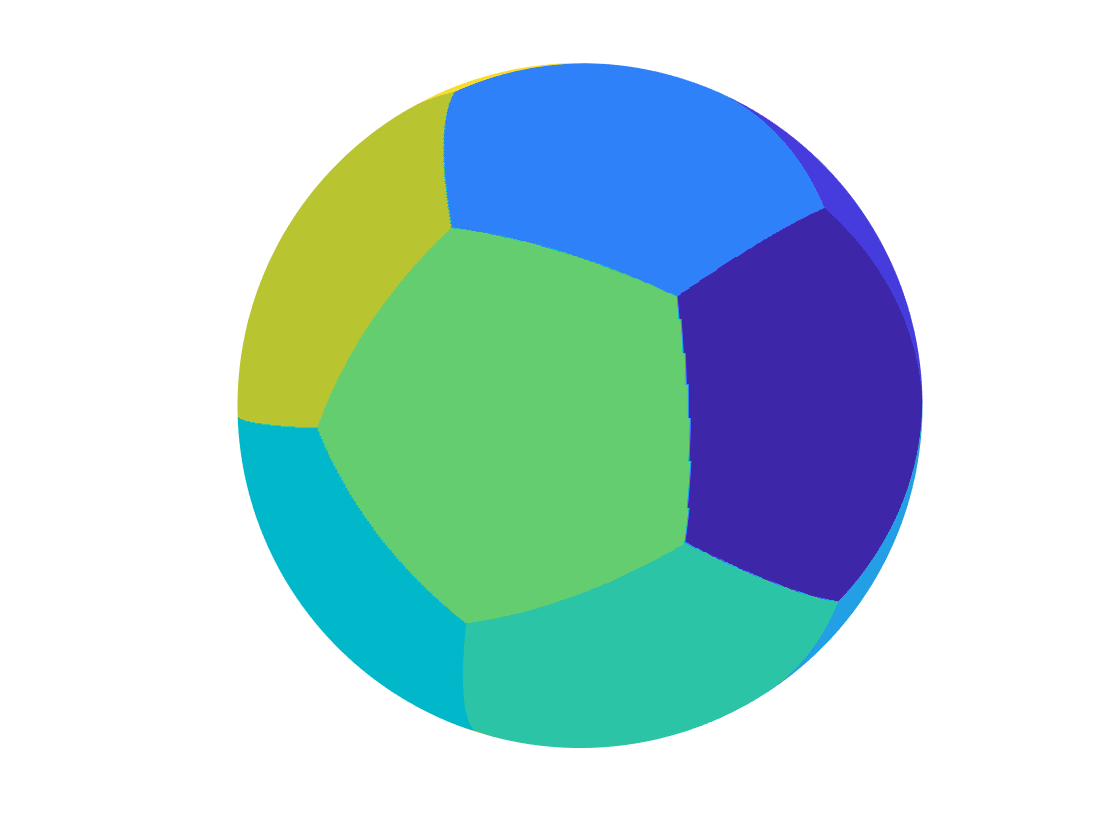}

\includegraphics[scale=0.15,clip,trim= 12cm 5cm 10cm 5cm]{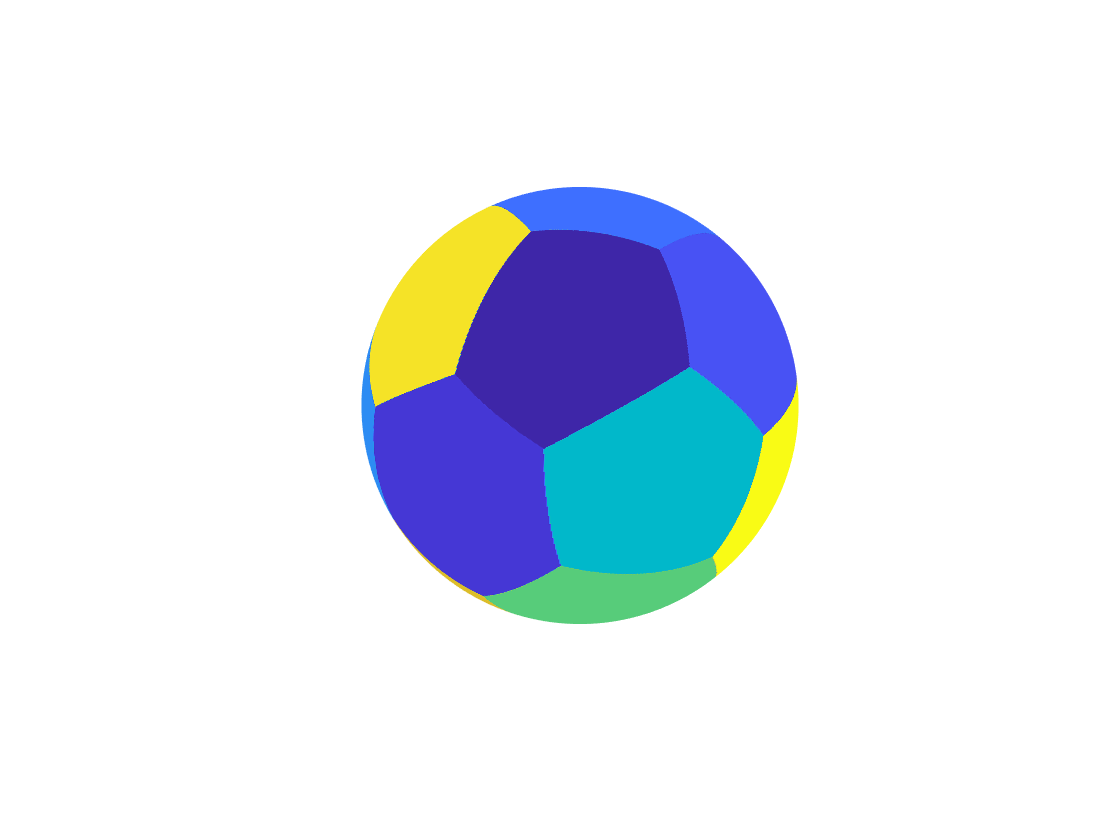} \  \  \ \ 
\includegraphics[scale=0.12,clip,trim= 8cm 2cm 7cm 2cm]{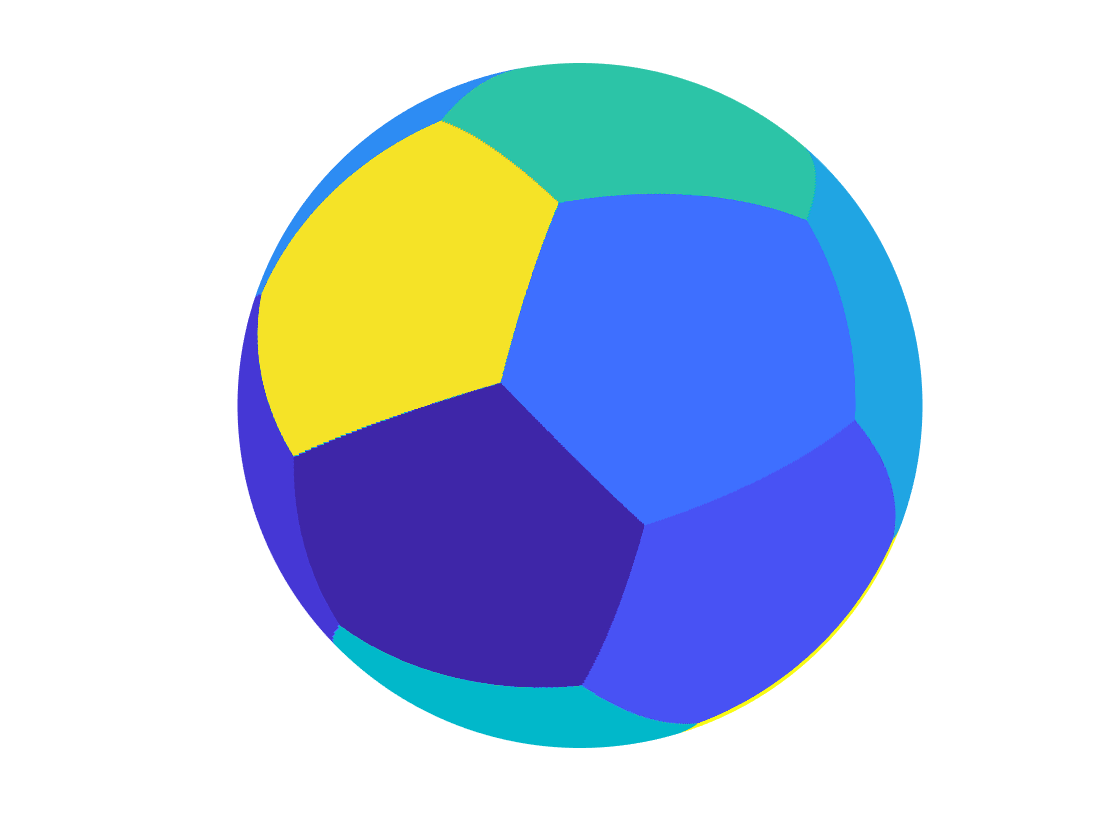}
\includegraphics[scale=0.12,clip,trim= 8cm 2cm 7cm 2cm]{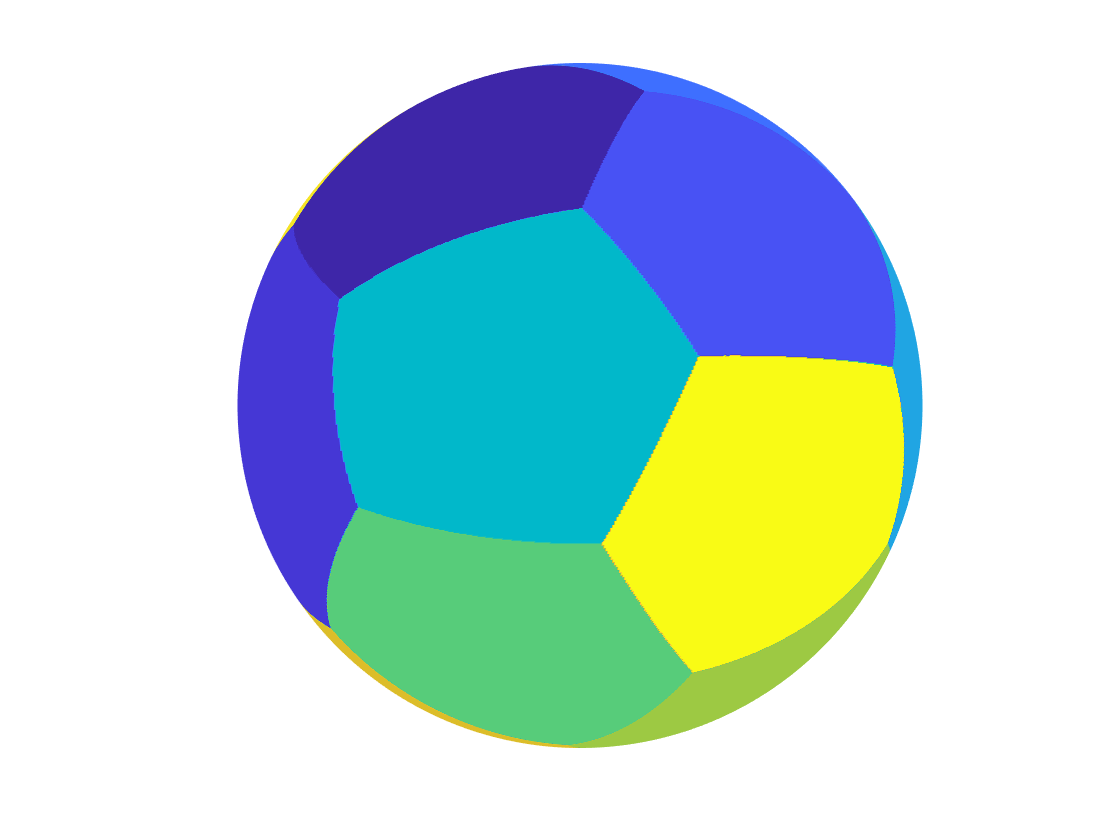}
\includegraphics[scale=0.12,clip,trim= 8cm 2cm 7cm 2cm]{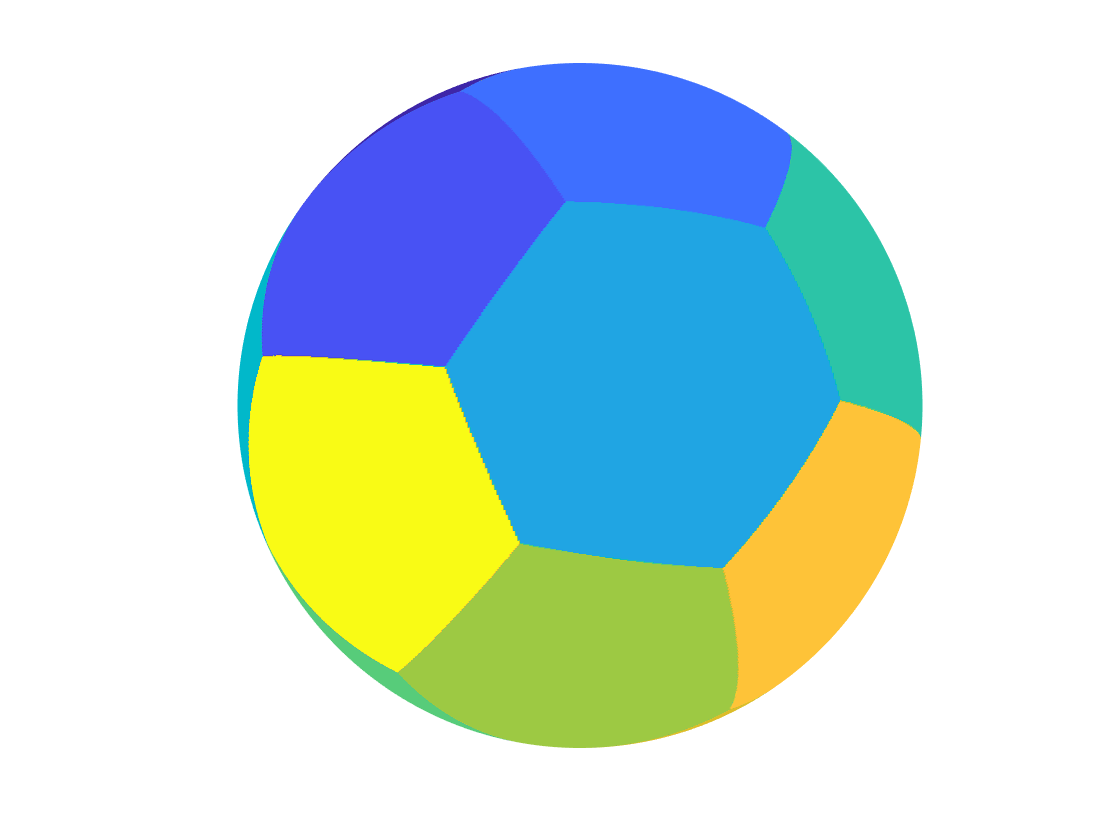}

\includegraphics[scale=0.15,clip,trim= 12cm 5cm 10cm 5cm]{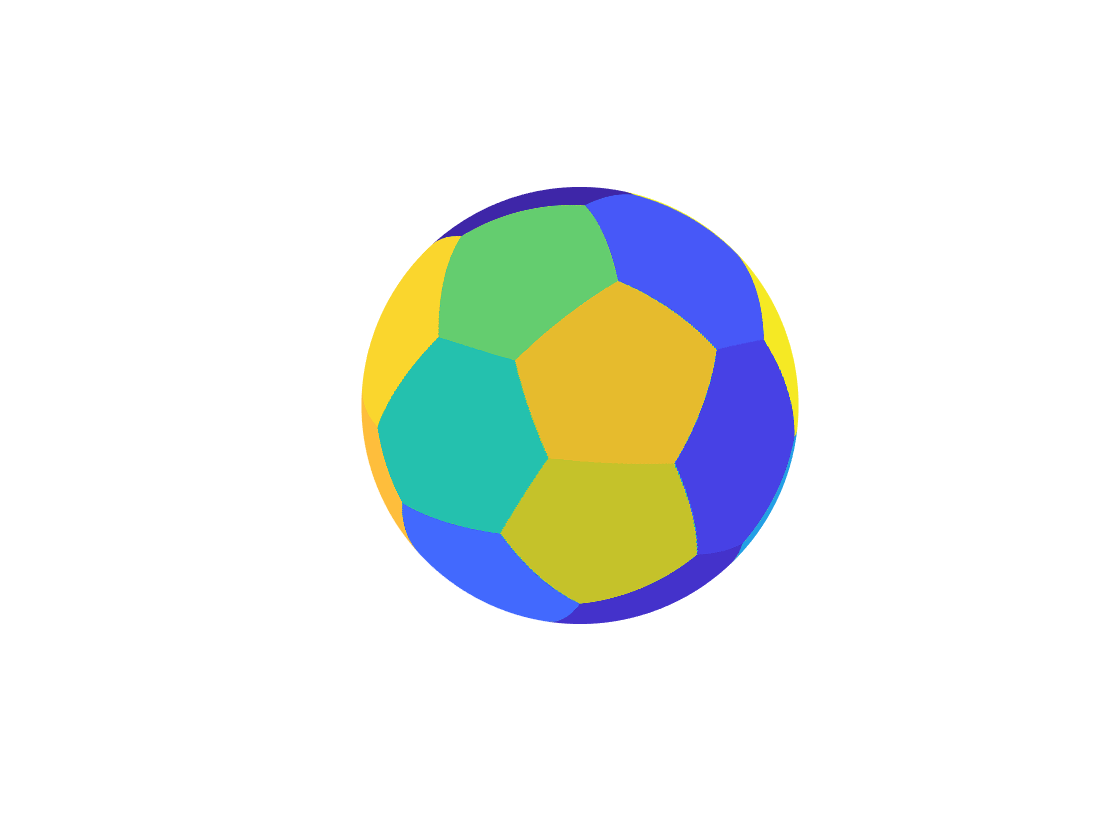} \  \  \ \ 
\includegraphics[scale=0.12,clip,trim= 8cm 2cm 7cm 2cm]{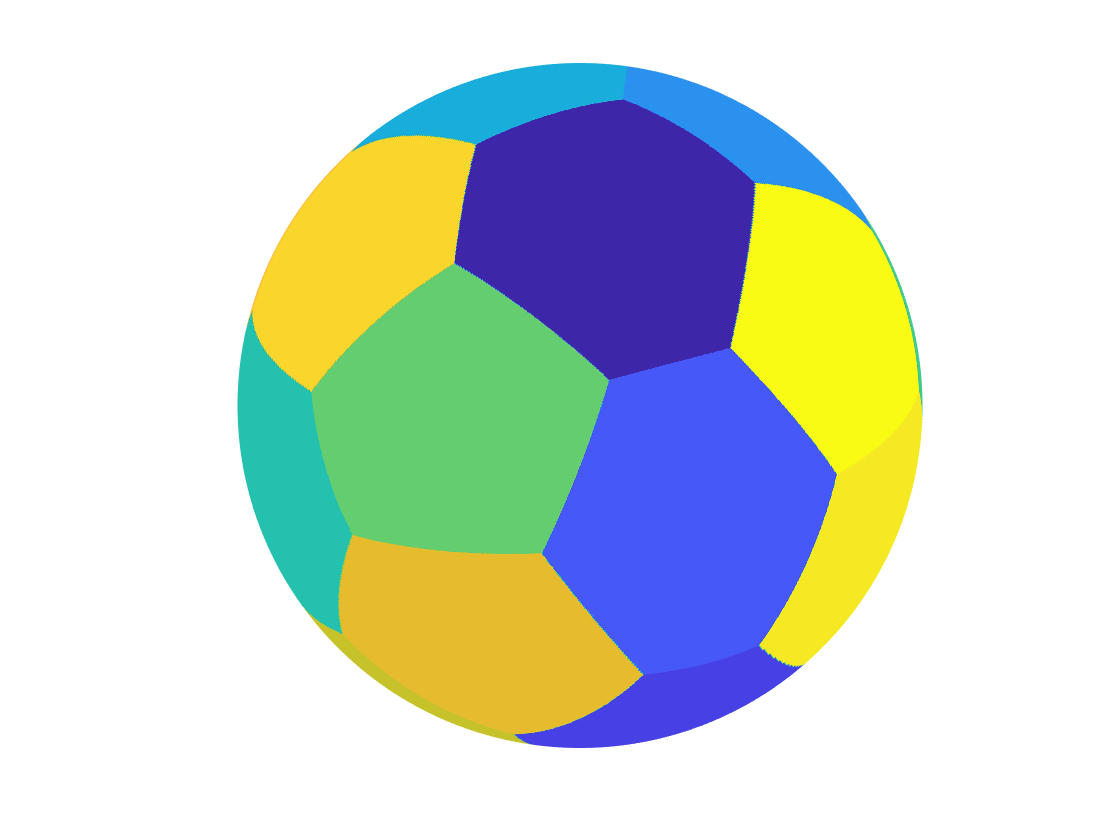}
\includegraphics[scale=0.12,clip,trim= 8cm 2cm 7cm 2cm]{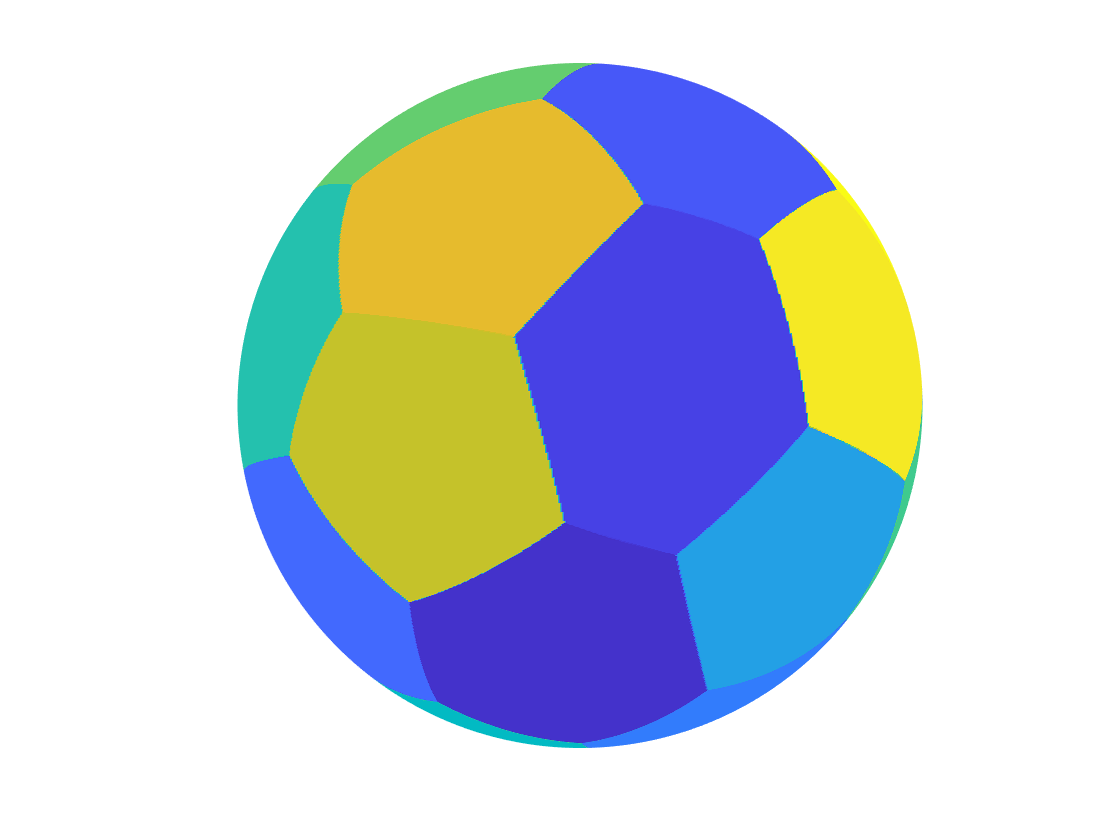}
\includegraphics[scale=0.12,clip,trim= 8cm 2cm 7cm 2cm]{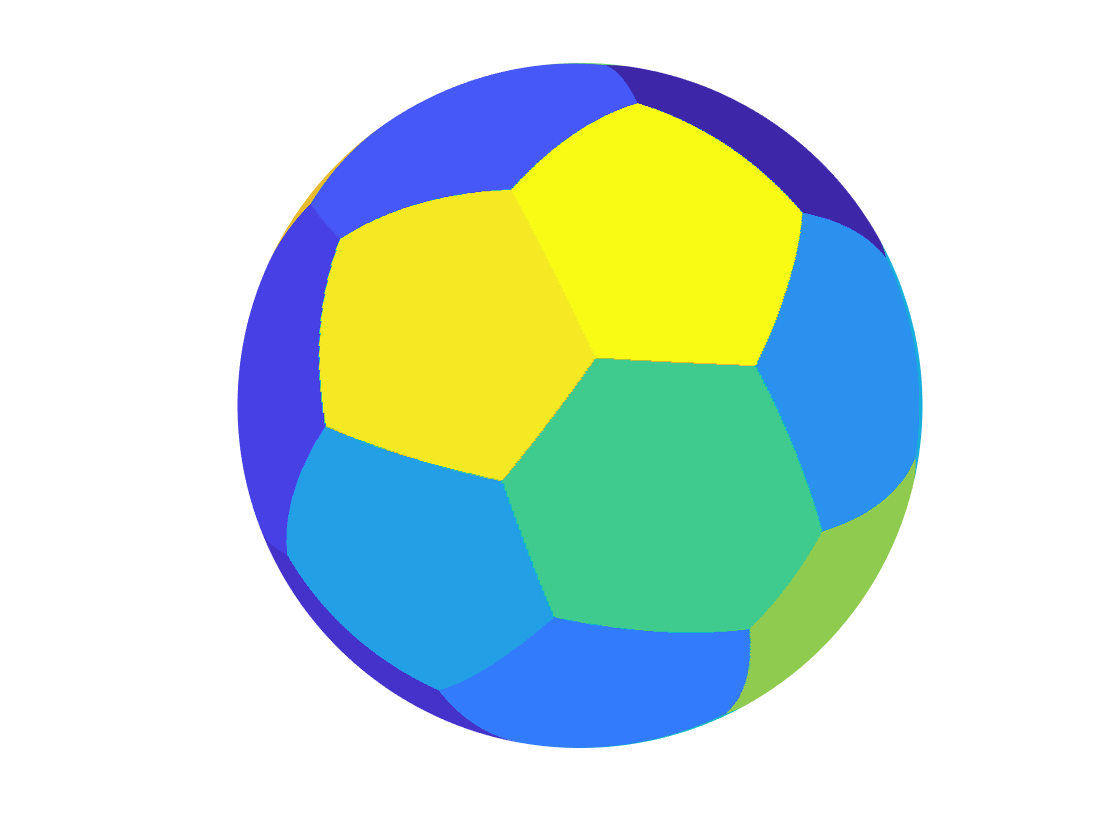}

\caption{Column 1: $k$-Dirichlet partitions of a sphere. Column 2: Vertical view. Column 3: Front view.  Column 4: Side view. 
From top to bottom: $k$-Dirichlet partitions of a sphere  with $k = 9$,10,12,14, and $20$. 
The CPU time for each case was $ 2040$, $2165$, $1631$, $1769$, and $9011$ seconds respectively.} 
\label{fig:sphere2}

\end{figure}

In Figures~\ref{fig:sphere1} and \ref{fig:sphere2}, we display Dirichlet partitions on the sphere  for $k=3$--7,9,10,12,14, and 20, obtained from an initialization using a random tessellation. 
In Table~\ref{tab:sphere}, the CPU times for each case are given. 
For parameterization as in \eqref{e:SphereCoord}, the inclination and azimuthal coordinates are discretized by $256^2$ uniform grid points and $\tau = 0.008$. 
Values of  $\tilde E$ in \eqref{e:aEnergy} for different values of $k$ are displayed in Table~\ref{tab:sphere}.

\begin{table}[ht]
\centering
\caption{Values of  $\tilde E$ in \eqref{e:aEnergy} and the average CPU time for different values of $k$.} \label{tab:sphere}
\begin{tabular}{|c|c|c|c|c|c|c|c|c|c|c|c|c|}
\hline
$k$ &  3 & 4 & 5 &6 &7\\
\hline 
$\tilde E$  & 13.49 & 13.64 & 14.16 & 13.73 & 13.96  \\
 \hline
 CPU time (s) &  180 & 485 & 727 & 901 & 1231  \\
 \hline
 \hline
$k$ &  9 & 10 & 12 & 14 &20\\
\hline 
$\tilde E$  & 13.65 & 13.54 & 13.08 & 12.95 & 12.20  \\
 \hline
 CPU time (s)  & 2040 & 2165 & 1631 & 1769 & 9011  \\
 \hline

\end{tabular}
\end{table}

{\small
\bibliography{SharpInterfaceDPart}
}

\end{document}